\documentclass[12pt]{article}
\usepackage{latexsym}
\usepackage{amsfonts}
\usepackage{amssymb}
\parskip=\medskipamount                 
\thicklines                     
\newcommand{\bimn}[7]{\bibitem{#1}#2,
{\em #3},
{ #4}\hspace{0.25em}{\bf
#5}\hspace{0.25em}(#6)\hspace{0.25em}{#7}.}

\def\inbar{\vrule height1.5ex width.4pt depth0pt}
\def\IC{\relax\,\hbox{$\inbar\kern-.3em{\rm C}$}}
\def\IN{\relax{\rm I\kern-.18em N}}
\def\IQ{\relax\,\hbox{$\inbar\kern-.3em{\rm Q}$}}
\def\IR{\relax{\rm I\kern-.18em R}}
\def\ZZ{\relax{\sf Z\kern-.4em Z}}
\def\a{\alpha} \def\b{\beta}   \def\e{\epsilon} \def\g{\gamma}
 \def\l{\lambda}

\def\cF{{\cal F}}   
 \def\cK{{\cal K}} \def\cL{{\cal L}}

\newtheorem{theorem}{Theorem}[section]
\newtheorem{proposition}[theorem]{Proposition}
\newtheorem{corollary}[theorem]{Corollary}

\newtheorem{lemma}[theorem]{Lemma}
\newtheorem{definition}[theorem]{Definition}
\newtheorem{remark}[theorem]{Remark}

\catcode`\@=11

\newif\if@fewtab\@fewtabtrue

\catcode`\@=11

\newif\if@fewtab\@fewtabtrue

{\count255=\time\divide\count255 by 60
\xdef\hourmin{\number\count255}
\multiply\count255 by-60\advance\count255 by\time
\xdef\hourmin{\hourmin:\ifnum\count255<10 0\fi\the\count255}}
\def\ps@draft{\let\@mkboth\@gobbletwo
    \def\@oddhead{}
    \def\@oddfoot
      {\hbox to 7 cm{\footnotesize {\em Draft of \jobname:} \draftdate
       \hfil}\hskip -7cm\hfil\rm\thepage \hfil}
    \def\@evenhead{}\let\@evenfoot\@oddfoot}

\def\ceqno{\global\@fewtabfalse
    \ifcase\@eqcnt \def\@tempa{& & &}\or \def\@tempa{& &}
      \or \def\@tempa{&}
      \or\def\@tempa{}\fi\@tempa
{\rm(\theequation)}}

\def\aeqno#1{\global\@fewtabfalse
    \ifcase\@eqcnt \def\@tempa{& & &}\or \def\@tempa{& &}
      \or \def\@tempa{&}
      \or\def\@tempa{}\fi\@tempa
{\rm(\theequation,#1)}}

\def\label#1{\ifnum\draftcontrol=1
 \global\def\draftnote{$\scriptstyle #1$}\fi
 \@bsphack\if@filesw {\let\thepage\relax
   \def\protect{\noexpand\noexpand\noexpand}%
\xdef\@gtempa{\write\@auxout{\string
      \newlabel{#1}{{\@currentlabel}{\thepage}}}}}\@gtempa
   \if@nobreak \ifvmode\nobreak\fi\fi\fi
  \@esphack}

\def\alabel#1#2{\label{#1}\global\@fewtabfalse
    \ifcase\@eqcnt \def\@tempa{& & &}\or \def\@tempa{& &}
      \or \def\@tempa{&}
      \or\def\@tempa{}\fi\@tempa
{\hbox to 3cm{\phantom{\rm(\theequation,#2)}
\draftnote \hfil}\hskip -3cm {\rm(\theequation,#2)}}}

\def\clabel#1{\label{#1}\global\@fewtabfalse
    \ifcase\@eqcnt \def\@tempa{& & &}\or \def\@tempa{& &}
      \or \def\@tempa{&}
      \or\def\@tempa{}\fi\@tempa
{\hbox to 3cm{\phantom{\rm(\theequation)}
\draftnote \hfil}\hskip -3cm{\rm(\theequation)}}}

\def\eqnarray{\def\draftnote{{}}\global\@fewtabtrue
\stepcounter{equation}\let\@currentlabel=\theequation
\global\@eqnswtrue
\global\@eqcnt\z@\tabskip\@centering\let\\=\@eqncr
$$\halign to \displaywidth\bgroup\@eqnsel\hskip\@centering\@eqcnt\z@
  $\displaystyle\tabskip\z@{##}$&\global\@eqcnt\@ne
  \hskip 1\arraycolsep \hfil$\displaystyle{##}$\hfil
  &\global\@eqcnt\tw@ \hskip 1\arraycolsep
$\displaystyle\tabskip\z@{##}$
\hfil  \tabskip\@centering&\global\@eqcnt\thr@@\llap{##}\tabskip\z@
\cr}

\def\endeqnarray{\@@eqncr\egroup
      \global\advance\c@equation\m@ne$$\global\@ignoretrue}

\def\@eqnnum{\hbox to 3cm{\phantom{\rm(\theequation)} \draftnote
                         \hfil}\hskip -3cm {\rm(\theequation)}}

\def\@@eqncr{\let\@tempa\relax
    \ifcase\@eqcnt \def\@tempa{& & &}\or \def\@tempa{& &}
      \or \def\@tempa{&}
      \or\def\@tempa{}
\fi\@tempa
\if@eqnsw
\if@fewtab\@eqnnum\fi
\stepcounter{equation}\fi\global
\@eqnswtrue\global\@eqcnt\z@\global\@fewtabtrue\cr}

\def\draftcite#1{\ifnum\draftcontrol=1#1\else{}\fi}

\def\@lbibitem[#1]#2{\item{}\hskip -3cm \hbox to 2cm
{\hfil$\scriptstyle\draftcite{#2}$}\hskip
1cm[\@biblabel{#1}]\if@filesw
     {\def\protect##1{\string ##1\space}\immediate
      \write\@auxout{\string\bibcite{#2}{#1}}}\fi\ignorespaces}

\def\@bibitem#1{\item\hskip -3cm \hbox to 2cm
{\hfil $\scriptstyle\draftcite{#1}$}\hskip 1cm
\if@filesw \immediate\write\@auxout
       {\string\bibcite{#1}{\the\value{\@listctr}}}\fi\ignorespaces}

\def\nsection#1{\section{#1}\setcounter{equation}{0}}

\def\nappendixsp#1#2{\def\thesection{A#1}\section*{Appendix #2}
\def\theequation{{A#1.\arabic{equation}}}
\def\thetheorem{{A#1.\arabic{theorem}}}
\setcounter{equation}{0}
\setcounter{theorem}{0}}

\def\draftdate{\number\month/\number\day/\number\year\ \ \ \hourmin }

\global\def\draftcontrol{0}
\catcode`\@=12

\def\theequation{{\thesection.\arabic{equation}}}

\def\qq{\begin{eqnarray}}
\def\qqq{\end{eqnarray}}
\def\rx#1{~(\ref{#1})}
\def\ex#1{eq.\hspace*{-3pt}\rx{#1}}
\def\eex#1{eqs.\hspace*{-3pt}\rx{#1}}
\def\cx#1{~\cite{#1}}
\def\rw#1{~\ref{#1}}

\newlength{\shiftwidth}
\def\shift#1{&&\hbox to \shiftwidth{\hfill $\displaystyle#1$}}
\newlength{\sshiftwidth}
\def\sshift#1{\lefteqn{\hbox to
\sshiftwidth{\hfill$\displaystyle#1$}}}

\def\cf{{\it cf.\ }}
\def\rhs{{\it r.h.s.\ }}
\def\lhs{{\it l.h.s.\ }}

\def\Rhs{$\IQ$HS\ }
\def\Rhsp{$\IQ$HS}

\def\ordH{{|H_1(M,\ZZ)|} }
\def\deg{ \mathop{{\rm deg}}\nolimits }
\def\p{^{\prime}}

\def\prosign{\mathop{{\rm sign}}\nolimits}

\def\intinf{ \int_{-\infty}^{+\infty} }

\def\lk{\mathop{{\rm lk}}\nolimits}

\def\pr#1#2{ \noindent{\em Proof of #1~\ref{#2}.} }
\def\proof{ \noindent{\em Proof.} }

\def\qed{ \hfill $\Box$ }
\def\const{ {\mbox{const}} }

\def\lrbc#1{ \left( #1 \right) }
\def\lrbs#1{ \left[ #1 \right] }
\def\lrbf#1{ \left\{ #1 \right\} }

\def\u#1{ \underline{#1} }
\def\uu#1{ \underline{\underline{#1}} }

\def\ux{ {\u{x}} }
\def\uy{ {\u{y}} }
\def\ua{ {\u{a}} }
\def\ual{ {\u{\a}} }
\def\ug{ {\u{g}} }
\def\um{ {\u{m}} }
\def\ut{ {\u{t}} }

\def\uux{ {\uu{x}} }

\def\sumn{ \sum_{0\leq |\um| \leq n} }
\def\smnt{ \sum_{0\leq m \leq n/2} }
\def\sjoL{ \sum_{j=1}^L }
\def\sjoN{ \sum_{j=1}^N }

\def\pjoL{ \prod_{j=1}^L }

\def\pjoN{ \prod_{j=1}^N }

\def\lk#1#2{\mathop{{\rm lk}(#1;#2)} \nolimits }

\def\Zh{ \ZZ\lrbs{\hlf} }

\def\Zh{ \ZZ[1/2] }

\def\ZZ{ \mathbb{Z} }
\def\IQ{ \mathbb{Q} }
\def\IC{ \mathbb{C} }

\def\qbezier{\bezier{120}}
\def\DottedCircle{
\bezier{4}(0.966,-0.259)(1.04,0)(0.966,0.259)
\bezier{4}(0.966,0.259)(0.897,0.518)(0.707,0.707)
\bezier{4}(0.707,0.707)(0.518,0.897)(0.259,0.966)
\bezier{4}(0.259,0.966)(0,1.04)(-0.259,0.966)
\bezier{4}(-0.259,0.966)(-0.518,0.897)(-0.707,0.707)
\bezier{4}(-0.707,0.707)(-0.897,0.518)(-0.966,0.259)
\bezier{4}(-0.966,0.259)(-1.04,0)(-0.966,-0.259)
\bezier{4}(-0.966,-0.259)(-0.897,-0.518)(-0.707,-0.707)
\bezier{4}(-0.707,-0.707)(-0.518,-0.897)(-0.259,-0.966)
\bezier{4}(-0.259,-0.966)(0,-1.04)(0.259,-0.966)
\bezier{4}(0.259,-0.966)(0.518,-0.897)(0.707,-0.707)
\bezier{4}(0.707,-0.707)(0.897,-0.518)(0.966,-0.259)
}
\def\Endpoint[#1]{
\ifcase#1
\put(1,0){\circle*{0.15}}
\or\put(0.866,0.5){\circle*{0.15}}
\or\put(0.5,0.866){\circle*{0.15}}
\or\put(0,1){\circle*{0.15}}
\or\put(-0.5,0.866){\circle*{0.15}}
\or\put(-0.866,0.5){\circle*{0.15}}
\or\put(-1,0){\circle*{0.15}}
\or\put(-0.866,-0.5){\circle*{0.15}}
\or\put(-0.5,-0.866){\circle*{0.15}}
\or\put(0,-1){\circle*{0.15}}
\or\put(0.5,-0.866){\circle*{0.15}}
\or\put(0.866,-0.5){\circle*{0.15}}
\fi}
\def\Arc[#1]{
\thicklines			
\ifcase#1
\bezier{25}(0.966,-0.259)(1.04,0)(0.966,0.259)
\or
\bezier{25}(0.966,0.259)(0.897,0.518)(0.707,0.707)
\or
\bezier{25}(0.707,0.707)(0.518,0.897)(0.259,0.966)
\or
\bezier{25}(0.259,0.966)(0,1.04)(-0.259,0.966)
\or
\bezier{25}(-0.259,0.966)(-0.518,0.897)(-0.707,0.707)
\or
\bezier{25}(-0.707,0.707)(-0.897,0.518)(-0.966,0.259)
\or
\bezier{25}(-0.966,0.259)(-1.04,0)(-0.966,-0.259)
\or
\bezier{25}(-0.966,-0.259)(-0.897,-0.518)(-0.707,-0.707)
\or
\bezier{25}(-0.707,-0.707)(-0.518,-0.897)(-0.259,-0.966)
\or
\bezier{25}(-0.259,-0.966)(0,-1.04)(0.259,-0.966)
\or
\bezier{25}(0.259,-0.966)(0.518,-0.897)(0.707,-0.707)
\or
\bezier{25}(0.707,-0.707)(0.897,-0.518)(0.966,-0.259)
\fi}
\def\DottedArc[#1]{
\ifcase#1
\bezier{4}(0.966,-0.259)(1.04,0)(0.966,0.259)
\or
\bezier{4}(0.966,0.259)(0.897,0.518)(0.707,0.707)
\or
\bezier{4}(0.707,0.707)(0.518,0.897)(0.259,0.966)
\or
\bezier{4}(0.259,0.966)(0,1.04)(-0.259,0.966)
\or
\bezier{4}(-0.259,0.966)(-0.518,0.897)(-0.707,0.707)
\or
\bezier{4}(-0.707,0.707)(-0.897,0.518)(-0.966,0.259)
\or
\bezier{4}(-0.966,0.259)(-1.04,0)(-0.966,-0.259)
\or
\bezier{4}(-0.966,-0.259)(-0.897,-0.518)(-0.707,-0.707)
\or
\bezier{4}(-0.707,-0.707)(-0.518,-0.897)(-0.259,-0.966)
\or
\bezier{4}(-0.259,-0.966)(0,-1.04)(0.259,-0.966)
\or
\bezier{4}(0.259,-0.966)(0.518,-0.897)(0.707,-0.707)
\or
\bezier{4}(0.707,-0.707)(0.897,-0.518)(0.966,-0.259)
\fi}
\def\Chord[#1,#2]{
\thinlines
\ifnum#1>#2\Chord[#2,#1]
\else\ifnum#1<#2
\ifcase#1
\ifcase#2
\or\qbezier(1,0)(0.516,0.138)(0.866,0.5)
\or\qbezier(1,0)(0.45,0.26)(0.5,0.866)
\or\qbezier(1,0)(0.327,0.327)(0,1)
\or\qbezier(1,0)(0.179,0.311)(-0.5,0.866)
\or\qbezier(1,0)(0.0536,0.2)(-0.866,0.5)
\or\put(1, 0){\line(-2, 0){2}}
\or\qbezier(1,0)(0.0536,-0.2)(-0.866,-0.5)
\or\qbezier(1,0)(0.179,-0.311)(-0.5,-0.866)
\or\qbezier(1,0)(0.327,-0.327)(0,-1)
\or\qbezier(1,0)(0.45,-0.26)(0.5,-0.866)
\or\qbezier(1,0)(0.516,-0.138)(0.866,-0.5)
\fi
\or\ifcase#2\or
\or\qbezier(0.866,0.5)(0.378,0.378)(0.5,0.866)
\or\qbezier(0.866,0.5)(0.26,0.45)(0,1)
\or\qbezier(0.866,0.5)(0.12,0.446)(-0.5,0.866)
\or\qbezier(0.866,0.5)(0,0.359)(-0.866,0.5)
\or\qbezier(0.866,0.5)(-0.0536,0.2)(-1,0)
\or\put(0.866, 0.5){\line(-5, -3){1.73}}
\or\qbezier(0.866,0.5)(0.146,-0.146)(-0.5,-0.866)
\or\qbezier(0.866,0.5)(0.311,-0.179)(0,-1)
\or\qbezier(0.866,0.5)(0.446,-0.12)(0.5,-0.866)
\or\qbezier(0.866,0.5)(0.52,0)(0.866,-0.5)
\fi
\or\ifcase#2\or\or
\or\qbezier(0.5,0.866)(0.138,0.516)(0,1)
\or\qbezier(0.5,0.866)(0,0.52)(-0.5,0.866)
\or\qbezier(0.5,0.866)(-0.12,0.446)(-0.866,0.5)
\or\qbezier(0.5,0.866)(-0.179,0.311)(-1,0)
\or\qbezier(0.5,0.866)(-0.146,0.146)(-0.866,-0.5)
\or\put(0.5, 0.866){\line(-3, -5){1}}
\or\qbezier(0.5,0.866)(0.2,-0.0536)(0,-1)
\or\qbezier(0.5,0.866)(0.359,0)(0.5,-0.866)
\or\qbezier(0.5,0.866)(0.446,0.12)(0.866,-0.5)
\fi
\or\ifcase#2\or\or\or
\or\qbezier(0,1.)(-0.138,0.516)(-0.5,0.866)
\or\qbezier(0,1.)(-0.26,0.45)(-0.866,0.5)
\or\qbezier(0,1.)(-0.327,0.327)(-1,0)
\or\qbezier(0,1.)(-0.311,0.179)(-0.866,-0.5)
\or\qbezier(0,1.)(-0.2,0.0536)(-0.5,-0.866)
\or\put(0, 1){\line(0, -2){2}}
\or\qbezier(0,1.)(0.2,0.0536)(0.5,-0.866)
\or\qbezier(0,1.)(0.311,0.179)(0.866,-0.5)
\fi
\or\ifcase#2\or\or\or\or
\or\qbezier(-0.5,0.866)(-0.378,0.378)(-0.866,0.5)
\or\qbezier(-0.5,0.866)(-0.45,0.26)(-1,0)
\or\qbezier(-0.5,0.866)(-0.446,0.12)(-0.866,-0.5)
\or\qbezier(-0.5,0.866)(-0.359,0)(-0.5,-0.866)
\or\qbezier(-0.5,0.866)(-0.2,-0.0536)(0,-1)
\or\put(-0.5, 0.866){\line(3, -5){1}}
\or\qbezier(-0.5,0.866)(0.146,0.146)(0.866,-0.5)
\fi
\or\ifcase#2\or\or\or\or\or
\or\qbezier(-0.866,0.5)(-0.516,0.138)(-1,0)
\or\qbezier(-0.866,0.5)(-0.52,0)(-0.866,-0.5)
\or\qbezier(-0.866,0.5)(-0.446,-0.12)(-0.5,-0.866)
\or\qbezier(-0.866,0.5)(-0.311,-0.179)(0,-1)
\or\qbezier(-0.866,0.5)(-0.146,-0.146)(0.5,-0.866)
\or\put(-0.866, 0.5){\line(5, -3){1.73}}
\fi
\or\ifcase#2\or\or\or\or\or\or
\or\qbezier(-1,0)(-0.516,-0.138)(-0.866,-0.5)
\or\qbezier(-1,0)(-0.45,-0.26)(-0.5,-0.866)
\or\qbezier(-1,0)(-0.327,-0.327)(0,-1)
\or\qbezier(-1,0)(-0.179,-0.311)(0.5,-0.866)
\or\qbezier(-1,0)(-0.0536,-0.2)(0.866,-0.5)
\fi
\or\ifcase#2\or\or\or\or\or\or\or
\or\qbezier(-0.866,-0.5)(-0.378,-0.378)(-0.5,-0.866)
\or\qbezier(-0.866,-0.5)(-0.26,-0.45)(0,-1)
\or\qbezier(-0.866,-0.5)(-0.12,-0.446)(0.5,-0.866)
\or\qbezier(-0.866,-0.5)(0,-0.359)(0.866,-0.5)
\fi
\or\ifcase#2\or\or\or\or\or\or\or\or
\or\qbezier(-0.5,-0.866)(-0.138,-0.516)(0,-1)
\or\qbezier(-0.5,-0.866)(0,-0.52)(0.5,-0.866)
\or\qbezier(-0.5,-0.866)(0.12,-0.446)(0.866,-0.5)
\fi
\or\ifcase#2\or\or\or\or\or\or\or\or\or
\or\qbezier(0,-1.)(0.138,-0.516)(0.5,-0.866)
\or\qbezier(0,-1.)(0.26,-0.45)(0.866,-0.5)
\fi
\or\ifcase#2\or\or\or\or\or\or\or\or\or\or
\or\qbezier(0.5,-0.866)(0.378,-0.378)(0.866,-0.5)
\fi\fi\fi\fi}
\def\FullChord[#1,#2]{
\Endpoint[#1]
\Endpoint[#2]
\Arc[#1]
\Arc[#2]
\Chord[#1,#2]
}
\def\EndChord[#1,#2]{
\Endpoint[#1]
\Endpoint[#2]
\Chord[#1,#2]
}
\def\Picture#1{
\begin{picture}(2,1)(-1,-0.167)
#1
\end{picture}
}
\def\DottedChordDiagram[#1,#2]{
\Picture{\DottedCircle \FullChord[#1,#2]}
}
\def\u#1{ \underline{#1} }
\def\uu#1{ \underline{\underline{#1}} }

\def\ux{ {\u{x}} }
\def\uy{ {\u{y}} }
\def\ua{ {\u{a}} }
\def\ual{ {\u{\a}} }
\def\ug{ {\u{g}} }
\def\um{ {\u{m}} }
\def\ut{ {\u{t}} }

\def\ut{ {\u{t}} }
\def\uo{ {\u{o}} }

\def\uux{ {\uu{x}} }

\def\unu{ {\uu{\nu}} }

\def\vx{ \vec{x} }
\def\vy{ \vec{y} }
\def\va{ \vec{a} }
\def\vb{ \vec{b} }
\def\vua{ \vec{\ua} }

\def\val{ \vec{\a} }
\def\vn{ \vec{n} }

\def\val{ \vec{\a} }

\def\tcL{ \tilde{\cL} }
\def\tL{ \tilde{L} }

\def\tux{ \tilde{\ux} }
\def\tx{ \tilde{x} }

\def\prb#1{ \{ #1 \} }

\def\pjoL{ \prod_{j=1}^L }
\def\pioL{ \prod_{i=1}^L }
\def\pjoN{ \prod_{j=1}^N }

\def\sjoL{ \sum_{j=1}^L }

\def\snzi{ \sum_{n=0}^\infty }

\def\sjoLp{ \sum_{j=1}^{L\p} }

\def\smti{ \sum_{m=2}^\infty }
\def\smtN{ \sum_{m=2}^N }
\def\snzi{ \sum_{n=0}^\infty }

\def\sioL{ \sum_{i=1}^L }
\def\sioLp{ \sum_{i=1}^{L\p} }

\def\sizL{ \sum_{i=0}^L }

\def\smzi{ \sum_{m=0}^\infty }

\def\Jbas#1#2#3{ J^{#2}_{#1}(#3;q) }

\def\Zbas#1#2#3{ Z^{#2}_{#1}(#3;K) }
\def\Ztrbas#1#2{ \Zbas{#1}{\tr}{#2} }

\def\ZtrM{ \Ztrbas{}{M} }

\def\Pbas#1#2#3{ \Pbasa{#1}{#2;#3} }

\def\yQ{ Q }
\def\Qlmnaz{ \yQ_{l,\um,n}(\etpaz,\etpua) }
\def\Qlmntz{ \yQ_{l,\um,n}(t_0,\ut) }
\def\Qnaz{ \yQ_{\ube;n}(\etpaz,\etpua) }
\def\Qntz{ \yQ_{\ube;n}(t_0,\ut) }
\def\Qnut{ \yQ_{\ube;n}(\ut) }
\def\Qnutp{ \yQ\p_{\ube;n}(\ut) }
\def\Qnua{ \yQ_{\ube;n}(\etpua) }
\def\Qnuao{ \yQ_{\ube;n}(\etpao) }
\def\Qnuto{ \yQ_{\ube;n}(t_1) }
\def\Qnuap{ \yQ\p_{\ube;n}(\etpua) }

\def\yq{ q }
\def\qmumn{ \yq_{m_0,\um;n}(\ube) }
\def\qumn{ \yq_{\um;n}(\ube) }

\def\ZZbas#1{ \ZZ[#1] }

\def\ZZutiz{ \ZZbas{t_0,t_0^{-1},\ut,\ut^{-1}} }

\def\ZZtoohom{ \ZZbas{t_1^{\pm 1/o_1},1/\ho{M}} }

\def\lk#1#2{ {\rm lk}(#1;#2) }

\def\DumnL{ \Dumn(\cL) }

\def\LmL#1{ \LmLbas{\cL}{#1} }

\def\Pi#1{ P_{#1} }

\def\brlst#1#2#3{ \{#1_{#2},\ldots,#1_{#3} \} }

\def\cLbas#1{ \cL_{[#1]} }
\def\cLk{ \cLbas{k} }

\def\cLp{ \cL\p }
\def\cLz{ \cL_{0} }
\def\cLs{ \cL_{s} }

\def\hi{ h^{-1} }

\def\hlfv{ {1\over 2} }

\def\lij{ l_{ij} }

\def\tpi{ 2\pi i }

\def\Ninf{ N\rightarrow \infty }
\def\Kinf{ K\rightarrow \infty }

\def\etpiKlkbas#1#2{ e^{\hpiKc\lk{#1}{#2}} }

\def\etpiKlkMLa{ \etpiKlkbas{M,\cL}{\ua} }

\def\lui#1#2{ l^{(#1,#2)} }
\def\laa{ \lui{a}{a} }
\def\lbb{ \lui{b}{b} }
\def\lab{ \lui{a}{b} }

\def\lcw{ \l_{\rm CW} }

\def\fhom{ H_1(M,\ZZ) }
\def\ordH{ |\fhom| }

\def\ust{ u^{\st} }
\def\uust{ \uuu^{\st} }

\def\etpi#1{ e^{\tpi#1} }

\def\skmn{ s_{k,m,n} }

\def\ordHM#1{ |H_1(#1,\ZZ)| }

\def\slim{ \mathop{\rm slim}\limits  }

\def\slimNinf{ \slim_{N\liminf} }

\def\ube{ \u{\b} }
\def\vube{ \vec{\ube} }

\def\uc{ \u{c} }
\def\vuc{ \vec{\uc} }
\def\vc{ \vec{c} }

\def\ug{ \u{\gamma} }

\def\ump{ \um\p }

\def\vual{ \u{\val} }
\def\val{ \vec{\a} }

\def\vux{ \vec{\ux} }
\def\xb{ \bar{x} }

\def\yp{ p }

\def\URC{$U(1)$-RC\ }
\def\TCC{TCC\ }

\def\ro#1{~(#1,I)}
\def\row#1{~#1(I)}

\def\liminf{ \rightarrow \infty }
\def\limz{ \rightarrow 0 }

\def\cLs{ \cL^s }
\def\Ls{ L^s }
\def\hcLs{ \hcL^s }
\def\cLLs{ \cL\cup\cLs }
\def\hcK{ \hat{\cK} }
\def\hcL{ \hat{\cL} }
\def\cLp{ \cL\p }
\def\cLLp{ \cL\cup\cLp }
\def\cLb{ \bar{\cL} }
\def\cLhcLs{ \cL\cup\hcLs }

\def\Pr{ P }

\def\oh{ o }

\def\HoQ#1{ H_1(#1;\IQ) }
\def\HoZ#1{ H_1(#1;\ZZ) }
\def\ordH#1{ |\HoZ{#1}| }

\def\ordHM{ \ordH{M} }
\def\ho#1{ h_1(#1) }
\def\hoe#1#2{ h_1^{#1}(#2) }
\def\hotM{ \hoe{-3/2}{M} }
\def\hoM{ \hoe{-1/2}{M} }
\def\hi{ h^{-1} }

\def\dtlks{ \det \lkbasm{\hcLs}{M} }
\def\dtlks{ \det(\hcLs|M) }
\def\sgnlks{ \sigl{\hcLs}{S^3} }

\def\ell{,\ldots,}
\def\ellp#1#2{ (#1\ell #2) }

\def\hlfv{ {1\over 2} }
\def\qrtv{ {1\over 4} }

\def\lij{ l_{ij} }
\def\lli#1{ l^{(#1)} }
\def\llim#1#2{ \lli{#1|#2} }
\def\lLM{ \llim{\cL}{M} }
\def\llimzo{ \llim{\cLzL}{M}_{01} }
\def\lzo{ l_{01} }

\def\llst{ \llim{\hcLs}{S^3} }

\def\lkw{ \mathop{\rm lk} }
\def\lkbas#1{ \lkw( #1 ) }
\def\lkbasm#1#2{ \lkbas{ #1|#2 } }

\def\lkbasvm#1#2#3{ \lkw(#1,#2;#3) }
\def\lkLMua{ \lkbasvm{\cL}{M}{\ua} }
\def\lkLMuaz{ \lkbasvm{\cLzL}{M}{a_0,\ua} }

\def\lkLLsSac{ \lkbasvm{\cLhcLs}{S^3}{\ua,\uc} }
\def\lkLLsSacz{ \lkbasvm{\cLzLhs}{S^3}{a_0,\ua,\uc} }

\def\lkLLsSacst{ \lkbasvm{\cLhcLs}{S^3}{\ua,\ucst(\ua)} }
\def\lkLsS{ \lkbasvm{\hcLs}{S^3}{\uc} }
\def\lkLzMa{ \lkbasvm{\hcLzL}{M}{a_0,\ua} }

\def\mtibas#1#2{ (#1_1,\ldots,#1_{#2}) }

\def\GCF#1{ \mathop{\rm GCF}\lrbc{#1} }
\def\LCM#1{ \mathop{\rm LCM}\lrbc{#1} }

\def\soiljL{ \sum_{1\leq i<j\leq L} }
\def\soiljLp{ \sum_{1\leq i<j\leq L\p} }
\def\soinjL{ \sum_{1\leq i\neq j\leq L} }

\def\somnLs{ \sum_{1\leq k,k\p\leq L} }
\def\skoLs{ \sum_{k=1}^{\Ls} }
\def\sjoL{ \sum_{j=1}^L }
\def\sioL{ \sum_{i=1}^L }

\def\snzi{ \sum_{n=0}^\infty }
\def\snoi{ \sum_{n=1}^\infty }
\def\snzN{ \sum_{n=0}^N }
\def\sumn{ \sum_{\um,n\geq 0 \atop |\um|\leq n} }
\def\summpn{ \sum_{\um,\um\p,n\geq 0 \atop |\um|+|\um\p|\leq n} }
\def\sumnN{ \sum_{0\leq \um,n\leq N \atop |\um|\leq n} }
\def\smti{ \sum_{m\geq 2} }
\def\smtN{ \sum_{2\leq m\leq N} }
\def\smnzi{ \sum_{m,n\geq 0} }
\def\sugoK{ \sum_{0 < \ug < K} }
\def\sjoLs{ \sum_{j=1}^{\Ls} }
\def\sumgz{ \sum_{\um\geq 0} }
\def\summzgz{ \sum_{m_0,\um\geq 0} }
\def\smu{ \sum_{\mu=\pm 1} }
\def\smumu{ \smu\mu }
\def\smnt{ \sum_{m\geq 0\atop n\geq - m/2} }
\def\szmn{ \sum_{0\leq m\leq n} }
\def\sioLnj{ \sum_{1\leq i\leq L\atop i\neq j} }
\def\sjoLni{ \sum_{1\leq j\leq L\atop j\neq i} }
\def\skmn{ \sum_{k,m\geq 0\atop n\geq \qrtv k - {3\over 4} m} }
\def\sklmn{ \sum_{k,l,m\geq 0\atop n\geq \qrtv k
   + \hlfv l- {3\over 4} m} }
\def\slmn{ \sum_{\um\geq 0\atop 0\leq l\leq {|\um|+1\over 2} } }
\def\slmmn{ \sum_{m\geq 0\atop 0\leq l\leq {m+1\over 2} } }
\def\sizL{ \sum_{i=0}^L }
\def\siLjLp{ \sum_{1\leq i\leq L\atop 1\leq j\leq L\p} }

\def\sojLp{ \sum_{j=1}^{L\p} }

\def\pioL{ \prod_{i=1}^L }
\def\pioLj{ \prod_{1\leq i\leq L\atop i\neq j} }
\def\pjoLp{ \prod_{j=1}^{L\p} }

\def\AP{ \Delta_{\rm A} }
\def\AF{ \nabla_{\rm A} }
\def\APbas#1#2{ \AP(#1;#2) }
\def\AFbas#1#2{ \AF(#1;#2) }
\def\APpbas#1#2#3{ \AP^{#1}(#2;#3) }
\def\AFpbas#1#2#3{ \AF^{#1}(#2;#3) }

\def\APpLzq{ \APpbas{2n+1}{\cLz,M}{q} }

\def\APpnLM{ \APpbas{2n}{\cL,M}{t_1} }
\def\APpnoLMa{ \APpbas{2n+1}{\cL,M}{\etpao} }
\def\AFpnLM{ \AFpbas{2n}{\cL,M}{\ut} }
\def\AFpnoLMaz{ \AFpbas{2n+1}{\cLzLhL,M}{\etpaz,\etpua,\etpi{\uc}} }
\def\AFpnoLSaz{ \AFpbas{2n+1}{\cLzLhL,S^3}{\etpaz,\etpua,\etpi{\uc}} }

\def\AFLMut{ \AFbas{\cL,M}{\ut} }
\def\AFLMput{ \AFbas{\cL,M\p}{\ut} }

\def\AFLzMua{ \AFbas{\cLzL,M}{\etpaz,\etpua} }
\def\AFLnozMua{ \AFpbas{2n+1}{\cLzL,M}{\etpaz,\etpua} }
\def\AFLnqmMua{ \AFpbas{2n+1}{\cLzL,M}{q^{\mu},\etpua} }

\def\AFLzMuto{ \AFbas{\cLzL,M}{1,\ut} }
\def\AFLzMuao{ \AFbas{\cLzL,M}{1,\etpua} }
\def\AFLMcut{ \AFbas{\cL,M\# M\p}{\ut} }
\def\AFLMua{ \AFbas{\cL,M}{\etpua} }
\def\AFLSut{ \AFbas{\cL,S^3}{\ut} }
\def\APKMt{ \APbas{\cK,M}{t} }

\def\AFLst{ \AFbas{\cLLs,S^3}{\etpi{\ua},\etpi{\ucst(\ua)}} }

\def\etpua{ \etpi{\ua} }
\def\etpaz{ \etpi{a_0} }
\def\etpao{ \etpi{a_1} }
\def\etpat{ \etpi{a_2} }

\def\ZZthopm{ \ZZ[t^{1/2\oh(\cK)}, t^{-1/2\oh(\cK)}] }
\def\ZZuthopm{ \ZZ[\ut^{1/2\uo}, \ut^{-1/2\uo} ] }
\def\ZZtopm{ \ZZ[t^{1/\oh(\cK)}, t^{-1/\oh(\cK)}] }
\def\ZZutopm{ \ZZ[\ut^{1/\uo}, \ut^{-1/\uo} ] }
\def\ZZtpmtoh{ \ZZ[\ut^{\pm {1\over 2\uo}}, {1\over 2\ho{M}} ] }
\def\ZZttpmtoh{ \ZZ[\ut^{\pm 1/\uo}, {1\over 2\ho{M}} ] }
\def\ZZtpmoh{ \ZZ[\ut^{\pm 1/\uo}, 1/\ho{M} ] }
\def\ZZtopmoh{ \ZZ[t^{\pm 1/o_1}, 1/\ho{M} ] }
\def\ZZtpmohz{ \ZZ[t_0^{\pm 1},\ut^{\pm 1/\uo},1/\ho{M} ] }
\def\ZZihoM{ \ZZ[1/\ho{M}] }
\def\ZZihoMuti{ \ZZ[1/\ho{M},\ut^{\pm 1/\uo} ] }
\def\ZZutus{ \ZZ[\ut^{\pm 1/\uo},(\uusr)^{\pm 1/\uo\p},1/\ho{M}] }
\def\ZZihoMh{ \ZZbas{1/\ho{M},1/2} }

\def\ZZihoMh{ \ZZihoM[[h]] }

\def\IQub{ \IQ[\ube] }
\def\IQubh{ \IQ[\ube][[h]] }

\def\cxx{ \IC[[x]] }

\def\fdhbas#1#2{ {#1}^{#2} - {#1}^{-#2} }
\def\fdh#1{ {#1}^{1/2} - {#1}^{-1/2} }
\def\fdhrbas#1#2{ {#1}^{#2/2} - {#1}^{-#2/2} }

\def\est{ ^{\rm (st)} }
\def\efr{ ^{\rm (fr)} }
\def\etr{ ^{\rm (tr)} }
\def\er{ ^{\rm (r)} }

\def\ust{ u\est }
\def\uuv{ \u{u} }
\def\uust{ \uuv\est }
\def\uo{ \u{\oh} }
\def\unu{ \u{\nu} }
\def\up{ \u{p} }
\def\us{ \u{s} }

\def\uusr{ \uuv^* }

\def\nup{ \nu\p }
\def\kapp{ \kappa }
\def\unup{ \u{\nu}\p }

\def\Phbas#1#2{ \Phi(#1;#2) }
\def\Phpbas#1#2#3{ \Phi^{#1}(#2;#3) }

\def\PhLsut{ \Phbas{\cL,\hcLs}{\ut} }

\def\Jbas#1#2{ J_{#1}(#2;q) }
\def\JuaL{ \Jbas{\ual}{\cL} }
\def\JuabLp{ \Jbas{\ual,\ube}{\cL,\cLp} }
\def\JubLp{ \Jbas{\ube}{\cLp} }

\def\Ibas#1#2#3{ I_{#1}(#2;#3;K) }
\def\IbasN#1#2#3{ I_{#1}(#2;#3;K|N) }
\def\IabLp{ \Ibas{\ube}{\cL,\cLp}{\ua} }
\def\IabLpM{ \Ibas{\ube}{\cL,\cLp,M}{\ua} }
\def\IabLpMN{ \IbasN{\ube}{\cL,\cLp,M}{\ua} }

\def\Ih{ \hat{I} }
\def\Irbas#1#2#3{ I\er_{#1}(#2;#3;K) }
\def\Ihrbas#1#2#3{ \hat{I}\er_{#1}(#2;#3;K) }

\def\ItrbasN#1#2#3{ \tilde{I}\er_{#1}(#2;#3;K|N) }
\def\IrbasN#1#2#3{ I\er_{#1}(#2;#3;K|N) }
\def\IrbasNN#1#2#3{ I\er_{#1}(#2;#3;K||N) }
\def\IrabLpM{ \Irbas{\ube}{\cL,\cLp,M}{\ua} }
\def\IrabLpMp{ \Irbas{\ube}{\cL,\cLp,M\p}{\ua} }
\def\IhrabLpM{ \Ihrbas{\ube}{\cL,\cLp,M}{\ua} }
\def\IhraobLpM{ \Ihrbas{\ube}{\cL,\cLp,M}{a_1} }
\def\IhrabLpMz{ \Ihrbas{\ube}{\cLzL,\cLp,M}{a_0,\ua} }
\def\IhrabLzpM{ \Ihrbas{\ube}{\cLzLs,\cLp,M}{a_0,\ua,\uc} }
\def\IrabLzpM{ \Irbas{\ube}{\cLzL,\cLp,M}{a_0,\ua} }

\def\IraobLpM{ \Irbas{\ube}{\cL,\cLp,M}{a_1} }
\def\IaobLpM{ \Ibas{\ube}{\cL,\cLp,M}{a_1} }
\def\IrabLpMNt{ \ItrbasN{\ube}{\cL,\cLp,M}{\ua} }
\def\IrabLpMN{ \IrbasN{\ube}{\cL,\cLp,M}{\ua} }
\def\IrabLpMNN{ \IrbasNN{\ube}{\cL,\cLp,M}{\ua} }
\def\IrmabLpMN{ \IrbasN{\ube}{\cL,\cLp,M}{a_1\ell\mu a_j\ell a_L} }
\def\IrabcLpsM{ \Irbas{\ube}{\cLLs,\cLp,M}{\ua,\uc} }
\def\IrabcLpsS{ \Irbas{\ube}{\cLLs,\cLp,S^3}{\ua,\uc} }

\def\IrabcLpsSN{ \IrbasN{\ube}{\cLLs,\cLp,S^3}{\ua,\uc} }
\def\IrabcLpsSNh{ \IrbasN{\ube}{\cLhcLs,\cLp,S^3}{\ua,\uc} }

\def\IrabLpMj{ \Irbas{\a_j,\ube}{\cLrj,\cLjLp,M}{\uarj} }
\def\IrabLpMNj{ \IrbasN{\a_j,\ube}{\cLrj,\cLjLp,M}{\uarj} }

\def\IrabLpS{ \Irbas{\ube}{\cL,\cLp,S^3}{\ua} }
\def\IrabLpSe{ I\er_{\ube}(\empt,\cLp,S^3;K) }
\def\IrabLpMe{ I\er_{\ube}(\empt,\cLp,M;K) }
\def\IhrabLpMe{ \Ih\er_{\ube}(\empt,\cLp,M;K) }

\def\Ihrasub{ \Ihrbas{\ube}{\cLz,\cLp,M}{1/K} }

\def\cJr{ \check{J}\er }
\def\cJrbas#1#2{ \cJr_{#1}(#2;h) }
\def\cJrLSa{ \cJrbas{\ube}{\cL,\cLp;\etpi{\ua}} }
\def\cJrLSt{ \cJrbas{\ube}{\cL,\cLp;\ut} }
\def\cJrLSe{ \cJrbas{\ube}{\empt,\cLp} }

\def\cLjLp{ \cL_j\cup\cLp }

\def\etpi#1{ e^{\tpi{#1}} }
\def\ehpi#1{ e^{\hlfv i\pi #1} }

\def\hpiK{ \hlfv\, i\pi K }
\def\hpiKc{ \hlfv i\pi K }
\def\hpi{ \hlfv\,i\pi }

\def\etpcl{ \exp \lrbc{\hpiK\sjoLs\llim{\hcLs}{S^3}_{jj}c_j^2 } }
\def\etpclM{ \exp \lrbc{\hpiK\sjoLs\llim{\hcLs}{M}_{jj}c_j^2 } }

\def\Zbas#1#2{ Z_{#1}(#2;K) }
\def\ZuaLM{ \Zbas{\ual}{\cL,M} }

\def\Ztrbas#1#2{ Z\etr_{#1}(#2;K) }
\def\ZtruaLM{ \Ztrbas{\ual}{\cL,M} }
\def\ZtruajLM{ \Ztrbas{\ual_{(j)}}{\cL_{(j)},M} }
\def\ZtruabLLpM{ \Ztrbas{\ual,\ube}{\cL,\cLp,M} }
\def\ZtrubLpM{ \Ztrbas{\ube}{\cLp,M} }

\def\Zh{ \hat{Z} }
\def\Zhr{ \Zh\er }
\def\Zhrbas#1#2{ \Zhr_{#1}(#2;h) }
\def\Zhrebas#1#2#3{ \Zhrbas{#1}{#2;#3} }
\def\ZhrubLLpMuth{ \Zhrebas{\ube}{\cL,\cLp,M}{\ut} }
\def\ZhrubLLpSuth{ \Zhrebas{\ube}{\cL,\cLp,S^3}{\ut} }
\def\Zhrem{ \Zhrbas{\ube}{\empt,\cLp,M} }

\def\empt{ \varnothing }

\def\frbas#1#2{ \phi\efr_{#1}(#2) }
\def\frbas#1#2{ \phi\efr(#2;#1) }
\def\qfrbas#1#2{ q^{\frbas{#1}{#2}} }
\def\qfrua{ \qfrbas{\ual}{\cL,\hcLs} }
\def\qefrua{ e^{(2\pi i/K)\, \frbas{\ual}{\cL,\hcLs} } }

\def\fruabpKx{ \frbas{K|\vua|,\ube}{\cLLp,\hcLs} }
\def\fruabpKy{ \frbas{K\ua,\ube}{\cLLp,\hcLs} }
\def\fruabpKyM{ \frbas{K\ua,\ube}{\cLLp,\hcLs,M} }
\def\fruabpKyz{ \frbas{Ka_0,K\ua,\ube}{\cLzLp,\hcLs} }

\def\exfruabpKx{ e^{\tpioK \,\fruabpKx} }
\def\exfruabpKy{ e^{\tpioK \,\fruabpKy} }
\def\exfruabpKyM{ e^{\tpioK \,\fruabpKyM} }
\def\frua{ \frbas{\ual}{\cL,\hcLs} }
\def\fruaM{ \frbas{\ual}{\cL,\hcLs,M} }

\def\phlk#1#2{ \phi_{\rm lk}(#1;#2) }

\def\phlkLpM{ \phlk{\cL,\cLp,M}{\ube} }
\def\phlkLpMe{ \phlk{\empt,\cLp,M}{\ube} }
\def\phlkLzpM{ \phlk{\cLz,\cLp,M}{\ube} }
\def\phlkLpzM{ \phlk{\cLzL,\cLp,M}{\ube} }
\def\phlkLpsS{ \phlk{\cLLs,\cLp,S^3}{\ube} }
\def\qphlkLpM{ q^{\phlkLpM} }
\def\qphlkLpMe{ q^{\phlkLpMe} }
\def\qphlkLzpM{ q^{\phlkLzpM} }
\def\qphlkLpzM{ q^{\phlkLpzM} }

\def\tpioK{ (2\pi i/K) }

\def\sigl#1#2{ \prosign(#1|#2) }

\def\Dbas#1#2{ D_{#1}(#2) }
\def\DumnL{ \Dbas{\um;n}{\cL} }
\def\DumnLM{ \Dbas{\um;n}{\cL,M} }
\def\DummpnL{ \Dbas{\um,\ump;n}{\cLLp} }

\def\Pbas#1#2{ P_{#1}(#2) }
\def\Ppbas#1#2{ P\p_{#1}(#2) }
\def\Ptbas#1#2{ \tilde{P}_{#1}(#2) }
\def\Pebas#1#2#3{ \Pbas{#1}{#2;#3} }
\def\Ppebas#1#2#3{ \Ppbas{#1}{#2;#3} }
\def\PubnLLpMut{ \Pebas{\ube;n}{\cL,\cLp,M}{\ut} }
\def\PpubnLLpMut{ \Ppebas{\ube;n}{\cL,\cLp,M}{\ut} }
\def\PpubnLLpMutz{ \Ppebas{\ube;n}{\cLzL,\cLp,M}{t_0,\ut} }
\def\PubnLLpMuaz{ \Pebas{\ube;n}{\cLzLhL,\cLp,M}{\etpaz,\etpua,
      \etpi{\uc}} }
\def\PpubnLLpMuaz{ \Ppebas{\ube;n}{\cLzL,\cLp,M}{\etpaz,\etpua,} }
\def\PpubnLLpSuaz{ \Ppebas{\ube;n}{\cLzLhL,\cLp,S^3}{\etpaz,\etpua,
      \etpi{\uc}} }
\def\PubnLLpMt{ \Pebas{\ube;n}{\cL,\cLp,M}{t_1} }
\def\PpubnLLpMt{ \Ppebas{\ube;n}{\cL,\cLp,M}{t_1} }
\def\PpubnLzLpMq{ \Ppebas{\ube;n}{\cLz,\cLp,M}{q} }
\def\PubnLp{ \Pbas{\ube;n}{\cLp} }
\def\PpubnLp{ \Ppbas{\ube;n}{\cLp} }
\def\PeLpMt{ \Pebas{\ube;0}{\cL,\cLp,M}{\ut} }
\def\PeLpMua{ \Pebas{\ube;0}{\cL,\cLp,M}{\etpua} }

\def\pbas#1#2{ p_{#1}(#2) }
\def\pebas#1#2#3{ \pbas{#1}{#2;#3} }
\def\pumnLLpMub{ \pebas{\um;n}{\cL,\cLp,M}{\ube} }
\def\ppbas#1#2{ p\p_{#1}(#2) }
\def\ppebas#1#2#3{ \ppbas{#1}{#2;#3} }
\def\ppumnLLpMub{ \ppebas{\um;n}{\cL,\cLp,M}{\ube} }

\def\Lbas#1#2#3 { L_{#1}(#2;#3) }

\def\LmL{ \Lbas{m}{\cL}{\vua} }

\def\LmLM{ \Lbas{m}{\cL,M}{\vua} }
\def\LmnLM{ \Lbas{m,n}{\cL,M}{\vua} }

\def\LtLal{ \Lbas{2}{\cL}{\vual} }
\def\LmLal{ \Lbas{m}{\cL}{\vual} }
\def\LmnLal{ \Lbas{m,n}{\cL}{\vual} }
\def\LtLMa{ \Lbas{2}{\cL,M}{\vua} }
\def\LmLMa{ \Lbas{m}{\cL,M}{\vua} }

\def\LmnLpalb{ \Lbas{m,n}{\cL,\cLp}{\vual;\ube} }
\def\LmnLpab{ \Lbas{m,n}{\cL,\cLp}{\vua;\ube} }

\def\LmnLMpab{ \Lbas{m,n}{\cL,\cLp,M}{\vua;\ube} }

\def\bE{ {\bf E} }
\def\bP{ {\bf P} }

\def\bEbas#1#2{ \bE(#1;#2) }
\def\bPbas#1#2{ \bP(#1;#2;K) }
\def\bPNbas#1#2{ \bP(#1;#2;K|N) }
\def\bPNbasn#1#2#3{ \bP_{#3}(#1;#2;K|N) }
\def\bPrNbasn#1#2#3{ \bP\er_{#3}(#1;#2;K|N) }
\def\bPNbasnx#1#2#3{ \bP_{#3}(#1;#2|N) }
\def\bPsNbasn#1#2#3{ \tilde{\bP}_{#3}(#1;#2;K|N) }
\def\bPsNbasnx#1#2#3{ \tilde{\bP}_{#3}(#1;#2|N) }

\def\bELa{ \bEbas{\cL}{\vua} }
\def\bELMa{ \bEbas{\cL,M}{\vua} }
\def\bELjMa{ \bEbas{\cLrj,M}{\vuarj} }
\def\bELMaN{ \bEbas{\cL,M}{\vua\,|N} }

\def\bELLsSacN{ \bEbas{\cLLs,S^3}{\vua,\vuc\,|N} }
\def\bELMajN{ \bEbas{\cLrj,M}{\vuarj\,|N} }
\def\bELaN{ \bEbas{\cL}{\vua\,|N} }
\def\bELsc{ \bEbas{\cLs}{\vuc} }
\def\bELLsc{ \bEbas{\cLLs}{\vua,\vuc} }
\def\bELLscN{ \bEbas{\cLLs}{\vua,\vuc\,|N} }
\def\bELscN{ \bEbas{\cLs}{\vuc\,|N} }
\def\etbELa{ e^{\hpiKc \bELa} }
\def\etbELMa{ e^{\hpiKc \bELMa} }
\def\etbELMaN{ e^{\hpiKc \bELMaN} }
\def\etbELMajN{ e^{\hpiKc \bELMajN} }
\def\etbA{ e^{i\pi K \bALj\cdot\va_j} }
\def\bPLpab{ \bPbas{\cL,\cLp}{\vua;\ube} }
\def\bPLpMab{ \bPbas{\cL,\cLp,M}{\vua;\ube} }
\def\bPLpMabN{ \bPNbas{\cL,\cLp,M}{\vua;\ube} }
\def\bPjLpMabN{ \bPNbasn{\cL,\cLp,M}{\vuarj;a_j;\ube}{(j)} }
\def\bPjLpMmabN{ \bPNbasn{\cL,\cLp,M}{\vuarj;\mu a_j;\ube}{(j)} }
\def\bPrjLpMabN{ \bPrNbasn{\cL,\cLp,M}{\vuarj;a_j;\ube}{(j)} }
\def\bPjLpMabNkmn{ \bPNbasnx{\cL,\cLp,M}{\vuarj;\ube}{(j)\,k,m,n} }
\def\bPsjLpMabN{ \bPsNbasn{\cL,\cLp,M}{\vua;\ube}{(j)} }
\def\bPsjmn{ \bPsNbasnx{\cL,\cLp,M}{\vua;\ube}{(j)\,m,n} }

\def\bPLpabN{ \bPNbas{\cL,\cLp}{\vua;\ube} }
\def\bPLsca{ \bPbas{\cLs,\cL}{\vuc;\ual} }
\def\bPLLsca{ \bPbas{\cLLs,\cLp}{\vua,\vuc;\ube} }
\def\bPLLscaN{ \bPNbas{\cLLs,\cLp}{\vua,\vuc;\ube} }
\def\bPLLsScaN{ \bPNbas{\cLLs,\cLp,S^3}{\vua,\vuc;\ube} }
\def\bPLscaN{ \bPNbas{\cLs,\cL}{\vuc;\ual} }

\def\vuarj{ \vua\xrem{j} }

\def\bvA{ {\bf \vec{A}} }
\def\bAbasN#1{ \bvA(#1|N) }
\def\bALj{ \bAbasN{\cLrj;\vuarj} }

\def\Ai#1{ A^{(#1)} }
\def\Asi#1{ \check{A}^{(#1)} }
\def\AiuaxN#1{ \Ai{#1}(\ua;\ux\xrem{#1},\uxb\xrem{#1}|N) }
\def\AouaxN{ \AiuaxN{1} }
\def\AnuaN#1{ \Ai{#1}(\ua|N) }
\def\AsoLuaN{ \Asi{1}(\cL,M;\ua|N) }
\def\AsoLjuaN{ \Asi{1}(\cLrj,M;\uarj|N) }
\def\AsjoLuaN{ \Asi{1;j}(\cL,M;\ua|N) }
\def\AniuaN#1#2#3{ \Ai{#1}_{#2}(#3|N) }
\def\AsniuaN#1#2#3{ \Asi{#1}_{#2}(#3|N) }
\def\AoikLuaN{ \AniuaN{1}{ik}{\cL,M;\ua} }
\def\AsoikLuaN{ \AsniuaN{1}{ik}{\cL,M;\ua} }
\def\AsojkLuaN{ \AsniuaN{1}{jk}{\cL,M;\ua} }
\def\AsojjLuaN{ \AsniuaN{1}{jj}{\cL,M;\ua} }
\def\AojjLuaN{ \AniuaN{1}{jj}{\cL,M;\ua} }
\def\AouaN{ \AnuaN{1} }
\def\Aouainv{ \tilde{A}^{(1)}(\ua;\uyo,\uybo|N) }

\def\uaxo{ \ua\xrem{1} }

\def\xF{ F }
\def\xFnuaN{ \xF_n(\ua|N) }

\def\Bmnuax{ B_{m,n}(\ua;\uxo,\uxbo |N) }
\def\Bzuax{ B_{0,0}(\ua;\uxo,\uxbo |N) }
\def\Bmnmuax{ B_{m,n-m}(\ua;\uxo,\uxbo |N) }
\def\Bmnmdel{ B_{m,n-m}(\ua;\del{\uyo},\del{\uybo} |N) }

\def\uxo{ \ux\xrem{1} }
\def\uxbo{ \uxb\xrem{1} }
\def\uxb{ \u{\xb} }
\def\uyo{ \uy\xrem{1} }
\def\uybo{ \uyb\xrem{1} }
\def\uyb{ \u{\yb} }
\def\yb{ \bar{y} }

\def\bspc{ \lrbf{ \sin \pi |\vuc| \over \pi |\vuc| } }
\def\ivc{ \lrbf{d\vuc} }

\def\eval#1#2{ \left. #1 \right|_{#2} }
\def\xrem#1{_{(#1)} }

\def\vuarj{ \vua\xrem{j} }
\def\vuaroj{ \vua\xrem{1,j} }
\def\uaroj{ \ua\xrem{1,j} }

\def\cLrj{ \cL\xrem{j} }

\def\del#1{ \partial_{#1} }

\def\intinf{ \int_{-\infty}^{+\infty} }

\def\iva{ \int_{ |\vua| = \ua } }
\def\ivast{ \int_{ |\vua| = \ua \atop [\vua = \ua \vn] } }
\def\ivasto{ \int_{ |\vua| = \ua
    \atop [\vua\xrem{1} = \ua\xrem{1} \vns] } }
\def\ivastog{ \int_{ {|\vua| = \ua
    \atop [\vua\xrem{1} = \ua\xrem{1} \vns]}\atop
    [\vuc=\vucst(\ua\vns)] } }
\def\icomb{ \int_{ |\vua| = \ua,\;\|\vuc| = \uc
    \atop [\vua\xrem{1} = \ua\xrem{1} \vns,\;\vuc=\uc\vns)]} }
\def\ivastoj{ \int_{ |\vuarj| = \uarj
    \atop [\vua\xrem{1,j} = \ua\xrem{1,j} \vns] } }
\def\ivastj{ \int_{ |\va_j| = \a_j
    \atop [\vua_j = \vajst] } }
\def\ivaj{ \int_{|\va_j| = a_j} }

\def\ivv{ \lrbf{ d\vua \over \ua } }
\def\ivvo{ \lrbf{ d\vua\xrem{1} \over \ua\xrem{1} } }
\def\ivvc{ \lrbf{ d\vuc \over \uc } }
\def\ivvoj{ \lrbf{ d\vuaroj \over \uaroj} }
\def\aex{ \sumn \LmnLpab\,K^{-n} }
\def\aexM{ \sumn \LmnLMpab\,K^{-n} }
\def\aexN{ \sumnN \LmnLpab\,K^{-n} }

\def\vajst{ \va\est_j(\vuarj) }

\def\incs{ \int_{[\vuc = 0]} }
\def\incst{ \int_{[\vuc = \vucstvua]} }
\def\incrs{ \int_{[\uc = \ucstua]} }

\def\ucst{ \uc\est }
\def\ucstua{ \ucst(\ua) }
\def\vucst{ \vuc\est }
\def\vucstvua{ \vucst(\vua) }

\def\piLt{ \lrbc{ \pioL t_i^{\llim{\cL,\cLp}{M}_{ij}  } } }
\def\piLte#1{ \piLt^{#1/2} - \piLt^{-#1/2} }

\def\Ki{ K^{-1} }

\def\IKx{ I^{(\uxst)}(K) }
\def\Ix{ I^{(\uxst)} }

\def\ftxz{ f_{2,\uxst}(\ux) }
\def\xst{ x\est }
\def\uxst{ \ux\est }

\def\Sd{ S_{\rm d} }

\def\Sn{ S }

\def\Sdv#1{ \Sd(#1) }
\def\SdvN#1{ \Sd(#1|N) }
\def\Sdve#1#2{ \Sd^{#1}(#2) }
\def\SdveN#1#2{ \Sd^{#1}(#2|N) }
\def\Snn#1{ \Sn_n(#1) }
\def\Snnp#1{ \Sn\p_n(#1) }
\def\Snz#1{ \Sn_0(#1) }
\def\Snzp#1#2{ \Sn^{#1}_0(#2) }
\def\SnnN#1{ \Sn_n(#1|N) }
\def\SnzN#1{ \Sn_0(#1|N) }
\def\SniN#1#2{ \Sn_{#1}(#2|N) }
\def\Sjklmn{ \SniN{(j)\,k,l,m,n}{\cL,\cLp,M;\uarj;\ube} }
\def\SdLMuaN{ \SdvN{\cL,M;\ua} }
\def\SdLMua{ \Sdv{\cL,M;\ua} }
\def\SdLSua{ \Sdv{\cL,S^3;\ua} }
\def\SdLzMua{ \Sdv{\cLzL,M;a_0,\ua} }
\def\SdLMuaj{ \Sdv{\cLrj,M;\uarj} }
\def\SdLMuajN{ \Sdv{\cLrj,M;\uarj|N} }
\def\SdLMuan{ \Sdve{2n+1}{\cL,M;\ua} }

\def\SdLMuaNn{ \SdveN{2n+1}{\cL,M;\ua} }
\def\SdLMuaNl{ \SdveN{}{\cLrj,M;\uarj} }
\def\SnnLMua{ \Snn{\cL,\cLp,M;\ua;\ube} }

\def\SpnnLMuaz{ \Snnp{\cLzL,\cLp,M;a_0,\ua;\ube} }
\def\SnzLMua{ \Snz{\cL,\cLp,M;\ua;\ube} }
\def\SnzLMuaz#1{ \Snzp{#1}{\cLzL,\cLp,M;a_0,\ua;\ube} }
\def\SnzLSua{ \Snz{\cL,\cLp,S^3;\ua;\ube} }
\def\SnzzLMua{ \Snz{\cLzL,\cLp,M;a_0,\ua;\ube} }
\def\SnnLMuaN{ \SnnN{\cL,\cLp,M;\ua;\ube} }
\def\SnzLMuaN{ \SnzN{\cL,\cLp,M;\ua;\ube} }
\def\SnnLMpuaN{ \Sn_{n\p}(\cL,\cLp,M;\ua;\ube) }

\def\uarj{ \ua\xrem{j} }

\def\Ouao{ O(\ua) }

\def\ucstaz{ \ucst(a_0,\ua) }

\def\vns{ \vn_* }

\def\bcLs{ \mathfrak{L}_* }
\def\Ci#1{ {\cal C}_{#1} }

\def\CL{ \Ci{\cL} }
\def\ec{ ^{(c)} }
\def\Fcnuab{ F\ec_n(\cL,\cLp,M;\ua;\ube) }
\def\Gcua{ G\ec(\cL;\ua) }

\def\cLz{ \cL_0 }
\def\cLzL{ \cL_0\cup\cL }
\def\cLzLs{ \cLz\cup\cL\cup\cLs }
\def\cLzLhs{ \cLz\cup\cL\cup\hcLs }
\def\cLzLhL{ \cLz\cup\cL\cup\cLs }
\def\cLzLp{ \cLzL\cup\cLp }

\def\hcLzL{ \hat{\cL}\cup\hat{\cL}_0 }

\def\sfrLsM{ \sjoLs \llim{\hcLs}{M}_{jj} }

\def\aldf{ i^{-L} (2\pi)^{L-2} }

\def\Lens#1{ L_{#1} }
\def\MK{ M_K }

\def\yF{ F }
\def\yFbas#1#2{ \yF_{#1}(#2) }
\def\yFnqaz{ \yFbas{n}{a_0,\ua;q} }

\def\eqpsgns{ e^{{i\pi\over 4}\sgnlks} }

\def\Phf#1#2#3{ \Phi_{#1}(#2;#3) }
\def\PhubLpMt{ \Phf{\ube}{\cL,\cLp,M}{\ut} }
\def\PhubLzpMa{ \Phf{\ube}{\cLzL,\cLp,M}{\etpaz,\etpua} }
\def\PhubLpMa{ \Phf{\ube}{\cL,\cLp,M}{\etpua} }
\def\PhubLpMao{ \Phf{\ube}{\cL,\cLp,M}{\etpao} }
\def\PhubLzpMao{ \Phf{\ube}{\cLz,\cLp,M}{q} }
\def\PhubLmpMa{ \Phf{\ube}{\cLzL,\cLp,M}{q^\mu,\etpua} }
\def\PhubLzpMt{ \Phf{\ube}{\cLzL,\cLp,M}{t_0,\ut} }
\def\PhubLpMto{ \Phf{\ube}{\cL,\cLp,M}{t_1} }
\def\PhubLspSaz{ \Phf{\ube}{\cLzLs,\cLp,S^3}{\etpaz,\etpua, \etpuc} }

\def\PhubLspStzst{ \Phf{\ube}{\cLz,\cLp,S^3}{t_0,\ut,\ust(t_0,\ut)}}

\def\etpuc{ \etpi{\uc} }

\def\cij{ c_{ij}(\cL,\hcLs) }
\def\cijM{ c_{ij}(\cL,\hcLs,M) }
\def\LtthLs{ \lrbc{K\over 2}^{{3\over 2}\Ls} }

\newtheorem{thm}{Theorem}[section]
\newtheorem{lem}[thm]{Lemma}
\newtheorem{cor}[thm]{Corollary}
\newtheorem{rem}[thm]{Remark}

\def\Z{ \ZZ }
\def\Q{ \IQ }

\begin{document}

\begin{titlepage}
\centerline{\hfill                 math.QA/9806066}
\vfill
\begin{center}

{\large \bf
A contribution of a $U(1)$-reducible connection to quantum invariants
of links II: Links in rational homology spheres} \\

\bigskip
\centerline{L. Rozansky\footnote{
This work was supported by NSF Grant DMS-9704893}
}


\centerline{\em Department of Mathematics, Yale University}
\centerline{\em 10 Hillhouse Ave., P.O. Box 208283}
\centerline{\em New Haven, CT 06520-8283}
\centerline{\em E-mail address: rozansky@math.yale.edu}

\vfill
{\bf Abstract}

\end{center}
\begin{quotation}

We extend the definition of the $U(1)$-reducible connection
contribution to the case of the Witten-Reshetikhin-Turaev invariant
of a link in a rational homology sphere. We prove that, similarly
ot the case of a link in $S^3$, this contribution is a
formal power series in powers of $q-1$ whose coefficients are
rational functions of $q^\a$, their denominators being the powers of
the Alexander-Conway polynomial. The coefficients of the polynomials
in numerators are rational numbers, the bounds on their
denominators are established with the help of the theorem proved by
T.~Ohtsuki in Appendix~2.


We derive a surgery formula for the $U(1)$-reducible connection
contribution, which relates it to the similar contribution into the
colored Jones polynomial of a surgery link in $S^3$.

Finally, we relate the $U(1)$-reducible connection contribution to a
contribution of a $U(1)$-invariant stationary phase point to the
Reshetikhin formula for the colored Jones polynomial in the
appropriate semi-classical limit.

\end{quotation} \vfill \end{titlepage}

\pagebreak
\tableofcontents

\nsection{Introduction}
\label{s1}
\hyphenation{Re-she-ti-khin}
\hyphenation{Tu-ra-ev}

\subsection{Motivation}

In this paper we continue to study a new invariant of knots
and links defined in\cx{Ro10} which we called \emph{a
$U(1)$-reducible connection contribution} to the colored Jones
polynomial, or simply a \URC invariant. In\cx{Ro10} we studied the
properties of this invariant for links in $S^3$. Our work was based
on the $R$-matrix formula for the colored polynomial of a link.

Our task in this paper is twofold. First of all, we will show that
the \URC invariant represents a particular stationary phase
contribution to the integral in the Reshetikhin
formula\cx{Re1},\cx{Re2},\cx{Ro2} for the colored Jones polynomial.
This formula presents the expansion of the Jones polynomial in powers
of $\log q$ as an integral over the adjoint orbits of the Lie
algebra $su(2)$. The orbits correspond to the $su(2)$ modules
assigned to the link components and the integrand itself contains
formal power series in $\log q$ whose coefficients are
$SU(2)$-invariant polynomial functions. If the Jones polynomial
considered in the limit
\qq
q\rightarrow 1, \qquad q^{\a_j}=\const,
\label{1.*1}
\qqq
where $\a_j$, $1\leq j\leq L$ are the dimensions of the $su(2)$
modules assigned to $L$ components of a link,
%
then the integral seems, at least formally, to be suitable for the
stationary phase approximation. We will demonstrate that a
contribution of a stationary point corresponding to a configuration
were all the $su(2)$ integration variables belong to the same Cartan
subalgebra, is well-defined and is equal to the \URC invariant.

The polynomials in the integrand of the Reshetikhin formula
can be derived\cx{A4} from 3-valent diagrams (sometimes referred to
as `chinese characters') which come from the Kontsevich integral of a
stringed link. Therefore, the relation between the Reshetikhin
formula and the \URC invariant will allow us\cx{Ro12}
to formulate the latter solely in the language of 3-valent diagrams
and Vassiliev invariants. However here we will use this new
relation only in order to accomplish the second task of this paper -
to extend the \URC invariant to links in rational homology spheres
(\Rhsp).

The Reshetikhin formula plays an important role in
defining\cx{Ro3},\cx{Ro6} the \emph{trivial connection contribution}
to the Witten-Reshetikhin-Turaev (WRT) invariant of \Rhs and, more
generally, of links in \Rhsp. We call this contribution \emph{the
\TCC invariant}. In case of a link in $S^3$, the \TCC invariant
coincides with the colored Jones polynomial. The \TCC invariant of a
link in a general \Rhs can be presented through an analog of the
Reshetikhin formula. If a \Rhs can be constructed by a surgery on a
link $\cLs\in S^3$, then the integrand of the Reshetikhin formula
for the \TCC invariant of a link $\cL$ in that \Rhs can be expressed
by a surgery formula in terms of the Reshetikhin integrand of the link
$\cLLs\in S^3$. As a result, there exists a similar surgery formula
which relates particular stationary phase contributions to both
Reshetikhin integrals. We will use it in order to define the \URC
invariant of a link $\cL$ in a \Rhsp. The form and propeties of this
\URC invariant is very similar to those of the \URC invariant of
links in $S^3$, except that the \URC invariant of links in \Rhs is
related to the \TCC invariant of the links rather than to their WRT
invariant, as in the case of $S^3$.

Since the Reshetikhin formula and the \TCC invariant play a central
role in our paper, we will have to review their properties prior to
formulating our results.

\subsection{Notations, definitions and basic properties}

\subsubsection{Multi-index and quantum notations}
We adopt the same multi-index notations as in\cx{Ro10}.

The quantum invariants of 3-manifolds and links depend on three
parameters: $q$, $h$ and $K$ which are related
\qq
q = \exp(2\pi i/K),\qquad h = q-1.
\label{1.1*}
\qqq
We treat $q$, $h$
and $K$ as formal variables unless stated otherwise (the
WRT invariant is defined for $K$ being a positive integer).
\begin{remark}
\label{r1.1*}
\rm
Since
\qq
K^{-1} = \log(1+h)/2\pi i = h/2 \pi i + O(h^2)\quad
 \mbox{as $h\limz$},
\label{1.21}
\qqq
then $\IC[[h]] = \IC[[K^{-1}]]$ and one can easily convert a power
series in $h$ into a power series in $K^{-1}$ and back.
\end{remark}

In this paper we will
use a standard notation for a ``quantum number''
\qq
[n] = { \fdhrbas{q}{n} \over \fdh{q} }.
\label{1.1*1}
\qqq

\subsubsection{Topological notations}
First, we set the linking matrix notations. If two oriented
knots $\cK_1$,
$\cK_2$ in a 3-manifold $M$ represent trivial elements of $\HoQ{M}$,
then one can define their linking number $\lkbasm{\cK_1,\cK_2}{M}$.

Following the notations of
C.~Lescop\cx{Les}, we denote the order of $\cK$ as an
element of $\HoZ{M}$ by $o(\cK)$ (for an $L$-component link $\cL$ in
a \Rhs $M$, $\uo=\mtibas{o}{L}$ denote the orders of its link
components). Note that
\qq
\GCF{o(\cK_1),o(\cK_2)}\lkbasm{\cK_1,\cK_2}{M} \in \ZZ,
\label{1.1}
\qqq
where $\GCF{m,n}$ is the greatest common factor of $m$ and $n$.

We denote the entries of the
linking matrix of an oriented link $\cL\in M$ as
\qq
\llim{\cL}{M}_{ij} = \lkbasm{\cL_i,\cL_j}{M}.
\label{1.2}
\qqq
For a pair of oriented links $\cL,\cLp\in M$ we denote
\qq
\llim{\cL,\cLp}{M}_{ij} = \lkbasm{\cL_i,\cLp_j}{M}.
\label{1.3}
\qqq
%
We use a shortcut notation for the quadratic
form associated with the linking matrix of an oriented $L$-component
link $\cL$:
\qq
\lkbasvm{\cL}{M}{\ux} = \soinjL \llim{\cL}{M}_{ij} x_i x_j,
\label{1.4}
\qqq
where $\ux=\mtibas{x}{L}$.


Let $\cK$ be a knot in a manifold $M$. A \emph{meridian} of $\cK$ is
a simple cycle on the boundary of the tubular neighborhood of $\cK$
which is contractible through that neighborhood. The meridian is
uniquely defined. A \emph{parallel} is a simple cycle on the boundary
of the tubural neighborhood of $\cK$ which has a unit intersection
number with the meridian. A parallel is not unique (one can always
add a meridian to it). A knot is called framed if a choice of a
particular parallel has been made. We denote a framed knot as $\hcK$.
A link is framed if all of its components are framed. We denote such
link as $\hcL$.

Suppose that an oriented framed knot
$\hcK$ is a trivial element in $\HoQ{M}$. Then
a linking number between $\cK$ and its parallel is well-defined. We
call this number a \emph{self-linking number} of $\hcK$. The
self-linking numbers of the link components of $\hcL$ form a diagonal
part of the linking matrix $\llim{\hcL}{M}_{ij}$. We denote the
signature of that matrix as $\sigl{\hcL}{M}$.
For a pair of an unframed $L$-component
link $\cL$ and a framed $\Ls$-component link $\hcLs$
we define the analog of\rx{1.4}
\qq
\lkbasvm{\cLhcLs}{M}{\ux,\uy} =
\lkbasvm{\cLLs}{M}{\ux,\uy}
+ \sfrLsM y_j^2,
\label{1.5}
\qqq
where $\ux=\mtibas{x}{L}$, $\uy=\mtibas{y}{\Ls}$.

Given a quadratic form $\lkbasvm{\hcLs}{M}{\uy}$, we denote its
determinant as $\dtlks$ and its signature simply as
$\sigl{\hcLs}{M}$.

The quantum invariants of knots and links that we will consider in
this paper, that is, the colored Jones polynomial, the WRT invariant
and the contributions of various connections to those invariants,
depend on the framing. However, this dependence is known to be quite
simple. Namely, if we change the framing of a link component $\hcL_j$
in such a way that its self-linking number increases by 1, then the
corresponding invariants of $\hcL$ are multiplied by a factor
$q^{(\a_j^2-1)/4}$, where $\a_j$ is a color assigned to $\cL_j$.
Therefore a quantum invariant of a link $\cL$ in a \Rhs $M$ can be
rendered framing-independent if we multiply it by an extra factor
\qq
q^{ - \sjoL \llim{\hcL}{M}_{jj} (\a_j^2-1) /4 }.
\label{1.5*}
\qqq
Thus whenever we write a quantum invariant of an unframed knot or
link, we assume that the factor\rx{1.5*} is included in the formula
for that invariant.

If a knot $\cK$ is a trivial element of $\HoZ{M}$ (that is, if
$\oh(\cK)=1$), then it can be endowed with a canonical framing for
which the self-linking number is equal to zero.

Let $\hcLs$ be a
framed $\Ls$-component link in a 3-manifold $M$. A
new 3-manifold $M\p$ is constructed by Dehn's surgery on $\hcLs$ in
the following way: we cut out the tubular neighborhoods of the link
components and then glue them back in such a way that the meridian on
the boundary of the tubular neighborhood of each link component
matches the parallel on the corresponding boundary of the link
complement while the parallel on the boundary of the tubular
neighborhood matches the meridian on the boundary of the complement
with opposite orientation.

Suppose that the original manifold $M$ contained an oriented
link $\cL$. Then
the linking numbers of the components of $\cL$ in $M\p$ are given by
the formula
\qq
\llim{\cL}{M\p}_{ij} = \llim{\cL}{M}_{ij} - \somnLs
\llim{\cL,\cLs}{M}_{ik}
\llim{\cL,\cLs}{M}_{ik\p}
(\llim{\cLs}{M})^{-1}_{kk\p},
\label{1.6}
\qqq
where $(\llim{\cLs}{M})^{-1}$ is the inverse of the matrix
$\llim{\cLs}{M}$.

We denote the order of $\fhom$ as $\ho{M}$. For a pair of
\Rhs related by a surgery on $\hcLs\subset M$
\qq
\ho{M\p} = |\dtlks|\,\ho{M}.
\label{1.6*}
\qqq
If $M$ is a connected sum of two \Rhs, $M=M_1\# M_2$, then
\qq
\ho{M} = \ho{M_1}\,\ho{M_2}.
\label{1.6*1}
\qqq

Let $\cL$ and $\hcLs$ be a pair of oriented links in $S^3$, the
second of them being framed. Suppose that Dehn's
surgery on $\hcLs$ produces a \Rhs $M$. Then the orders of the link
components of $\cL$ as elements of $\HoZ{M}$ are
\qq
o_i = \LCM{\cij| 1\leq j\leq \Ls}
\label{1.7}
\qqq
where $\LCM{m,n}$ denotes the least common multiple of $m$ and $n$,
while
\qq
\cij = \skoLs \llim{\cL,\cLs}{S^3}_{ik}
(\llim{\cLs}{S^3})^{-1}_{kj}.
\label{1.8}
\qqq

\subsubsection{Alexander-Conway polynomial}
In subsection~1.2.4 of\cx{Ro10} we set our notations and established
some properties of the Alexander-Conway polynomial and
Alexander-Conway function of knots and links in $S^3$. Here we have
to extend these properties to knots and links in \Rhsp.

Let $\cL$ and $\hcLs$ be a pair of oriented links in $S^3$. Suppose
that a \Rhs $M$ is constructed by Dehn's surgery on $\hcLs$. Then we
define the Alexander-Conway function of $\cL$ in $M$ by a surgery
formula
\qq
&
\AFLMut =
(-1)^{\Ls}\sgnlks\,
{\AFbas{\cLLs,M}{\ut,\uust} \over \{ \fdh{(\uust)} \} },\quad
\label{1.9}
\\
&
\mbox{where $\ust_j = \pioL t_i^{-\cij}$},
\nonumber
\qqq
and $\cij$ are defined by \ex{1.8}. The Alexander-Conway function
has the following properties:
\qq
\AFbas{\cL,M}{\ut^{-1}} = (-1)^L \AFLMut
\label{1.10}
\qqq
and if we change an orientation of a component $\cL_j$, then for the
new oriented link $\cLp$
\qq
\AFbas{\cLp,M}{\ut} = - \AFbas{\cL,M}{t_1\ell t_j^{-1}\ell t_L}.
\label{1.11}
\qqq
%
The Alexander function satisfies the Torrence formula
\qq
\AFLzMuto = \
\lrbc{ \pjoL t_j^{\hlfv\llim{\cLzL}{M}_{0j}} -
\pjoL t_j^{-\hlfv\llim{\cLzL}{M}_{0j} }  }\AFLMut,
\label{1.11*1}
\qqq
where $\cL$ is an $L$-component link and $\cLz$ is a knot.
If we take a connected sum of a \Rhs $M$ containing a link $\cL$,
with another \Rhs $M\p$, then
%
\qq
\AFLMcut = \ho{M\p}\,\AFLMut.
\label{1.11*}
\qqq

If $\cK$ is a knot in $M$, then its Alexander-Conway polynomial is
defined as
\qq
\APKMt = (\fdhbas{t}{1/2o(\cK)}) \AFbas{\cK,M}{t}.
\label{1.12}
\qqq
It satisfies the property
\qq
\APbas{\cK,M}{1} = \ordHM / \oh(\cK).
\label{1.13}
\qqq

Both Alexander-Conway polynomial of a knot and Alexander-Conway
function of a link with at least 2 components are Laurent polynomials
\qq
&
\APKMt \in \ZZthopm,
\label{1.14*}
\\
&
\AFLMut \in \ZZuthopm\qquad\mbox{if $L\geq 2$}.
\label{1.14}
\qqq

\subsubsection{Lie algebra notations}
In this paper we will be dealing with the Lie algebra $su(2)$,
although it will be clear that most of our considerations can be
easily extended to other simple Lie algebras. We denote the elements
of $su(2)$ as $\vx$. The Lie algebra $su(2)$ has a Killing form. We
denote the corresponding scalar product in $su(2)$ as $\vx\cdot\vy$.
A polynomial function $\Pr\ellp{\vx_1}{\vx_n}$ on $su(2)$ is called
$SU(2)$-invariant if it is invariant under the \emph{simultaneous}
adjoint action of $SU(2)$ on all of its variables.

\subsubsection{The Reshetikhin Formula}

Let $\cL$ be an $L$-component link in $S^3$. We refer the readers to
subsection~1.2.7 of\cx{Ro10} for the definition of the colored Jones
polynomial of $\cL$. P.~Melvin and H.~Morton
considered an expansion of that polynomial in
negative powers of
\qq
K =  2\pi i / \log q.
\label{1.19}
\qqq
They found that the coefficients of this expansion are polynomials of
colors $\ual$ assigned to the components of $\cL$
\qq
\JuaL = \sumn \DumnL\,\ual^{2\um+1} K^{-n}.
\label{1.20}
\qqq
Note that
in view of Remark\rw{r1.1*}
the expansion\rx{1.20} is equivalent to that of~(1.46) of\cx{Ro10}.
%

A few years ago N.~Reshetikhin suggested\cx{Re1},\cx{Re2} to present
the expansion\rx{1.20} as an integral over adjoint orbits of the
Lie algebra $su(2)$ which correspond to $su(2)$ modules of dimensions
$\ual$. We implemented this program in\cx{Ro2}. The only shortcoming
of the arguments of\cx{Ro2} is that we based them on deriving the
Melvin-Morton expansion\rx{1.20} from a generic quantum Chern-Simons
perturbation theory, which is not a rigorous mathematical tool yet.
We easily corrected this in\cx{Ro6} by showing how to use Kontsevich
integral as a foundation of the calculations which lead to the
Reshetikhin formula.

\begin{theorem}
\label{t1.1}
For an $L$-component link $\cL\in S^3$ there exist homogeneous
$SU(2)$-invariant polynomials $\LmLal$, $\LmnLal$ such that
\qq
& \deg \LmLal = \deg \LmnLal = m,
\label{1.22}\\
& \LtLal =  \soinjL
\llim{\cL}{S^3}_{ij}
\, \val_i \cdot \val_j,
\label{1.23}\\
& \Lbas{m}{\cL}{\ual\vn} = \eval{ \del{\vual}\LmLal }{\vual=\ual\vn}
= 0,\quad \mbox{for any $\vn\in su(2)$ if $m\geq 3$},
\label{1.24}\\
&\Lbas{0,0}{\cL}{\vual} = 1
\label{1.25}
\qqq
and the Melvin-Morton expansion\rx{1.20} is equal to the
following integral of a formal power series in $K^{-1}$ and $\ual$
%
\qq
\lefteqn{
\JuaL =
\sumn \DumnL\,\ual^{2\um+1} K^{-n}
}
\label{1.26}\\
& = &
(4\pi)^{-L}
\int_{|\vual|= \ual}\lrbf{d\vual \over \ual}
\exp \lrbc{ \hpi K^{1-m} \smti \LmLal} \smnzi \LmnLal\, K^{-m-n}.
\nonumber
\qqq
\end{theorem}

\begin{remark}
\label{r1.1}
\rm
The polynomials $\LmLal$ and $\LmnLal$ are not invariants of $\cL$.
In\cx{Ro2},\cx{Ro6} we showed how these polynomials may be derived
from the Kontsevich integral of a stringed link (in\cx{A4} we will
translate these arguments into the language of Vassiliev invariants
and 3-valent graphs). We may expect that all the polynomials
$\LmLal$ and $\LmnLal$ which satisfy the conditions of
Theorem\rw{t1.1} come from the Kontsevich integral of a properly
stringed link. In any case,
from now on we will use in our formulas only the polynomials
$\LmLal$, $\LmnLal$ which come from the Kontsevich integral for a
stringed link as described in\cx{Ro2} and\cx{Ro6}.
%
\end{remark}

We will modify the integral\rx{1.26} in two ways. First, consider the
Reshetikhin formula for a link $\cLLp$ whose components carry the
colors $\ual$ and $\ube$. Let us expand the exponential into the
powers of all monomials which depend on the variables $\vube$ and
then integrate over those variables. As a result, we will get a
formula
\qq
\lefteqn{
\JuabLp =
\summpn \DummpnL\,\ual^{2\um+1}\ube^{2\ump+1} K^{-n}
}
\label{1.27}\\
& = &
(4\pi)^{-L}
\int_{|\vual|= \ual}\lrbf{d\vual \over \ual}
\exp \lrbc{ \hpi K^{1-m} \smti \LmLal} \smnzi \LmnLpalb\, K^{-m-n},
\nonumber
\qqq
where $\LmnLpab$ are homogeneous polynomials of $\vual$ and $\ube$,
homogeneous in $\vual$
%
\qq
\deg_{\vual} \LmnLpab = m,\qquad \Lbas{0,0}{\cL,\cLp}{\vual,\ube} =1.
\label{1.28}
\qqq
Now a substitution 
\qq
\ual = K \ua,\qquad \vual = K \vua
\label{1.29}
\qqq
turns \ex{1.27} into
\qq
&\hspace{0.3in}\JuabLp = \Ibas{\ube}{\cL,\cLp}{\ual/K},
\label{1.30*}
\vspace{10pt}
\\
&\IabLp = \lrbc {K\over 4\pi}^L
\iva \ivv
\etbELa \, \bPLpab,
\label{1.30}
\qqq
where $\IabLp$ is a formal power series in $\ua$, $\ube$ and $K$
\qq
\IabLp =
\summpn \DummpnL\,\ua^{2\um+1}\ube^{2\ump+1} K^{2m-n+1}.
\label{1.32}
\qqq
and $\bELa$, $\bPLpab$ denote the formal power series
\qq
\bELa = \smti \LmL,\qquad
\bPLpab = \aex.
\label{1.32*}
\qqq
%

Since for any $j$, $1\leq j\leq L\p$
\qq
\Jbas{\ual,\ube\xrem{j}}{\cL,\cLp\xrem{j}} = \eval{\JuabLp}{\b_j=1},
\label{1.33}
\qqq
then we may conclude from \ex{1.30} that
\qq
\Lbas{m,n}{\cL,\cLp\xrem{j}}{\vua,\ube\xrem{j}} =
\eval{ \LmnLpab }{\b_j=1}.
\label{1.34}
\qqq
Indeed, in view of \ex{1.33}, the polynomials
$\Lbas{m,n}{\cL,\cL\xrem{j}}{\vua,\ube\xrem{j}}$ as defined by
\ex{1.34} would suite the formula\rx{1.27} for
$\Jbas{\ube\xrem{j}}{\cL,\cLp\xrem{j}}$. However, it follows
easily from the calculations of\cx{Ro2} and\cx{Ro6} that we can make
an even stronger statement: the stringing of $\cLp\xrem{j}$
which yields
the \lhs of \ex{1.34} is induced by the stringing of $\cLp$ which
yields the \rhs of that same equation (in fact, a more careful
analysis would even show that the polynomials $\LmnLpab$ do not
depend on the stringing of $\cLp$).

\subsubsection{Formal power series and stationary phase
approximation}
\label{ss1.fs}
In this paper we will perform a lot of calculations with formal power
series. The basic operations with these series, such as addition,
multiplication and term-by-term integration used in \ex{1.27},
are well-known.
However, we will need a potentially riskier procedure
Here we will need another procedure
- a stationary phase integration of an
exponential of a formal power series. Therefore we have to explain
which manipulations with formal power series we consider well
defined.

Let $S(x;N)\in \IC[[x]]$ be a formal power series in $x$ whose
coefficients depend on an integer variable $N$
\qq
S(x;N) = \snzi s_n(N)\,x^n.
\label{1.34*}
\qqq

\begin{definition}
\label{d1.1}
The series\rx{1.34*} has a \emph{stable limit} as $N\liminf$
\qq
\slim_{N\liminf} S(x;N) = S(x),\qquad S(x) = \snzi s_n x^n,
\label{1.34*1}
\qqq
if for any $m >0$ there exists $m\p >0$ such that if $n < m$ and
$N > m\p$, then $s_n(N) = s_n$.
\end{definition}

A map $\cF:\,\cxx \rightarrow \cxx$ is
\emph{well-defined} for us if for any $S(x)\in \cxx$
\qq
\cF [S(x)] = \slimNinf
\cF[S(x|N)],\qquad\mbox{where $S(x|N) = \snzN s_n x^N$}.
\label{1.34*2}
\qqq

For two smooth functions $f(\ux)$, $g(\ux,K)$, $\ux=\mtibas{x}{L}$
$g(\ux,K)$ being
analytical in the vicinity of $K\liminf$, consider an integral
\qq
I(K) = K^{-L/2} \int e^{iK f(\ux)} g(\ux, K)\, d\ux.
\label{1.34*3}
\qqq
Suppose that $\uxst$ is a stationary point of $f(\ux)$
\qq
\eval{\del{\ux} f(\ux) }{\ux=\uxst} = 0.
\label{1.34*4}
\qqq
Then
\qq
\IKx =
K^{-L/2}
\int_{[\ux=\uxst]}
e^{iK f(\ux)} g(\ux, K)\, d\ux
\label{1.34*5}
\qqq
denotes the contribution of the point $\uxst$ to the integral $I(K)$
in the stationary phase approximation
as $K\liminf$. $\IKx$ is almost a formal power
series in $\Ki$
\qq
\IKx = e^{iKf(\uxst)}\snzi \Ix_n K^{-n},
\label{1.34*6}
\qqq
whose coefficients are determined through a special procedure.
We expand the following function in powers of $\ux-\uxst$:
\qq
\exp iK\lrbc{ f(\ux) - f(\uxst) - \ftxz } g(\ux,K) =
\sum_{\um\geq 0 \atop n \geq -|m|/3} p_{\um,n} \,\ux^{\um} K^{-n}
\label{1.34*7}
\qqq
where
\qq
\ftxz = \hlfv \sum_{i,j=1}^L \eval{\del{x_i}\del{x_j}
f(\ux)}{\ux=\uxst} (x_i-\xst_i)(x_j - \xst_j),
\label{1.34*8}
\qqq
so that
\qq
\def\umh{ {|\um| \over 2} }
\Ix_n = K^{-L/2}\sum_{\um\geq 0 \atop |\um|\leq 6n}
p_{\um,\umh-n}\,K^{\umh-n}
\int e^{iK\ftxz} \ux^{\um}\,d\ux.
\label{1.34*9}
\qqq
The gaussian integrals with polynomial prefactors can be calculated
with the help of the following formula
\qq
\intinf e^{iK A(\ux,\ux)}\, \ux^{\um}\, d\ux = (i\pi)^{L/2}
(\det A)^{-1/2}\,(2iK)^{-|\um|}\,\eval{\del{\uy}^{\um}\,
e^{-iK A^{-1}(\uy,\uy)} }{\uy=0},
\label{1.34*10}
\qqq
where $A(\ux,\ux)$ is a non-degenerate quadratic form and
$A^{-1}(\uy,\uy)$ is the quadratic form with the inverse matrix.

\subsubsection{The TCC invariant of links in rational homology
spheres}

First, let us recall the definition of the WRT invariant of links in
3-manifolds. For a link $\cL\in S^3$ and for a positive integer $K$
we define
\qq
\Zbas{\ual}{\cL,S^3} = \JuaL,
\label{1.35}
\qqq
where $q$ and $K$ are related by \ex{1.1*}. If $\cL$ is a framed link
in a \Rhs $M$, which can be constructed by Dehn's surgery on a
framed link $\hcLs\in S^3$, then
\qq
\ZuaLM = \qfrua \sugoK q^{ \qrtv\sjoLs \llst_{jj}\g_j^2 }
\Jbas{\ual,\ug}{\cLLs}
\lrbf{ \fdhrbas{q}{\ug} },
\label{1.36}
\qqq
where $\qfrua$ is a factor which removes the dependence of $\ZuaLM$
on the framing of $M$ and $\cL$
\qq
\frua
& = & - {3\over 8}\,(K-2)\sigl{\hcLs}{S^3}
-\qrtv\sjoLs\llim{\hcLs}{S^3}_{jj}
- {1\over 4} \sjoL \llim{\hcL}{M}_{jj}(\a_j^2-1)
\label{1.37}
\qqq
and $\hcL$ in the last sum is the link $\cL\in M$ endowed with the
framing induced by the 0-framing of $\cL\in S^3$ so that
according to \ex{1.6}
\qq
\llim{\hcL}{M}_{jj} =
- \somnLs
\llim{\cL,\cLs}{S^3}_{jk}
\llim{\cL,\cLs}{S^3}_{jk\p}
(\llim{\cLs}{S^3})^{-1}_{kk\p},
\label{1.38}
\qqq

The TCC invariant of a link $\cL$ in a \Rhs $M$ is defined by a
surgery formula which is similar to\rx{1.36}.
\begin{theorem}[see \cx{Ro6}]
\label{t1.2}
For an $L$-component link $\cL$ in a \Rhs $M$ there exists a unique
set of invariants $\DumnLM\in\IC$ such that the formal power series
\qq
\ZtruaLM =
\sumn \DumnLM\,\ual^{2\um+1} K^{-n},
\label{1.39}
\qqq
satisfies the following two conditions:
\begin{enumerate}
\item[(1)]
If $M=S^3$, then
\qq
\Ztrbas{\ual}{\cL,S^3} = \JuaL,
\label{1.40}
\qqq
where $\JuaL$ is represented by its Melvin-Morton expansion\rx{1.20}.

\item[(2)]
If $M$ is constructed by Dehn's surgery on a framed link
$\hcLs\in S^3$, then
\qq
\lefteqn{
\ZtruaLM = \qefrua \LtthLs
}
\label{1.41}\\
&&\times
\incs
\ivc
\bspc
\exp\hpiK\lrbc{ \sjoLs \llim{\cLs}{S^3}_{jj} \vc_j^{\,2}
+ \bELsc}\,\bPLsca,
\nonumber
\qqq
where $\incs$ denotes a contribution of the point $\vuc=0$ to the
integral over $\vuc$ in the stationary phase approximation
at $K\liminf$.
\end{enumerate}
The first of these invariants is
\qq
\Dbas{0,0}{\cL,M} = \lrbf{\ual}\hotM.
\label{1.41*}
\qqq
\end{theorem}

\rm
Since $\bELsc$ and $\bPLsca$ are formal power series, then \ex{1.41}
is just a short-cut for a more precise statement in the spirit of
\ex{1.34*2}
\qq
\lefteqn{
\ZtruaLM = \slimNinf\qefrua \LtthLs
}
\label{1.42}\\
&&\times
\incs
\ivc
\bspc
\exp\hpiK\lrbc{ \sjoLs \llim{\cLs}{S^3}_{jj} \vc_j^{\,2}
+ \bELscN}\,\bPLscaN.
\nonumber
\qqq
Here we introduced truncated series
\qq
\bELaN = \smtN \LmL,\quad
\bPLpabN = \aexN.
\label{1.43}
\qqq
%

Similarly to $\JuaL$,
the TCC invariant $\ZtruaLM$ has an integral representation.
\begin{theorem}[see\cx{Ro6}]
\label{t1.3}
For an $L$-component link $\cL$ and an $L\p$-component link $\cLp$
in a \Rhs $M$ there exist
$SU(2)$-invariant polynomials $\LmLM$, $\LmnLM$ which satisfy the
properties
\qq
& \deg \LmLMa = \deg_{\vua} \LmnLMpab = m,
\label{1.44}\\
& \LtLMa = \soiljL
\llim{\cL}{M}_{ij}
\, \va_i \cdot \va_j,
\label{1.45}\\
& \Lbas{m}{\cL,M}{\ua\vn} =
\eval{ \del{\vua}\LmLMa }{\vua=\ua\vn}
= 0,\quad \mbox{for any $\vn\in su(2)$ if $m\geq 3$},
\label{1.46}\\
&\Lbas{0,0}{\cL,M}{\vua;\ube} = \ordHM^{-3/2} \lrbf{\ube}
\label{1.47}
\qqq
and such that the formal power series
\qq
\hspace{-10pt}\bELMa = \smti \LmLM,\quad
\bPLpMab = \aexM
\label{1.48}
\qqq
combine into an integral representation
\qq
&\hspace{0.3in}\ZtruabLLpM = \Ibas{\ube}{\cL,\cLp,M}{\ual/K},
\label{1.49}
\\
&\IabLpM = \lrbc {K\over 4\pi}^L
\iva \ivv
\etbELMa \, \bPLpMab.
\label{1.50}
\qqq
\end{theorem}

There are various ways to derive the series $\bELMa$, $\bPLpMab$
from the series $\bELLsc$ and $\bPLLsca$. For our purposes here we
will use one of the formulas of\cx{Ro6}:
\qq
\lefteqn{
\etbELMa\,\bPLpMab
=
\exfruabpKx
\LtthLs
}
\nonumber\\
&&
\times
\incst
\ivc
\bspc
\exp\hpiK\lrbc{ \sjoLs \llim{\cLs}{S^3}_{jj} \vc_j^{\,2}
+ \bELLsc}\,
\label{1.51}
\\
&&\hspace{4in}\times\bPLLsca,
\nonumber
\qqq
%
where
$\vucstvua$ is the stationary point of
the exponential in the integral\rx{1.51}, which is a formal power
series in $\vua$
\qq
\vucstvua =
- \sioL \cij\, \va_i + O(\vua^2),
\label{1.52*1}
\qqq

Similarly to the case of $S^3$, the polynomials $\LmLMa, \LmnLMpab$
themselves are not invariants of $\cL$. Equation\rx{1.51}
demonstrates that they depend on the stringing of the links
$\cL,\cLp,\cLs$ in $S^3$ (a more carefull analysis\cx{A4} will show
that they depend only on the stringing on $\cL$ in $S^3$, but we
will not use this fact here).

Applying \ex{1.34} to the \rhs of \ex{1.41} we find that for any $j$,
$1\leq j\leq L$
\qq
\ZtruajLM = \eval{ \ZtruaLM}{\a_j=1}.
\label{1.52**}
\qqq
Applying \ex{1.34} to the \rhs of
\ex{1.51} we conclude that \ex{1.34} holds also
for links in $M$
\qq
\Lbas{m,n}{\cL,\cLp\xrem{j},M}{\vua,\ube\xrem{j}} =
\eval{ \LmnLMpab }{\b_j=1}.
\label{1.52*}
\qqq
if the polynomials in both sides of this equation originate
through the surgery formula\rx{1.51} from the same stringing of $\cL$
in $S^3$.

The multiplicativity of Reshetikhin's integrand under an operation of
disconnected sum of two links in $S^3$ leads to the following
\begin{proposition}
\label{p1.1}
Consider a connected sum of two \Rhs $M\# M\p$. If
$\cL,\cLp\subset M$, then
\qq
&\bEbas{\cL,M\# M\p}{\vua} = \bELMa,
\label{1.52*2}\\
&
\bPbas{\cL,\cLp,M\# M\p}{\vua;\ube} =
\Ztrbas{}{M\p}\,\bPLpMab,
\label{1.52*3}
\qqq
if we assume that the \lhs and \rhs of these equations originate from
the same stringings of links in $S^3$ through the surgery
formula\rx{1.51}.
\end{proposition}

\subsection{Results}
\label{ssres}
We prove three theorems which describe the \URC invariant of links in
\Rhs and its properties.
\begin{theorem}[\cf Main Theorem of\cx{Ro10}]
\label{t1.4}
Let $\cL$ be an oriented $L$-component link in a \Rhs $M$, whose
Alexander-Conway function is not identically equal to zero
\qq
\AFLMut \not\equiv 0,\quad
\ut = \mtibas{t}{L}.
\label{1.53}
\qqq
Let $\cLp$ be an $L\p$-component link in $M$ with positive integers
$\ube=\mtibas{\b}{L\p}$ assigned to its components. Then there exists
a unique sequence of polynomial invariants
\qq
\PubnLLpMut \in \ZZtpmtoh, \quad n\geq 0,
\label{1.54}
\qqq
such that an expansion of $\PubnLLpMut$ in powers of $(\ut-1)$ has
the form
\qq
\PubnLLpMut = \sumgz \, \pumnLLpMub\,(\ut-1)^\um,
\qquad \pumnLLpMub \in \IQub
\label{1.55}
\qqq
and the formal power series in $h$
\qq
&
\lefteqn{
\ZhrubLLpMuth
}
\nonumber
\\
& = &
\left\{
\begin{array}{cl}   \displaystyle
\hotM \snzi \PubnLp\,h^n
&\mbox{if $L=0$,}
\vspace{5pt}
\\
\displaystyle
\hi\,\hoM\,
{
1 \over \AFbas{\cL,M}{t_1}}
\snzi {\PubnLLpMt  \over \APpnLM}
\,h^n
&\mbox{if $L=1$,}
\vspace{5pt}
\\
\displaystyle
\hi\,\hoM\,{
1\over \AFLMut}
\snzi
{\PubnLLpMut \over \AFpnLM}
\,h^n
&\mbox{if $L\geq 2$}
\end{array}
\right.
\label{1.56}
\qqq
%
satisfies the following properties:
\begin{itemize}
\item[(1)]
If the link $\cL$ is empty (that is, if $L=0$), then
\qq
\Zhrem = \Ztrbas{\ube}{\cLp,M}\quad
\mbox{in $\IQubh$}.
\label{1.58}
\qqq

\item[(2)]
For any $j$, $1\leq j\leq L\p$,
\qq
\eval{ \ZhrubLLpMuth }{\b_j=1} =
\Zhrebas{\ube\xrem{j}}{\cL,\cLp\xrem{j},M}{\ut}.
\label{1.59}
\qqq

\item[(3)]
Let $\cL$ be a non-empty link. Suppose that for a number $j$,
$1\leq j\leq L$, the link component $\cL_j$ is algebraically
connected to $\cL\xrem{j}$ and
$\AFbas{\cL\xrem{j}}{\ut\xrem{j}}\not \equiv 0$. Then
\qq
\lefteqn{
\Zhrebas{\b_0,\ube}{\cL\xrem{j},\cL_j\cup\cLp}{\ut\xrem{j}}
}
\label{1.60}\\
& = &
\left\{
\begin{array}{cl}
\displaystyle
\smumu
\lrbc{ \pioLj t_i^{\llim{\cL}{M}_{ij}} }^{\mu\b_0/2}
\eval{\ZhrubLLpMuth}{t_j = q^{\mu\b_0}}
& \mbox{if $L\geq 2$,}\\
\displaystyle
\Zhrebas{\ube}{\cL,\cLp}{q^{\b_0}},
& \mbox{if $L=1$.}
\end{array}
\right.
\nonumber
\qqq
\end{itemize}

Thus defined, the polynomials $\PubnLLpMut$ have additoinal
properties:
\qq
\PubnLLpMut \in \ZZttpmtoh \quad \mbox{if $\ube$ are odd},
\label{1.61}
\qqq
\qq
&\Pebas{\ube;0}{\cL,\cLp,M}{\ut} =
\left\{
\begin{array}{cl}
\displaystyle
\lrbf{\ube} &\mbox{if $L$=0},
\vspace{10pt}
\\
\pjoLp {\piLte{\b_j} \over \piLte{1} }
& \mbox{if $L\geq 1$},
\end{array}
\right.
\label{1.62}
\qqq
\qq
&\Pebas{\ube;n}{\cL,\cLp,M}{\ut^{-1}} = \PubnLLpMut,
\label{1.63}
\qqq
and if we reverse the orientation of a link component $\cL_j$, then
for the new oriented link $\cLb$
\qq
\Pebas{\ube;n}{\cLb,\cLp,M}{\ut} =
\Pebas{\ube;n}{\cL,\cLp,M}{t_1\ell t_j^{-1} \ell t_L}.
\label{1.64}
\qqq
The polynomials $\pumnLLpMub$ contain only odd powers of $\b_j$
\qq
\pebas{\um,n}{\cL,\cLp,M}{\b_1\ell-\b_j\ell \b_L} =
-\pumnLLpMub.
\label{1.65}
\qqq
%
%
Finally, if we select an orientation for the components of $\cLp$,
then the formal power series\rx{1.56} can be rewritten as
\qq
&
\lefteqn{
\ZhrubLLpMuth
}
\label{1.66**}
\\
&&
\hspace{-10pt}
=
\left\{
\begin{array}{cl}   \displaystyle
\hotM\,
\qphlkLpMe
\snzi \PpubnLp\,h^n
&\mbox{if $L=0$,}
\vspace{5pt}
\\
\displaystyle
\hi\,\hoM\,
\qphlkLpM\,{\PhubLpMto \over \AFbas{\cL,M}{t_1}}
\snzi {\PpubnLLpMt  \over \APpnLM}
\,h^n
&\mbox{if $L=1$,}
\vspace{5pt}
\\
\displaystyle
\hi\,\hoM\,\qphlkLpM\, {\PhubLpMt\over \AFLMut}
\snzi
{\PpubnLLpMut \over \AFpnLM}
\,h^n
&\mbox{if $L\geq 2$,}
\end{array}
\right.
\nonumber
\qqq
where
\qq
\lefteqn{
\phlkLpM  =  3\lcw(M)
+ \hlfv\Bigg(
L + \soiljL \llim{\cL}{M}_{ij}
}
\label{1.66**1}
\\
&&
\hspace{1in}
+ \soiljLp \llim{\cLp}{M}_{ij}(\b_i-1)(\b_j-1)
- \siLjLp \llim{\cL,\cLp}{M}_{ij}(\b_j-1) -
\sojLp (\b_j-1)\Bigg),
\nonumber
\qqq
\qq
\PhubLpMt = \pioL t_i^{\hlfv\sjoLp \llim{\cL,\cLp}{M}_{ij}(\b_j-1)},
\label{1.66**2}
\qqq
while
\qq
\PpubnLLpMut \in \ZZtpmoh
\label{1.66**3}
\qqq
for any integer $\ube$ and
\qq
\PpubnLLpMut = \sumgz \, \ppumnLLpMub\,(\ut-1)^\um,
\qquad \ppumnLLpMub \in \IQub.
\label{1.66**4}
\qqq

\end{theorem}

\begin{remark}
\label{r1.4}
\rm
Obviously, if $M=S^3$, then this theorem is reduced to Main Theorem
of\cx{Ro10} and
\qq
\ZhrubLLpSuth = \cJrLSt.
\label{1.66*}
\qqq
%
\end{remark}

\begin{remark}
\rm
It follows from \eex{1.58} and\rx{1.60} that for a knot $\cK$ and a
link $\cLp$ in a \Rhs $M$
\qq
\Ztrbas{\a,\ube}{\cK\cup\cLp,M} = \Zhrbas{\ube}{\cK,\cLp,M;q^{\a}}.
\qqq
\end{remark}

To formulate two other theorems, we substitute
\qq
\ut = e^{\tpi \ua}
\label{1.67}
\qqq
into $\ZhrubLLpMuth$ and convert it into a formal power series
in $\Ki$ (\cf Remark\rw{r1.1*}) which we call the \URC invariant of
$\cL$ and $M$
\qq
\IrabLpM = \ehpi{ K\lkLMua }\,
\Zhrebas{\ube}{\cL,\cLp,M}{\etpi{\ua}}.
\label{1.68}
\qqq

\begin{theorem}
\label{t1.5}
Let $\cL$ be a non-empty
link in a \Rhs $M$ which satisfies\rx{1.53}. Suppose
that $M$ is constructed by Dehn's surgery on a link
$\hcLs\in S^3$ such that neither of $\Ls$ functions
\qq
c\est_j(\ua) = - \sioL \cij a_i,
\label{1.70}
\qqq
where the coefficients $\cij$ are defined by \ex{1.8}, is
identically equal to zero. Then
\qq
\lefteqn{
\IrabLpM  =
\exfruabpKy
(2K)^{\Ls/2}
}
\label{1.69}\\
&&
\qquad\times
\incrs \lrbf{d\uc}
\lrbf{ \sin(\pi\uc) }
\etpcl
\IrabcLpsS.
\nonumber
\qqq
%
\end{theorem}
Note, that
\qq
\lefteqn{\etpcl\IrabcLpsS}
\label{1.69*}
\\
&&\hspace{2.5in}=\ehpi{ K\lkLLsSac }
\Zhrebas{\ube}{\cLLs,\cLp,S^3}{\etpi{\ua},\etpi{\uc}},
\nonumber
\qqq
and $\ucstua$ is the stationary point of the phase
$\lkLLsSac$.
\begin{theorem}
\label{t1.6}
Let $\cL$ be a link in a \Rhs $M$ which satisfies\rx{1.53}. Then
\qq
\IrabLpM =
\lrbc {K\over 4\pi}^L
\ivast \ivv
\etbELMa \, \bPLpMab.
\label{1.71}
\qqq
Here $\ivast$ is a contribution of the stationary phase points of the
form $\vua = \ua\vn$, $\vn$ being a unit vector in $su(2)$, to the
integral\rx{1.50}.
\end{theorem}

Let us sketch the ideas of the proofs. We will prove the theorems in
the reverse order. First of all, we show that for the links
satisfying condition\rx{1.53}, the integral in the \rhs of \ex{1.71}
is well-defined. This allows us to regard \eex{1.71} and\rx{1.68}
as a \emph{definition} of the power series $\ZhrubLLpMuth$. We prove
that thus defined, $\ZhrubLLpMuth$ satisfies the
properties\rx{1.58}--(\ref{1.60}). Then we use the uniqueness theorem
of\cx{Ro10} to prove that $\ZhrubLLpMuth$ is the topological
invariant of $\cL,\cLp,M$ and that it coincides with the \URC
invariant defined in\cx{Ro10} in the case when $M=S^3$.

Next, we prove that the invariant $\ZhrubLLpMuth$ defined by
\eex{1.71} and\rx{1.68} satisfies the surgery formula\rx{1.69}. The
idea of the proof is to substitute the formula\rx{1.51} into
\ex{1.71}. Then the stationary phase integral over $\vuc$ can be
split into two integrals: a stationary phase integral over the
directions of $\vuc$ which together with the integral over the
directions of $\vua$ yields $\IrabcLpsS$, and a stationary phase
integral over the absolute values $|\vuc|$ which becomes the
stationary phase integral over $\uc$ in \ex{1.69}.

Finally, since we have already established that $\ZhrubLLpSuth$
defined by \eex{1.71} and\rx{1.68} coincides with the \URC invariant
of\cx{Ro10} and therefore admits the presentation\rx{1.56}, we use
the surgery formula\rx{1.69} to show that $\ZhrubLLpMuth$ also admits
that presentation for a general \Rhs $M$.

\nsection{The \URC contribution to the Reshetikhin integral}
\subsection{The definition of the \URC contribution}

The definition of the \URC invariant of a link $\cL$ depends on the
number of its components $L$.
If $L=0$, that is, if $\cL$ is an empty link, then we define
\qq
\IrabLpMe = \ZtrubLpM.
\label{2.1*}
\qqq
If $\cL$ is non-empty, then we
consider a truncated version of a general Reshetikhin
formula\rx{1.50}
\qq
\IabLpMN = \lrbc {K\over 4\pi}^L
\iva \ivv
\etbELMaN \, \bPLpMabN,
\label{2.1}
\qqq
where $N$ is a positive integer.
For all the unit vectors $\vn\in su(2)$, the points
\qq
\vua = \ua \vn, \qquad |\vn|=1
\label{2.2}
\qqq
form a 2-dimensional sphere inside the integration space
\qq
|\vua| = \ua.
\label{2.3}
\qqq
%
\begin{proposition}
\label{p2.1}
The points\rx{2.2} are stationary points of $\bELMaN$.
\end{proposition}
\proof
Since the condition\rx{2.3} implies that the linear variations of
$\va_j$ have to be orthogonal to $\va_j$, then the points\rx{2.2} are
stationary for $\LtLMa$ of\rx{1.45}. Conditions\rx{1.46} show that
these points are also stationary for all other polynomials $\LmLMa$.
\qed

We denote the stationary contribution of the sphere\rx{2.2} to
the integral\rx{2.1} in the limit of $\Kinf$ as
\qq
\IrabLpMN =
\lrbc {K\over 4\pi}^L
\ivast \ivv
\etbELMaN \, \bPLpMabN.
\label{2.4}
\qqq
Let us calculate this integral following the general procedures of
subsection\rw{ss1.fs}.

First of all, we have to reduce the
sphere\rx{2.2} to a single point. We can do this with the help of the
$SU(2)$ symmetry of the integrand. Consider the integral\rx{2.1}.
Suppose that we integrate over all the vectors except $\va_1$. The
resulting expression is a function of $\va_1$, however, due to the
$SU(2)$ symmetry, it does not depend on the direction of $\va_1$.
Therefore the remaining integral over the directions of $\va_1$ is
trivial and contributes only a factor $4\pi a_1$.
Thus we conclude that
\qq
\lefteqn{
\IabLpMN}
\label{2.5}
\\
&&\hspace{1in} =  \lrbc {K\over 4\pi}^{L-1}
Ka_1
\iva \ivvo
\etbELMaN \, \bPLpMabN.
\nonumber
\qqq
If $a_1\neq 0$, then the stationary sphere\rx{2.2} is reduced to
the stationary point
\qq
\vua\xrem{1} = \ua\xrem{1} \vns, \qquad \vns = \va_1 / a_1
\label{2.6}
\qqq
and
\qq
\lefteqn{
\IrabLpMN}
\label{2.7}
\\
&& \hspace{0.8in} = \lrbc {K\over 4\pi}^{L-1}
Ka_1
\ivasto \ivvo
\etbELMaN \, \bPLpMabN.
\nonumber
\qqq

\begin{remark}
\label{r2.1}
\rm
Fixing the direction of $\va_1$ in the integral\rx{2.5} is similar to
``breaking a strand'' in the closure of a braid in\cx{Ro10}.
\end{remark}

Now we have to consider two separate cases. If $L=1$, that is, if the
link $\cL$ is just a knot, then, in fact,
there is no integration left in
\ex{2.7}. Therefore we define
\qq
\IrabLpM = \IabLpM\qquad\mbox{if $L=1$.}
\label{2.7*}
\qqq
Furthermore, in this case $\bELMa=0$ because of \ex{1.46}, and
therefore
\qq
\IabLpM & = & Ka_1 \bPbas{\cL,\cLp,M}{a_1\vn;\ube}
\nonumber\\
& = &
Ka_1\sum_{m,\um,n\geq 0\atop |\um|\leq n}
\Dbas{m,\um;2m+n}{\cLLp,M}\,a_1^{2m}\ube^{2\um+1}\,K^{-n}.
\label{2.7*1}
\qqq
(\cf \ex{1.39}).

If $L>1$, then we proceed
to expand the exponent and
preexponential factor of the integrand in the vicinity of the
stationary point. For a fixed unit vector $\vn$
we introduce small vectors $\vux=\mtibas{x}{L}$
such that
%
$\vux \cdot \vn = 0$
and
\qq
\va_j & = & a_j \lrbc{ \cos(|\vx_j|)\,\vn +
{\sin(|\vx_j|) \over |\vx_j| }\, \vx_j}
\label{2.8}
\\
& = &
a_j\lrbs{ \vn + \vx_j +
\lrbc{ \snoi {(-1)^n\over (2n)!}\,(\vx_j\cdot\vx_j)^n} \vn +
\lrbc{ \snoi {(-1)^n\over (2n+1)!}\,(\vx_j\cdot\vx_j)^n } \vx_j }.
\nonumber
\qqq
Actually, it is convenient to identify a 2-dimensional space
orthogonal to $\vn$ with $\IC$ (\cf the complex structure on the
co-adjoint orbits), so that $\vux$ become complex variables
$\ux\in\IC$ and
\qq
\vx_i\cdot\vx_j = (x_i \xb_j + \xb_i x_j)/2.
\label{2.9}
\qqq
The integration measure for the variables $\ux$ is
\qq
{d^2\va_j \over a_j} =
{i\over 2}\,a_j\, {\sin(|\vx_j|) \over |x_j|}\,
dx_j d\xb_j
= {i\over 2}\,a_j \lrbc{
\snzi {(-1)^n\over (2n+1)!}\,(x_j\xb_j)^n } dx_j d\xb_j.
\label{2.10}
\qqq
Of course, the $SU(2)$ symmetry reduction removes the integration
over $x_1$.

Let us substitute the series\rx{2.8} into $\bELMaN$
and keep only the constant and quadratic terms in $\ux,\uxb$, while
setting $x_1=0$.
Only the term $\LtLMa$ contributes to the constant term because of
condition\rx{1.46}. Therefore
\qq
\bELMaN = \lkLMua + \AouaxN + O(\uxo^2\uxbo^2),
\label{2.11}
\qqq
where $\AouaxN$ is a bilinear form in $\ux$ and $\uxb$ whose
coefficients are polynomials of $\ua$ (the index $(1)$ indicates
that we set $x_1=0$). The terms quadratic in $\ux$ or $\uxb$ are
absent due to the $U(1)$ symmetry of the points\rx{2.2}. This $U(1)$
rotates the vectors around the axis of $\vn$ and therefore multiplies
each $x_j$ by a phase factor $e^{i\phi}$ and each $\xb_j$ by a phase
factor $e^{-i\phi}$ (note that as a result, there are no linear terms
in $\ux,\uxb$ as we expected, since the points\rx{2.2} are stationary
for the function \rx{2.11}).

Suppose that the form $\AouaxN$ is non-degenerate. Then we can
proceed with calculating the integral\rx{2.7} in the stationary phase
approximation. We keep the constant and bilinear terms in the
exponent, while expanding it in terms of higher order in $\ux,\uxb$,
so that
\qq
\lefteqn{
\etbELMaN \, \bPLpMabN
}
\label{2.12}\\
&&\hspace{1in}=
\etpiKlkMLa\,e^{\hpiKc \AouaxN}
\smnt \Bmnuax\, K^{-n},
\nonumber
\qqq
where $\Bmnuax$ are polynomials of $\ux$ and $\uxb$ of homogeneous
degree $m$, whose coefficients are polynomials of $\ua$. Then
\qq
\IrabLpMN =
(i/2\pi)^{L-1} K\lrbf{\ua}
\etpiKlkMLa
\snzi
\xFnuaN\,K^{-n}.
\label{2.13}
\qqq
Here
\qq
\lefteqn{
\xFnuaN
}
\label{2.14}\\
& = &
4^{1-L}
\szmn K^{L+m-1} \intinf \lrbf{ d\uxo d\uxbo}\,
\ehpi{K\AouaxN} \Bmnmuax\,
\nonumber
\\
& = &
{1 \over \det \AouaN} \szmn K^{m}
\eval{\Bmnmdel
\,e^{(2i/\pi K)\Aouainv}
}{\uyo=0\atop\uybo=0}
\nonumber
\\
& = &
{\SnnLMuaN \over \SdLMuaNn},
\nonumber
\qqq
where $\AouaN$ is the matrix of the quadratic form $\AouaxN$,
$\Aouainv$ is the quadratic form whose matrix is the inverse of
$\AouaN$, $\SnnLMuaN$ are polynomials in $\ua$ and we defined
\qq
\SdLMuaN =
\left\{
\begin{array}{cl}
\displaystyle
1&\mbox{if $L=1$,}\\
\displaystyle
\det \AouaN &\mbox{if $L\geq 2$}.
\end{array}
\right.
\label{2.15}
\qqq
%

\begin{remark}
\label{r2.2}
\rm
We use the notation\rx{2.15}, because it follows from the $SU(2)$
invariance of $\bELMaN$ that $\det \AnuaN{j}$ does not depend on $j$.
\end{remark}

It follows from \ex{1.47} that
\qq
\Bzuax = \Lbas{0,0}{\cL,M}{\ua\vn;\ube}
+ \Ouao = \hotM \lrbf{\ube} + \Ouao,
\label{2.15*}
\qqq
and then, according to \ex{2.14},
\qq
\SnzLMuaN = \hotM  \lrbf{\ube} + \Ouao.
\label{2.15*1}
\qqq

\begin{proposition}
\label{p2.2}
If
\qq
\SdLMuaN \not\equiv 0,
\label{2.16}
\qqq
then the stationary phase contribution\rx{2.7} is well-defined and
\qq
\IrabLpMN =
(i/2\pi)^{L-1} K\lrbf{\ua}
\etpiKlkMLa
\snzi
{\SnnLMuaN \over \SdLMuaNn}\,K^{-n}.
\label{2.17}
\qqq
The polynomials $\SdLMuaN$ and $\SnnLMuaN$ have stable limits
\qq
&\slimNinf \SdLMuaN = \SdLMua,\quad \slimNinf \SnnLMuaN = \SnnLMua,
\nonumber\\
&\SdLMua\in\IC[[\ua]],\qquad \SnnLMua\in\IC[\ube][[\ua]]
\label{2.18}
\qqq
and
\qq
\SnzLMua =  \hotM  \lrbf{\ube} + \Ouao.
\label{2.18*}
\qqq

\end{proposition}
\proof
Since we have just derived \ex{2.17}, it remains only to establish
the limits\rx{2.18}. They follow easily from the fact that a
polynomial $\LmLM$ contributes only to the terms of order
$\ua^{\um\p}$, $|\um\p|\geq m$ in $\SdLMuaN$ and $\SnnLMuaN$, while a
polynomial $\LmnLMpab$ contributes only to the terms of order
$\ua^{\um\p}$, $|\um\p|\geq m$ in $\SnnLMpuaN$, $n\p\geq n$. \qed

\begin{definition}
\label{d2.1}
Let $\cL$ be a link in a \Rhs $M$. If for a sringing of $\cL$, the
corresponding determinant\rx{2.15} is not identically equal to zero
or if $\cL$ is an empty link, then the \URC contribution to the
Reshetikhin integral\rx{1.30}, which we denote formally as
\qq
\IrabLpM =
\lrbc {K\over 4\pi}^L
\ivast \ivv
\etbELMa \, \bPLpMab,
\label{2.19*}
\qqq
is defined as
\qq
\lefteqn{\IrabLpM}
\label{2.19}\\
&& =
\left\{
\begin{array}{cl}
\displaystyle
\ZtrubLpM & \mbox{if $L=0$,}\\
\displaystyle
\IabLpM & \mbox{if $L=1$,}\\
\displaystyle
(i/2\pi)^{L-1} K\lrbf{\ua}
\etpiKlkMLa
\snzi
{\SnnLMua \over \SdLMuan}\,K^{-n}
& \mbox{if $L\geq 2$}.
\end{array}
\right.
\nonumber
\qqq

\end{definition}
%
It follows from \ex{1.52*}, that
\qq
\eval{\IrabLpM}{\ube_j=1}
= \Irbas{\b\xrem{j}}{\cL,\cLp\xrem{j},M}{\ua},
\label{2.19*1}
\qqq
if both sides of this equation originated from the same stringing of
links
in $S^3$.
Also it follows from Proposition\rw{p1.1} that
\qq
\Irbas{\ube}{\cL,\cLp,M\# M\p}{\ua} = \Ztrbas{}{M\p}\,\IrabLpM,
\label{2.19*2}
\qqq
where $M\p$ is another \Rhsp.

\subsection{Expansion of the \URC invariant}

We are going to study the properties of the expansion of the formal
power series\rx{2.19} in powers of $\ua$. First of all, we have to
establish the relation between the determinants\rx{2.15} of links and
sublinks.
%
\begin{proposition}
\label{p2.3}
For a number $j$, $1\leq j\leq L$,
if the formal power series
$\bELMa$ and $\bELjMa$ for $\cL$ and its sublink $\cLrj$
come from the same stringing of the link
$\cL$, then
\qq
\SdLMua = -a_j\,\SdLMuaj
\sioLnj\lLM_{ij} a_i
 + O(a_j^2).
\label{2.20}
\qqq
\end{proposition}
\proof
In view of Remark\rw{r2.2}, we can assume that $j\neq 1$. Consider the
matrix elements of the matrix $\AouaN$. Since the coefficients of the
powers of $\vx_i$ in the expression\rx{2.8} for $\va_i$ are
proportional to $a_i$, then it is easy to see that $\AoikLuaN$ is
proportional to $a_i a_k$ if $i\neq k$ and to $a_i$ if $i=k$.
Therefore we can introduce a matrix
\qq
\AsoikLuaN = \AoikLuaN/a_k,
\label{2.21}
\qqq
whose coefficients are also polynomials of $\ua$.
Obviously,
\qq
\SdLMua = \lrbf{\uaxo} \det \AsoLuaN.
\label{2.21*}
\qqq
We can expand $\det \AsoLuaN$ in $j$-th row. Since all the
coefficients $\AsojkLuaN$, $k\neq j$ are still proportional to $a_j$,
then
\qq
\det \AsoLuaN & = &\AsojjLuaN \det \AsjoLuaN + O(a_j)
\nonumber\\
& = &
\eval{\AsojjLuaN}{a_j=0} \det \AsjoLuaN + O(a_j),
\label{2.22}
\qqq
where $\AsjoLuaN$ is the $(j,j)$-minor of $\AsoLuaN$.

Since
$\bELMa$ and $\bELjMa$ come from the same stringing of $\cL$, then
\qq
\bELjMa = \eval{\bELMa}{a_j=0}
\label{2.23}
\qqq
and, as a result,
\qq
\AsoLjuaN = \eval{\AsjoLuaN}{a_j=0}.
\label{2.24}
\qqq
Now consider the diagonal matrix element $\AojjLuaN$. The part of it,
which is linear in $a_j$, comes from the terms in the polynomials
$\LmLMa$ which are linear in $\va_j$. However, in view of \ex{1.46},
only the polynomial $\LtLMa$ gives a non-zero contribution, so that
%
\qq
\AsojjLuaN = -  \sioLnj \lLM_{ij} a_i + O(a_j).
\label{2.25}
\qqq
Now we can conclude from \ex{2.20} written for both $\cL$ and $\cLrj$
as well as from \eex{2.22}, (\ref{2.24}) and\rx{2.25} that
\qq
\SdLMuaN = -a_j\,\SdLMuajN
\sioLnj\lLM_{ij} a_i
+ O(a_j^2).
\label{2.26}
\qqq
Equation\rx{2.20} follows if we take the stable limits of both sides
of \ex{2.26} as $\Ninf$.\qed


Let us establish a relation between the formal power series
$\IrabLpM$ and similar series coming from the sublinks of $\cL$.
\begin{proposition}
\label{p2.4}
Suppose that a link component $\cL_j$ of an $L$-component link $\cL$,
$L>1$, is algebraically linked to
$\cLrj$ (that is, there is a number $i$, $i\neq j$,
such that $\lLM_{ij}\neq 0$)
and $\SdLMuaj\not\equiv 0$. Then
\begin{itemize}
\item[(1)]
$\SdLMua\not\equiv 0$ and, moreover,
the ratios
$\lrbf{\ua}\SnnLMua/ \SdLMuan$ are non-singular at $a_j=0$;
\item[(2)]
if we expand ${\SnnLMua / \SdLMuan}$ in powers of $a_j$ and
substitute $a_j=\a_j/K$, then
\qq
\smumu \eval{\IrabLpM}{a_j = \mu\a_j/K} =
\IrabLpMj
\label{2.27}
\qqq
\end{itemize}
\end{proposition}
\proof
The fact that $\SdLMua\not\equiv 0$ follows easily from \ex{2.20}.
Thus the stationary phase integral
$\IrabLpMN$ is well-defined for sufficiently large $N$. Let us
study its expansion in powers of $a_j$. Instead
of expanding the series\rx{2.13} directly, we will work with the
integral\rx{2.7}. Let us assume without the loss of generality that
$j\neq 1$. We can split the integral\rx{2.7} into two integrals
\qq
\lefteqn{
\IrabLpMN = \lrbc{K\over 4\pi}^{L-2}
}
\label{2.28}
\\
&&
\hspace{0.5in}\times Ka_1
\ivastoj \ivvoj
\etbELMajN \,\bPrjLpMabN,
\nonumber
\qqq
\qq
\bPrjLpMabN
& = & {K\over 4\pi} \ivastj
{d\va_j \over a_j}\,
\etbA
\label{2.29}
\\
&&\hspace{2in}\times
\bPsjLpMabN,
\nonumber
\qqq
\qq
\vajst = {\bALj \over |\bALj|}\,a_j.
\label{2.30}
\qqq
We derived this formula by splitting the exponent $\bELMaN$ of
\ex{2.7} in 3 pieces: the terms which do not depend on $\va_j$
(that is, $\bELMajN$), the terms which are linear in $\va_j$
(that is, $\bALj\cdot\va_j$)
and the terms of higher order in $\va_j$. We
expanded the exponential $\etbELMaN$ in the third group of terms,
combining this expansion with preexponential series $\bPLpMabN$ and
thus obtaining $\bPsjLpMabN$:
\qq
&\bPsjLpMabN = \smnt \bPsjmn\,K^{-n},
\label{2.30*}
\qqq
where $\bPsjmn$ are polynomials of $\vua$ and $\ube$ of homogeneous
degree $m$ in $\va_j$.
Then we fixed the vectors $\vuarj$ near
the stationary point\rx{2.6} and performed the stationary phase
integral\rx{2.29} over $\va_j$. Note that
\qq
\bALj = \sioLnj \lLM_{ij}\va_i + O(\vuarj^2).
\label{2.31}
\qqq
Since $\cL_j$ is algebraically connected to $\cLrj$, then
\qq
\sioLnj \lLM_{ij}\va_i \not \equiv 0,\qquad
|\bALj| \not\equiv 0,
\label{2.32}
\qqq
and the stationary point\rx{2.30} is well-defined. The
integral\rx{2.29} is calculated on a term-by-term basis with the help
of the formula
\qq
\ivastj
{d\va_j \over a_j}\,
\etbA\,P(\va_j)
=
-\tpi\,P(-i\del{\vb})\,
\eval{ e^{ia_j |\vb|} \over |\vb|}{\vb=\pi K\bALj},
\label{2.33}
\qqq
which follows from a simple relation
\qq
\ivaj
{d\va_j \over a_j}\,
\etbA\,P(\va_j)
=
-\tpi\,P(-i\del{\vb})\,
\eval{ e^{ia_j |\vb|} - e^{-ia_j |\vb|}
\over |\vb|}{\vb=\pi K\bALj},
\label{2.34}
\qqq
where $P(\va_j)$ is a polynomial of $\va_j$. As a result,
\qq
\lefteqn{
\bPrjLpMabN
}
\label{2.35}
\\
&&
\hspace{1in}
 = e^{i\pi K a_j |\bALj|} \skmn
{\bPjLpMabNkmn \over |\bALj|^{k+1}}\,a_j^m\,K^{-(n+1)},
\nonumber
\qqq
where
\qq
\bPjLpMabNkmn \in \IC[\vuarj,\ube].
\label{2.36}
\qqq
To complete the calculation of the expansion of $\IrabLpMN$ in powers
of $a_j$, it remains to calculate the integral\rx{2.28}
with\rx{2.35} on a term-by-term basis. The integration is performed
in the same way as in\rx{2.4}, the only difference being the presence
of denominators in the preexponential factor\rx{2.35}. However, it
follows from \ex{1.46} that
\qq
\eval{\bALj}{\vuarj=\uarj\vn} =
\vn\sioLnj \lLM_{ij}\a_i,
\label{2.37}
\qqq
and, as a result,
\qq
\eval{|\bALj|}{\vuarj=\uarj\vn} =
\sioLnj \lLM_{ij}a_i  \not\equiv 0.
\label{2.38}
\qqq
Therefore the denominators are non-zero in the vicinity of the
stationary phase point $\vuarj=\uarj\vn$, so we can expand them in
powers of $\vux\xrem{1,j}$ after the substitution\rx{2.8} without
encountering singularities. Thus the the integration
in\rx{2.28} will yield the following expression
\qq
\lefteqn{\IrabLpMN}
\label{2.39}
\\
&& = K\lrbf{\uarj}\etpiKlkMLa
\sklmn{\Sjklmn\,a_j^m \over \lrbc{\SdLMuaNl}^{l+1}
\lrbc{\sioLnj \lLM_{ij}a_i }^{k+1}
}\,K^{-n},
\nonumber
\qqq
\qq
\Sjklmn\in\IC[\uarj,\ube]
\nonumber
\qqq
which represents the expansion of $\IrabLpMN$ in powers of $a_j$.
Note that we did not touch $a_j$ in the exponential $\etpiKlkMLa$.
It is clear from \ex{2.39} that the ratios
$\lrbf{\ua}{\SnnLMuaN / \SdLMuaNn}$ are non-singular at $a_j=0$. The
first claim of Proposition\rw{p2.4} follows from this by taking the
stable limit of those ratios at $\Ninf$.

To prove the second claim of Proposition\rw{p2.4}, we
recall the relation between the integrands for $\IrabLpM$ and
$\IrabLpMj$. Namely, the integrand for the stationary phase
Reshetikhin integral of $\IrabLpMj$ can be obtained by
expanding the exponential of the integrand for $\IrabLpM$ in terms
which contain the powers of $\va_j$ and then by performing the
integration over $\va_j$, $|\va_j| = \a_j/K$. This means that
\qq
\IrabLpMNj = \eval{\IrabLpMNt}{a_j=\a_j/K},
\label{2.40}
\qqq
\qq
\lefteqn{
\IrabLpMNt = \lrbc{K\over 4\pi}^{L-2}
}
\label{2.41}
\\
&&
\hspace{0.5in}\times Ka_1
\ivastoj \ivvoj
\etbELMajN \,\bPjLpMabN,
\nonumber
\qqq
\qq
\bPjLpMabN
& = & {K\over 4\pi} \ivaj
{d\va_j \over a_j}\,
\etbA
\label{2.42}
\\
&&\hspace{2in}\times
\bPsjLpMabN.
\nonumber
\qqq
The latter integral can be calculated on a term-by-term basis with
the help of \ex{2.34} (note that we do not have to expand $\etbA$ in
powers of $\va_j$ in the \lhs of this equation, because its \rhs is
non-singular at $a_j=0$). Comparing \eex{2.33} and\rx{2.34}, we
conclude that
\qq
\bPjLpMabN = \smumu\, \bPjLpMmabN.
\label{2.43}
\qqq
The similarity between \eex{2.28} and\rx{2.41} implies that
\qq
\IrabLpMNt = \smumu \,\IrmabLpMN
\label{2.44}
\qqq
and \ex{2.27} follows by applying the stable limit $\Ninf$ in view of
\ex{2.40}. \qed

\subsection{Uniqueness arguments}
\label{s2.3}
We will use the uniqueness theorem of\cx{Ro10} in order to prove
that the formal power series $\IrabLpM$ as defined by \ex{2.19}, is
the invariant of the links $\cL,\cLp$, that is, it does not depend on
their stringing. First of all, we prove that $\IrabLpM$ is
well-defined for algebraically connected links.

\begin{proposition}
\label{c2.1}
If $\cL$ is algebraically connected, then for \emph{any} stringing
\qq
\SdLMua\not\equiv 0.
\label{2.26*}
\qqq
\end{proposition}
\proof
We prove this proposition
by induction in the number of link components $L$. If $L=1$,
then, obviously, $\SdLMua = 1$. Suppose that \ex{2.26*} holds for all
$L-1$-component algebraically connected links. Let $\cL$ be an
$L$-component algebraically connected link. Then it is easy to see
that there exists a link component $\cL_j$ such that $\cLrj$ is also
algebraically connected. According to the assumption of the
induction, this means that $\SdLMuaj\not\equiv 0$. Since $\cL_j$ is
algebraically connected to $\cLrj$, then
\qq
\sioLnj\lLM_{ij} a_i \not\equiv 0
\label{2.26*1}
\qqq
and \ex{2.26*} follows from
\ex{2.26}.\qed

Let $\CL$ denote a set of stringings of a link $\cL$ is a \Rhs $M$
such that for those stringings $\SdLMua\not\equiv 0$. Let $\bcLs$
denote the set of all links for which the set $\CL$ is non-empty. For
each stringing $c\in\CL$ we cancel all common powers of $\ua$
in the ratios $$\lrbf{\ua}{\SnnLMua / \SdLMuan}$$ and denote the
result as $\Fcnuab/\lrbc{\Gcua}^{2n+1}$. Finally, let $\e=\Ki$. Now
we are ready to apply the uniqueness theorem of\cx{Ro10}.

\begin{proposition}
\label{p2.5}
The formal power series $\IrabLpM$ is a topological invariant of an
oriented link $\cL$ and a link $\cLp$ in a \Rhs $M$, that is,
$\IrabLpM$ does not depend on a stringing $c\in\CL$ used in order to
define the integrand of\rx{2.19*}.
\end{proposition}
\proof
This result follows from Theorem\row{3.1}, which is applicable in
view of \ex{2.7*} (coupled with \ex{1.49}),
\ex{2.19*1}
and
Proposition\rw{p2.4}. \qed

\begin{proposition}
\label{p2.6}
Let $\cL$ be a link in $S^3$ such that
\qq
\SdLSua\not\equiv 0,\qquad \AFLSut \not\equiv 0.
\label{2.45}
\qqq
Then
\qq
\IrabLpS = \cJrLSa.
\label{2.46}
\qqq
where $\cJrLSt$ is defined by eq.\ro{1.52}.
\end{proposition}
\proof
This proposition follows from the relation
\qq
\IrabLpSe = \cJrLSe = \JubLp
\label{2.47}
\qqq
and Theorem\row{3.1}.\qed

Comparing the leading terms of the formal power series in both sides
of \ex{2.46} we come to the following
\begin{corollary}
\label{c2.2}
If a link $\cL$ satisfies the conditions\rx{2.45}, then
%
\qq
\lrbf{\ua} {\SnzLSua\over \SdLSua}
=\aldf\;
{\Pebas{\ube;0}{\cL,\cLp,S^3}{\ut} \over \AFLSut},
\label{2.47*}
\qqq
where $\Pebas{\ube;0}{\cL,\cLp,S}{\ut}$ is given by \ex{1.62}.
\end{corollary}
\begin{proposition}
\label{p2.7}
If $\cL\in\bcLs$ then $\SdLMua\not\equiv 0$ for \emph{any} stringing
of $\cL$.
\end{proposition}
\proof
Let us add an extra component $\cLz$ to $\cL$ in such a way that
\qq
\lLM_{0j} = 1\qquad\mbox{for all $j$, $1\leq j\leq L$}.
\label{2.48}
\qqq
Select a stringing of $\cLzL$ which is compatible with the stringing
of $\cL$. Then, according to \ex{2.20},
\qq
\eval{a_0^{-1}\SdLzMua}{a_0=0}
= - \SdLMua \sjoL a_j
\label{2.49}
\qqq
This means that, as functions of $\ua$,
\qq
\SdLMua\equiv 0\qquad
\mbox{iff}\qquad \eval{a_0^{-1}\SdLzMua}{a_0=0}\equiv 0.
\label{2.50}
\qqq
On the other hand, in view of \ex{2.18*},
\qq
\eval{\SdLzMua \over a_0\, \SnzzLMua}{a_0=0}\equiv 0\qquad
\mbox{iff}\qquad
\eval{a_0^{-1}\SdLzMua}{a_0=0}\equiv 0.
\label{2.51}
\qqq
As a result,
\qq
\SdLMua\equiv 0\qquad
\mbox{iff}\qquad
\eval{\SdLzMua \over a_0\, \SnzzLMua}{a_0=0}\equiv 0.
\label{2.52}
\qqq

Obviously, the link $\cLzL$ is algebraically connected, so in view
of Proposition\rw{c2.1}, it belongs to $\bcLs$. Therefore the
ratio $\SdLzMua/ (a_0\,\SnzzLMua)$, as a part of $\IrabLzpM$, is
the topological invariant of the link $\cLzL$, that is, it does not
depend on its stringing. Then\rx{2.52} proves
Proposition\rw{p2.7}.\qed
\begin{proposition}
\label{p2.8}
For the links in $S^3$,
$\SdLSua\not\equiv 0$
iff $\AFLSut\not\equiv 0$.
\end{proposition}
\proof
Throughout the proof we assume that $M=S^3$ (we want to use $M$,
because later we will use the same proof in order to prove this
proposition for any \Rhs $M$). The proof is similar to the previous
one. Let $\cL$ be an $L$-component link in $M$. We supplement it
with an extra component $\cLz$ satisfying \ex{2.48} and thus come
to\rx{2.52}. On the other hand, \ex{2.47*} together with \eex{1.62}
and\rx{2.18*} imply that
\qq
\eval{a_0^{-1}\SdLzMua}{a_0=0}\equiv 0
\qquad\mbox{iff}\qquad
\AFLSut\equiv 0.
\label{2.53}
\qqq
A combination of\rx{2.52} and\rx{2.53} proves the proposition.\qed

A combination of Propositions\rw{p2.6} and\rw{p2.8} leads to the
following
\begin{proposition}
\label{p2.9*}
If a link $\cL\subset S^3$ satisfies a condition
$\AFLSut\not\equiv 0$, then the stationary phase integral\rx{2.19*}
is well-defined and
\qq
\IrabLpS = \cJrLSa,
\label{2.53*1}
\qqq
where $\cJrLSt$ is defined by eq.\ro{1.52}.
%
\end{proposition}


\nsection{Surgery formula and rational expressions}
\subsection{A surgery formula for the \URC invariant}
So far, we have defined the \URC invariant $\IrabLpM$ as a
particular stationary phase contribution to the Reshetikhin integral
and established its relation to the \URC invariant of links in $S^3$,
as defined in\cx{Ro10}. Our next step is to prove a surgery formula
which would express the \URC invariant of links in rational homology
spheres in terms of the \URC invariant of links in $S^3$.
\begin{proposition}
\label{p2.9}
Let $\cL$ be a link in a \Rhs $M$ such that
\qq
\SdLMua\not\equiv 0, \qquad \AFLMut\not\equiv 0.
\label{2.53*}
\qqq
Suppose that $M$ can be constructed by Dehn's
surgery on a framed link $\hcLs\in S^3$ such that
\qq
\ucst(\ua)\not\equiv 0,
\label{2.54}
\qqq
where the functions $\ucst(\ua)$ are defined by \ex{1.70}. Then
for any link $\cLp\in M$, the \URC invariant of $\cL,\cLp$ is given
by the surgery formula\rx{1.69}.
\end{proposition}
\proof
First of all, we note that
\qq
\AFLst\not\equiv 0,
\label{2.54*}
\qqq
because, according to \ex{1.9},
\qq
\AFLst = \sgnlks
\lrbf{\fdhbas{e}{i\pi\ucst(\ua)} }
\AFLMua.
\label{2.54*1}
\qqq
Now let us substitute the truncated version of the
stationary phase surgery formula\rx{1.51} into
the stationary phase integral\rx{2.7}
and combine both integrals into one
\qq
\lefteqn{
\IrabLpMNN
}
\label{2.55}
\\
& = & \lrbc{K\over 4\pi}^{L-1} \LtthLs Ka_1
\,\exfruabpKy
\nonumber\\
&&\qquad\times
\ivasto \ivvo
\incst \ivc \bspc
\nonumber
\\
&&\qquad\qquad\times
\exp\hpiK\lrbc{ \sjoLs \llim{\cLs}{S^3}_{jj} \vc_j^{\,2}
+ \bELLscN}\,\bPLLscaN
\nonumber\\
& = & \lrbc{K\over 4\pi}^{L-1} \LtthLs Ka_1
\,\exfruabpKy
\nonumber\\
&&\qquad\times
\ivastog \ivvo
\ivc \bspc
\nonumber
\\
&&\qquad\qquad\times
\exp\hpiK\lrbc{ \sjoLs \llim{\cLs}{S^3}_{jj} \vc_j^{\,2}
+ \bELLscN}\,\bPLLscaN,
\nonumber
\qqq
where in view of \eex{1.46}
\qq
\vucst(\ua\vns) =
- \sioL \cij a_i\vns.
\label{2.56}
\qqq
Note that, generally speaking, $\IrabLpMNN\neq\IrabLpMN$, but it is
easy to see that $\IrabLpMNN$, as defined by \ex{2.55}, is expressed
in the same form as\rx{2.17} and it has the same stable limit
\qq
\slimNinf \IrabLpMNN = \IrabLpM
\label{2.56*}
\qqq
in the sense of \ex{2.18}. We leave the details to the reader.

We can split the integral over $\vuc$ into the integral over the
orbits $|\vuc|=\uc$ and the integral over the radii $\uc$ of the
orbits
\qq
\lefteqn{ \IrabLpMNN =
\exfruabpKy
 (2K)^{\Ls/2}
\lrbc{K\over 4\pi}^{L+\Ls-1}Ka_1
}
\label{2.57}
\\
&&\times
\incrs
\lrbf{d\uc}
\lrbf{ \sin(\pi\uc) }
\etpcl
\nonumber\\
&&\qquad
\times
\icomb
\hspace{-0.15pt}
\ivvo \ivvc e^{\hpiKc \bELLsSacN}
\bPLLsScaN.
\nonumber
\qqq
We can rearrange the integral in this way, because
\eex{2.54} and\rx{2.54*}
guarantee in view of Proposition\rw{p2.8} that
the inner integral of\rx{2.57} is well-defined. In fact, according
to \ex{2.5} the inner integral (together with the appropriate
normalization factor) yields $\IrabcLpsSN$, and so we come to the
following surgery formula for the truncated invariants
\qq
\lefteqn{
\IrabLpMNN  =
\exfruabpKy
(2K)^{\Ls/2}
}
\label{2.58}\\
&&
\qquad
\times\ \incrs \lrbf{d\uc}
\lrbf{ \sin(\pi\uc) }
\etpcl\IrabcLpsSNh.
\nonumber
\qqq
We recover \ex{1.69} by taking the stable limits of both
sides of this equation at $\Ninf$.\qed

Now we want to prove that two conditions\rx{2.53*} are equivalent
to each other. This will bring us directly to Theorem\rw{t1.5}. The
proof requires a few steps.

\begin{proposition}
\label{p2.10}
Under the conditions of Proposition\rw{p2.9}
\qq
\lrbf{\ua} {\SnzLMua\over \SdLMua}
= \aldf\,\hoM\,
{\PeLpMua \over \AFLMua},
\label{2.59}
\qqq
where $\PeLpMt$ is defined by \ex{1.62}.
\end{proposition}
\proof
Compare the leading contribution to the stationary phase integral of
\ex{1.69}
\qq
\lefteqn{
\IrabLpM
}
\label{2.60}
\\
&& = {K\over\tpi}\,
\exp\hpiK\lrbc{\lkLLsSacst -
\sjoL \llim{\hcL}{M}_{jj}a_j^2 }
\nonumber
\\
&&
\qquad\times
\Bigg(
{\sin(\pi\ucst(\ua))\over \AFLst}\,
e^{ -{3\over 4}\,i\pi\sigl{\hcLs}{S^3} }
(2K)^{\Ls/2}
\nonumber
\\
&&\hspace{3in}\times\intinf d\uc\, \ehpi{K \lkLsS}
\,+\,O(\Ki)
\Bigg).
\nonumber
\qqq
%
%
with \ex{2.19}. Then relation\rx{2.59} follows from the following
equations:
\qq
\lkLLsSacst - \sjoL \llim{\hcL}{M}_{jj}a_j^2
= \lkLMua
\label{2.61}
\qqq
(which follows from \ex{1.6}),
\qq
{\sin(\pi\ucst(\ua))\over \AFLst}
=
{\sigl{\hcLs}{S^3}\, (i/2)^{L_s} \over \AFLMua}
\label{2.62}
\qqq
(which follows from \ex{1.9}),
\qq
&\intinf d\uc\, \ehpi{K \lkLsS} = (K/2)^{-\hlfv\Ls}e^{{i\pi\over 4}
\sigl{\hcLs}{S^3}} |\dtlks|^{-\hlfv},
\label{2.63}\\
&e^{-{i\pi\over 2}\sgnlks} i^{\Ls} = \sgnlks,
\label{2.64}
\qqq
and
\qq
|\dtlks| = \ho{M}
\label{2.65}
\qqq
(which is a particular case of \ex{1.6*}).\qed

\begin{proposition}
\label{p2.11}
If \ex{2.59} holds for a \Rhs $M\# M\p$ such that
$\cL,\cLp\subset M$, then it also holds for $M$.
\end{proposition}
\proof
This proposition follows easily from \ex{2.59} for $M\# M\p$ and from
\ex{2.19*2}, from \rx{1.41*} applied to $\Ztrbas{}{M\p}$ and from
\eex{1.11*} and\rx{1.6*1}.\qed

\begin{proposition}
\label{p2.12}
Let $\cL\subset M$ be a link
satisfying the conditions\rx{2.53*}.
If there exists a
knot $\cLz\subset M$ such that it is algebraically connected to $\cL$
and \ex{2.59} holds for a link $\cLzL$, then this equation also holds
for the original link $\cL$.
\end{proposition}
\proof
Let us compare the expressions for both sides of \ex{2.59} in case of
a link $\cL$ and in case of a link $\cLzL$ at $a_0=0$. First of all,
a combination of \eex{2.27} and\rx{2.19*1} indicates that
\qq
\IrabLpM = \smumu\, \Irbas{\ube}{\cLzL,M}{\mu/K,\ua}.
\label{2.66}
\qqq
Comparing the leading terms on both sides of this equation we find
that
\qq
\lrbf{\ua} {\SnzLMua\over \SdLMua} & = &
\eval{
a_0\lrbf{\ua}
{\Snz{\cLzL,\cLp,M;a_0,\ua;\ube} \over \Sdv{\cLzL,M;a_0,\ua} }
}{a_0=0}
\label{2.67}
\\
&&\qquad\times
\lrbc{
\fdhbas{e}{i\pi\sjoL\llim{\cLzL}{M}_{0j}a_j}
}.
\nonumber
\qqq
On the other hand, since  according to \ex{1.62},
\qq
\Pebas{\ube;0}{\cLzL,\cLp,M}{1,\ut} =
\Pebas{\ube;0}{\cL,\cLp,M}{\ut},
\label{2.68}
\qqq
and in view of
\ex{1.11*1}
\qq
\lefteqn{{\PeLpMua \over \AFLMua}}
\label{2.69}
\\
&&
\hspace{1in}=
\lrbc{\fdhbas{e}{i\pi\sjoL\llim{\cLzL}{M}_{0j}a_j}}
{\Pebas{\ube;0}{\cLzL,\cLp,M}{1,\etpua}\over \AFLzMuao}.
\nonumber
\qqq
Equations\rx{2.67} and\rx{2.69} demonstrate that if we set $a_0=0$
in \ex{2.59} written for the $\cLzL$, then we get the same
equation for $\cL$.\qed

\begin{proposition}
\label{p2.13}
Equation\rx{2.59} holds for any link $\cL$ in a \Rhs $M$
which satisfies the conditions\rx{2.53*}.
\end{proposition}
\proof
According to a slightly generalized theorem of\cx{murakami}, for any
link $\cL$ in a \Rhs $M$ there exists an algebraically split link
$\hcLs\subset S^3$ such that Dehn's surgery on $\hcLs$ produces a
\Rhs $M\# M\p$ with $\cL\subset M$. Let us pick a knot
$\cLz\subset S^3$ which has the following linking numbers with $\cL$
and $\cLz$:
\qq
\llim{\cLzL}{S^3}_{0i} = \delta_{1i} + \sjoLs
\llim{\cLzL,\cLs}{S^3}_{ij},\qquad
\llim{\cLzL,\cLs}{S^3}_{0j} = \llim{\hcLs}{S^3}_{jj}.
\label{2.70}
\qqq
Then, according to \ex{1.6},
\qq
\llim{\cLzL}{M}_{0i} = \delta_{1i}.
\label{2.71}
\qqq
Thus $\cLz$ is algebraically connected to $\cL$ in $M$ and therefore
the link $\cLzL$ satisfies the conditions\rx{2.53*}. It also
satisfies the condition\rx{2.54}, because according to
\eex{2.70},\rx{1.70} and\rx{1.8},
\qq
c\est_j(a_0,\ua) = -a_0 + c\est_j(\ua),\qquad
e\est_j(\ua) = - {1\over \llim{\hcLs}{S^3}_{jj} }
\sioL\llim{\cL,\cLs}{S^3}_{ij}\,a_i.
\label{2.72}
\qqq
Therefore, according to Proposition\rw{p2.10}, \ex{2.59} holds for
$\cLzL$ and $M\# M\p$. Then, according to Proposition\rw{p2.11},
\ex{2.59} holds also for $\cLzL$ and $M$, and therefore, according to
Proposition\rx{p2.12}, it holds for $\cL$ and $M$.\qed

Now, repeating the same steps that led us from Corollary\rw{c2.2} to
Proposition\rw{p2.8}, we come to the following
\begin{proposition}
\label{p2.14}
For a link $\cL$ in a \Rhs $M$, $\SdLMua\not\equiv 0$ iff
$\AFLMut\not\equiv 0$.
\end{proposition}

A combination of this proposition with Propositions\rw{p2.2}
and\rw{p2.9} sums up to the following
\begin{proposition}
\label{p2.15}
If a link $\cL$ in a \Rhs $M$ satisfies the condition
$\AFLMut\not\equiv 0$, then for any stringing of $\cL$ the
integral\rx{2.19*} is well-defined and it is a topological invariant
of $\cL$. If $M$ can be constructed by Dehn's surgery on a link
$\hcLs\in S^3$ which satisfies the condition\rx{2.54}, then the
surgery formula\rx{1.69} holds.
\end{proposition}

Now it remains to show that the ratios $\SdLMua/\SnnLMua$ of
\ex{2.19} are rational functions of $\etpua$, as suggested by
\eex{1.56} and\rx{1.68}.


\subsection{Rational structure and integrality}

\subsubsection{Surgery on algebraically split links}

In this subsection we will complete the proof of our results by
proving the following
\begin{proposition}
\label{p3.1}
Let $\cL$ be an oriented link in a \Rhs $M$ such that
$\AFLMut\not\equiv 0$. Then for any link $\cLp\subset M$ there exist
the polynomials\rx{1.66**3}
%
%
such that they satisfy\rx{1.66**4}
and the formal power series $\ZhrubLLpMuth$ defined by \ex{1.66**} is
related to
the modified \URC invariant
\qq
\IhrabLpM & = & e^{-\hlfv i\pi K\lkLMua}\, \IrabLpM
\label{3.2}
\qqq
as in \ex{1.68}:
\qq
\Zhrebas{\ube}{\cL,\cLp,M}{\etpi{\ua}} = \IhrabLpM.
\label{3.3}
\qqq
\end{proposition}
A proof of this proposition will require us to prove a few lemmas. We
begin with a modification of Lemma~2.3 of\cx{murakami} which was
suggested by T.~Ohtsuki.
\begin{lemma}[H.~Murakami, T.~Ohtsuki]
\label{l3.1}
For a \Rhs $M$ and a prime number $K$ which does not divide $\ho{M}$,
there exists an algebraically split link $\hcLs\subset S^3$, such
that Dehn's surgery on $\hcLs$ produces a \Rhs $M\# \MK$, where
\qq
\MK = \Lens{p_1,1}\#\cdots \#\Lens{p_N,1}.
\label{3.1}
\qqq
Here $\Lens{p_j,1}$ are lens spaces
and the numbers $\up=\mtibas{p}{N}$ are not divisible by $K$.
\end{lemma}
We can use the proof of\cx{murakami} if we observe that for even
values of $\ho{M}$ the factors $E^k_0$ and $E^k_1$ do not appear in
the decomposition of the linking pairing of $M$.\qed

Now we can prove the following
\begin{lemma}
\label{l3.2}
If Proposition\rw{p3.1} holds for any \Rhs which can be constructed
by Dehn's surgery on an algebraically split link in $S^3$, then it
holds for any \Rhsp.
\end{lemma}
\proof
Let $M$ be a \Rhsp. According to our assumption, Proposition\rw{p3.1}
holds for the \Rhs $M\# \MK$ described in the previous lemma.
It also holds for a lens space $\Lens{p,1}$, because this manifold
can be constructed by Dehn's surgery on an unknot in $S^3$ with
self-linking number $p$. Since $\ho{\Lens{p,1}}=|p|$, then according
to \ex{1.66**} (case of $L=0$),
\qq
&\Ztrbas{}{\Lens{p,1}} = |p|^{-3/2}\,q^{3\lcw(\Lens{p,1})}
\snzi \Ppbas{n}{\Lens{p,1}}\,h^n,
\label{3.3*2}\\
&\Ppbas{n}{\Lens{p,1}}
\in \ZZ[1/p],\qquad\Ppbas{0}{\Lens{p,1}} =1
\nonumber
\qqq
and
\qq
{1\over \Ztrbas{}{\Lens{p,1}}} = |p|^{-3/2} q^{-3\lcw(\Lens{p,1})}
\snoi \Ptbas{n}{\Lens{p,1}}\,h^n,\qquad\Ptbas{n}{\Lens{p,1}}
\in \ZZ[1/p].
\label{3.3*1}
\qqq
According to \ex{2.19*2},
\qq
\IhrabLpM =
{\Zhrebas{\ube}{\cL,\cLp,M\# \MK}{\etpi{\ua}} \over
\pjoN
\Ztrbas{}{\Lens{p_j,1}} }.
\label{3.3*}
\qqq
Since
\qq
\ho{M\# \MK} = \ho{M}\pjoN |p_j|,\qquad
\lcw(M\#\MK) = \lcw(M) + \sjoN \lcw(\Lens{p_j,1}),
\label{3.3*3}
\qqq
then it follows from \eex{3.3*} and \rx{3.3*1} that $\IhrabLpM$ can
be presented in the \rx{1.66**}
with $\ut=\etpua$ if we extend the
ring\rx{1.66**3} by $1/\lrbf{\up}$.
Since $\lrbf{\up}$ is not divisible by $K$, then the intersection of
all the extended rings for all prime $K$ which do not divide $\ho{M}$
is the original ring. This proves the lemma.\qed

\noindent
\emph{Another proof of Lemma\rw{l3.2}}
According to Corollary~A2.2 of T.~Ohtsuki's Lemma~A2.1 of Appendix,
for any \Rhs $M$ there exists an algebraically split link
$\hcLs\subset S^3$ such that a surgery on $\hcLs$ produces a \Rhs
$M\# M\p$, where
\qq
M\p = \Lens{p_1,q_1}\#\cdots\#\Lens{p_N,q_N},
\label{3.4}
\qqq
and $1/\lrbf{\up}\in \ZZihoM$. A lens space $\Lens{p_j,q_j}$ can be
constructed by a rational surgery on an unknot in $S^3$, therefore it
satisfies Proposition\rw{p3.1} in view of\cx{Ro8} (note that for an
empty lens space
\ex{2.1*} sets $\Zhrbas{}{M} = \ZtrM$). We can also see
this explicitly, because there is an analytic expression for
$\Ztrbas{}{\Lens{p,q}}$ which can be found in\cx{Je} (see
also\cx{Ro8} and references therein)
\qq
\Ztrbas{}{\Lens{p,q}} =
\prosign(p)\,|p|^{-1/2} \,q^{3s(p,q)}\,
{\fdhbas{q}{1/(2p)}\over \fdh{q}},
\label{3.5}
\qqq
where $s(p,q)$ is the Dedekind sum.
According to \ex{2.19*2},
\qq
\IhrabLpM =
{\Zhrebas{\ube}{\cL,\cLp,M\# M\p}{\etpi{\ua}} \over
\pjoN
\Zhrbas{}{\Lens{p_j,q_j}} }.
\label{3.6}
\qqq
and the claim of the lemma follows from this expansion if we perform
a division of formal power series in its \rhs in the same way as
we did in \ex{3.3*}.\qed

\subsubsection{Integrality properties of the Alexander polynomial and
its expansion}

Before we proceed with the proof of Proposition\rw{p3.1}, we have to
strengthen the integrality properties\rx{1.14*} and\rx{1.14}. For a
link $\cL$ in a \Rhs $M$ we define
\qq
\nu_i(\cL,M) = \sjoLni \llim{\cL}{M}_{ij} + 1.
\label{1.17*}
\qqq
Suppose that $M$ can be constructed by a surgery on an algebraically
split link $\hcLs\subset S^3$. Then we define
\qq
\nup_i(\cL,\hcLs) = \sjoLs
{\llim{\cL,\cLs}{S^3}_{ij}
(p_j - \llim{\cL,\cLs}{S^3}_{ij})
\over p_j },
\label{1.17*1}
\qqq
where $p_j$ denote the self-linking numbers of the link components of
$\hcLs$
\qq
p_j = \llim{\hcLs}{S^3}_{jj}.
\label{1.17*2}
\qqq
We also assume that if $M=S^3$ (that is, if $\hcLs$ is an empty
link), then $\nup(\cL,\empt) = 0$. Now we define
\qq
\PhLsut = \ut^{-\hlfv\lrbc{
\unu(\cL,M) + \unup(\cL,\hcLs)
} }.
\label{1.15}
\qqq
The factor $\PhLsut$ absorbs the factors $1/2$ in the powers of $t$
and $\ut$ in\rx{1.14*} and\rx{1.14}
\qq
&t^{1/2}\,\Phbas{\cK,\hcLs}{t}\APKMt \in \ZZtopm,
\label{1.16}\\
&\PhLsut\AFLMut \in \ZZutopm\qquad
\mbox{if $L\geq 2$}
\label{1.17}
\qqq
If $M=S^3$, then these relations can be found in\cx{Tu}. Otherwise
they follow easily from \ex{1.9}.

\begin{lemma}
\label{l3.*1}
If $\cL$ is a link in a \Rhs
$M$ is constructed by a surgery on an algebraically split link
$\hcLs\subset S^3$, then
\qq
\sjoLs\hlfv\,{p_j+1\over p_j},\;  \hlfv\,\unup(\cL,\hcLs) \in
\ZZihoM.
\label{3.6*}
\qqq
\end{lemma}
\proof
Since $\ho{M}=\lrbf{\up}$, then $\ZZbas{1/p_j}\subset \ZZihoM$ and it
is enough to prove that
\qq
\hlfv\,{p_j + 1\over p_j},\;
\hlfv\,
{\llim{\cL,\cLs}{S^3}_{ij}
(p_j - \llim{\cL,\cLs}{S^3}_{ij})
\over p_j }
\in \ZZbas{1/p_j}.
\label{3.6*1}
\qqq
Consider two cases. If $p_j$ is even, then
$\ZZbas{1/p_j}=\ZZbas{1/2,1/p_j}$ and\rx{3.6*1} is obvious. If $p_j$
is odd, then $p_j+1$ and
$\llim{\cL,\cLs}{S^3}_{ij}(p_j - \llim{\cL,\cLs}{S^3}_{ij})$ are even
and\rx{3.6*1} follows from there.\qed

\begin{lemma}
\label{l3.*2}
Let $\cL$, $\cLp$ be a pair of links in a \Rhs $M$ constructed by a
surgery on an algebraically split link $\hcLs\subset S^3$. Then we
can expand the negative odd powers of the Alexander-Conway function
around $\uuv=\uusr$ in the following way
\qq
\lefteqn{
{1\over \AFpbas{2n+1}{\cLLp,M}{\ut,\uuv}}
}
\label{3.10*}
\\
&&
\hspace{1in}
= \pjoLp \lrbc{u_j\over u_j^*}^{\hlfv\lrbc{
\nu_j(\cLLp,M) + \nup_j(\cLp,\hcLs)}}
\slmn {w_{l,\um,n}(\ut,\uusr)\,(\uuv/\uusr-1)^{\um} \over
\AFpbas{2n+2l+1}{\cLLp,M}{\ut,\uusr} },
\nonumber\\
&&\hspace{1.5in}
w_{l,\um,n}(\ut,\uusr)\in \ZZutus.
\nonumber
\qqq
Also if $\cK$ is a knot in $M$, then
\qq
&{1\over \APpbas{2n+1}{\cK,M}{t} }
= t^{\hlfv\nup(\cK,\hcLs)} \smzi w_m\,(t-1)^m,
\qquad
w_m\in\ZZihoM.
\label{3.10**1}
\qqq
\end{lemma}
\proof
Let us first prove \ex{3.10*}.
Consider the expansion
\qq
\lefteqn{
\Phbas{\cLLp,\hcLs}{\ut,\uuv}\,\AFbas{\cLLp,M}{\ut,\uuv}
}
\label{3.6*2}\\
&&
\hspace{1.5in}
= \Phbas{\cLLp,\hcLs}{\ut,\uusr}\,\AFbas{\cLLp,M}{\ut,\uusr}
+ w(\ut,\uusr,\uuv/\uusr-1),
\nonumber\\
&&
\hspace{1in}
w(\ut,\uusr,\us)\in
\ZZutus[[s]],
\qquad w(\ut,\us) = O(\us).
\nonumber
\qqq
Then
\qq
{1\over \AFpbas{2n+1}{\cLLp,M}{\ut,\uuv}}
&=&{\Phpbas{2n+1}{\cLLp,\hcLs}{\ut,\uuv} \over
\lrbc{\Phbas{\cLLp,\hcLs}{\ut,\uuv}\,\AFbas{\cLLp,M}{\ut,\uuv}
}^{2n+1} }
\label{3.6*3}
\\
& = & {\Phbas{\cLLp,\hcLs}{\ut,\uuv}\over
\Phbas{\cLLp,\hcLs}{\ut,\uusr}}
\slmn {w_{l,\um,n}(\ut,\uusr)\,(\uuv/\uusr-1)^{\um} \over
\AFpbas{2n+2l+1}{\cLLp,M}{\ut,\uusr} }
\nonumber\\
&&w_{l,\um,n}(\ut,\uusr)\in
\ZZutus.
\nonumber
\qqq
We kept the power in denominators odd by multiplying both deniminator
and numerator by an extra factor of $\Phi\nabla_{\rm A}$ if needed.
Also, we absorbed even powers of $\Phi$ in the polynomials $w$.
Equation\rx{3.10*} follows easily from \ex{3.6*3}.

The proof of \ex{3.10**1} is similar to that of \ex{3.10*}, except we
have to use \ex{1.13} which determines the value of denominators.
We leave the details to the reader.\qed

\subsubsection{Extra link component}

\begin{lemma}
\label{l3.3}
Let $M$ be a \Rhs which is constructed by a surgery on a link
$\hcLs\subset S^3$. Suppose that a link $\cL\subset M$, $L\geq 1$
 satisfies the
condition $\AFLMut\not\equiv 0$. If we add a new component $\cLz$ to
$\cL$ in such a way that its linking numbers with the components of
$\cLs$ and $\cL$ in $S^3$ are given by \ex{2.70}, then the claims of
Proposition\rw{p3.1} hold for $\cLzL$.
\end{lemma}
\proof
As we showed in the proof of Proposition\rw{p2.15}, the
\URC invariant $\IrabLzpM$ can be expressed by the surgery
formula\rx{1.69} where the link $\cL$ is substituted by $\cLzL$. On
the other hand, Proposition\rw{p2.9*} establishes that
Proposition\rw{p3.1} holds for links in $S^3$, which allows us to
present $\IhrabLzpM$ in the form\rx{1.66**}. As a result
\qq
\IhrabLpMz = \snzi \yFnqaz\,h^{n-1}
\label{3.7}
\qqq
where
\qq
\yFnqaz & = &
q^{\fruabpKyz+
\phlkLpsS
}\,e^{-\hlfv i\pi K\lkLMuaz}\,
i^{-\Ls}
\label{3.8}\\
&&\quad\times(K/2)^{\Ls/2}
\incrs \lrbf{d\uc}
\lrbf{ e^{\pi i\uc} - e^{-\pi i\uc} }
\ehpi{ K\lkLLsSacz }
\nonumber
\\
&&\hspace{1.5in}\times
\PhubLspSaz
\nonumber\\
&&\hspace{2in}\times
{\PpubnLLpSuaz \over \AFpnoLSaz}.
\nonumber
\qqq
Now we have to change the integration variable into
\qq
\ux = \uc - \ucst(\ua),
\label{3.9}
\qqq
expand the preexponential factor in\rx{3.8} in powers of
$\etpi{\ux}-1$ and integrate each term in this expansion
individually. We use \ex{3.10*} in which we substitute $\cL$ by
$\cLzL$, $\cLp$ by $\cLs$, $\ut$ by $\etpaz, \etpua$ and $\uusr$ by
$\etpi{\uust(a_0,\ua)}$.
Then the expansion of the preexponential factor in \ex{3.8}
takes the form
\qq
\lefteqn{
\lrbf{ e^{i\pi \uc} - e^{-i\pi \uc} }
\PhubLspSaz}
\nonumber\\
&&\hspace{2in}
\times{\PubnLLpMuaz \over \AFpnoLMaz}
\label{3.11}\\
&&
\lefteqn{
=
\PhubLzpMa
\e^{i\pi\sjoLs \kapp_j x_j}
}
\nonumber\\
&&\hspace{1.5in}\slmn{\Qlmnaz\; (\etpi{\ux}-1)^{\um}
\over \lrbc{ \AFLzMua
\lrbf{\etpi{\ucst(a_0,\ua)}-1}
}^{2n+2l +1}},
\nonumber
\qqq
where
\qq
&\kapp_j = - \sizL \llim{\cLzL,\cLs}{S^3}_{ij}
+ \sioLp \llim{\cLp,\cLs}{S^3}_{ij}(\b_i - 1),
\label{3.12}\\
&
\qquad
\Qlmntz \in \ZZutiz.
\label{3.12*}
\qqq
and we used the fact that
\qq
\PhubLspStzst = \PhubLzpMt.
\label{3.12*1}
\qqq
The substitution\rx{3.9} turns the exponent of \ex{3.8} into
\qq
\lkLLsSacz = \lkLzMa +
\sjoLs \yp_j x_j^2
\label{3.13}
\qqq
where
\qq
\llim{\hcLzL}{M}_{ii} = -\sjoLs{(\llim{\cLzL,\cLs}{S^3}_{ij})^2
\over
\yp_j
},\qquad 0\leq i\leq L
\label{3.14}
\qqq
and the numbers $\up$ denote the diagonal linking numbers of
the algebraically split surgery link $\hcLs$ according to
\ex{1.17*2}.
Now the gaussian integrals are easy to calculate: for any
$\Ls$ integers $\um$
\qq
\lefteqn{
(K/2)^{\Ls/2}\intinf \lrbf{d\ux} \exp \hpi
\sjoLs \lrbc{ Kp_j\,x_j^2 + 2(\kapp_j + 2m_j)x_j}
}
\label{3.16}\\
&&\hspace{2.5in}=
\eqpsgns\, |\lrbf{\up}|^{-1/2}
\,
q^{-\sjoLs{1\over \yp_j}\lrbc{ m_j^2 + m_j\kapp_j + \qrtv(\kapp_j)^2}},
\nonumber
\qqq
and therefore
\qq
\lefteqn{
(K/2)^{\Ls/2}\intinf\lrbf{d\ux} e^{\hpiKc\sjoLs p_jx_j^2}\,
e^{i\pi\sjoLs\kapp_j x_j}\,\lrbf{\etpi{\ux}-1}^{\um}
}
\label{3.17}\\
&&\hspace{2.5in}=
\eqpsgns\,\hoM\,q^{-\qrtv\sjoLs {(\kapp_j)^2\over \yp_j}}
\sum_{n\geq \hlfv|\um|} C_n\,h^n,
\nonumber\\
&&\hspace{1in}\qquad C_n\in\ZZihoM.
\nonumber
\qqq
Assembling all phase factors from \eex{3.8},\rx{3.13} and\rx{3.17},
we can present their product as
\qq
\lefteqn{
q^{\fruabpKyz+
\phlkLpsS
}\,e^{-\hlfv i\pi K\lkLMuaz}\,
i^{-\Ls}
\ehpi{K\lkLzMa}
}
\label{3.18}\\
&&\hspace{4in}
\times
\eqpsgns\,q^{-\qrtv\sjoLs {(\kapp_j)^2\over \yp_j}}
\nonumber\\
&& = (-1)^{\Ls}\sgnlks\,
\qphlkLpzM\,
q^{{3\over 4}\sgnlks-\qrtv\sjoLs \lrbc{p_j-{2\over p_j}} -3\lcw(M)}
q^{\hlfv\sjoLs {p_j+1\over p_j}}
\nonumber   \\
&&
\hspace{3.5in}
\times q^{\hlfv\lrbc{\sizL\nup_i(\cL,\hcLs) -
\sioLp\nup_i(\cLp,\hcLs)\, (\b_i-1)} }
\nonumber
\qqq
It follows from the surgery formula for $\lcw(M)$ that
\qq
{3\over 4}\sgnlks-\qrtv\sjoLs \lrbc{p_j-{2\over p_j}} -3\lcw(M)
\in \ZZihoM.
\label{3.20}
\qqq
Combining this with Lemma\rw{l3.*1} we conclude that
%
\qq
\lefteqn{
(1+h)^{{3\over 4}\sgnlks-\qrtv\sjoLs \lrbc{p_j-{2\over p_j}}
-3\lcw(M)}\,
(1+h)^{\hlfv\sjoLs {p_j+1\over p_j}}
}
\label{3.21}
\\
&&\hspace{1.5in}
\times
(1+h)^{\hlfv\lrbc{\sizL\nup_i(\cL,\hcLs) -
\sioLp\nup_i(\cLp,\hcLs)\, (\b_i-1)} }
\in\ZZbas{1/\ho{M},\ube}[[h]].
\nonumber
\qqq
%
Expanding the expression\rx{3.21} in powers of $h$ as it appears in
the product of phase factors\rx{3.18}, and combining
this expansion with the integrals\rx{3.8} calculated with the help of
\eex{3.11} and\rx{3.17}, we come to the following formula
\qq
\lefteqn{\IhrabLpMz
= \hi\, \hoM\,
\qphlkLpzM\,
\PhubLzpMa
}
\nonumber
\\
&&
\hspace{1in}
\times
\snzi
{\Qnaz\over
\lrbc{ \AFLzMua\lrbf{\etpi{\ucst(a_0,\ua)}-1}}^{2n+1}}\,h^n,
\label{3.23}
\qqq
where
\qq
&\Qntz\in
\ZZtpmohz\quad\mbox{if $\ube\in\ZZ_+$}.
\label{3.24}
\qqq
Thus the lemma is proved if we show that the factors
$\lrbf{\etpi{\ucst(a_0,\ua)}-1}$ can be canceled in the denominators
of the \rhs of \ex{3.23}.


Let us rewrite \ex{2.19} for the link $\cLzL$ as a power series in
$h$, while substituting $\SdLzMua$ with its expression in terms
of $\AFLzMua$ which can be derived from \ex{2.59}
\qq
\lefteqn{
\IhrabLpMz
}
\label{3.25}\\
&&=
\hi \snzi {\SpnnLMuaz \over
\lrbc{
\lrbf{\ua} \SnzLMuaz{}\,\AFLzMua
}^{2n+1}}\,h^n,
\nonumber
\qqq
where
\qq
\SpnnLMuaz \in \IQ[\ube][[a_0,\ua]].
\label{3.26}
\qqq
Note that we dropped the factor of $a_0^{2n+1}$ in the denominators,
because, according to claim~1 of Proposition\rw{p2.4}, each term in
the sum\rx{3.26} is non-singular at $a_0=0$, and therefore the
powers of $a_0$ can be canceled against the numerators.

Compare the series\rx{3.23} and\rx{3.25}.
Obviously, all $\Ls$ factors $\etpi{\ucstaz}-1$ are divisible by
the corresponding functions $\ucstaz$, but according to \ex{2.18*},
$\SnzLMuaz{}$ is not divisible by either of $\ucst(a_0,\ua)$.
Therefore $\Qnaz$ is divisible by $\lrbf{\ucstaz}^{2n+1}$. Since the
expressions
\qq
\uust(t_0,\ut) -1 = t_0^{-1}\,\uust(\ut) -1
\label{3.27}
\qqq
do not factor over $\ZZtpmohz$, then the divisibility of $\Qnaz$
means that $\Qntz$ is divisible by $\lrbf{\uust(t_0,\ut)-1}^{2n+1}$:
\qq
&\Qntz = \lrbf{\uust(t_0,\ut)-1}^{2n+1}
\PpubnLLpMutz,
\label{3.28}\\
&\PpubnLLpMutz\in\ZZtpmoh.
\nonumber
\qqq
Substituting this expression into \ex{3.23} we get
\qq
\lefteqn{\IhrabLpMz
}
\label{3.29}\\
&&
\hspace{1.2in}
= \hi\, \hoM\,
\qphlkLpzM\,
\PhubLzpMa
\nonumber
\\
&&
\hspace{2.3in}
\times
\snzi
{\PpubnLLpMuaz\over \AFLnozMua},
\nonumber
\qqq
which is equivalent to \ex{1.66**} for $\cLzL$.

It remains to check \ex{1.66**4} for $\cLzL,\cLp$. We substitute the
expansion\rx{1.66**4} for the polynomials $\PpubnLLpSuaz$ in \ex{3.8}
and integrate term by term. Then an easy power counting indicates that
\qq
\Qntz = \summzgz  \qmumn\,(t_0-1)^{m_0}(\ut-1)^\um,
\qquad \qmumn \in \IQub.
\label{3.29*1}
\qqq
According to \ex{3.28}, this series is divisible by
$\lrbf{\uust(t_0,\ut)-1}^{2n+1}$. Thus performing this division at
the level of formal power series in $t_0-1,\ut-1$ we come
to\rx{1.66**4} for $\cLzL,\cLp$.\qed


\begin{lemma}
\label{l3.5}
Proposition\rw{p3.1} holds for a link $\cL$ in a \Rhs $M$ if
$L\geq 2$.
\end{lemma}
\proof According to Lemma\rw{l3.2}, we can assume that $M$ is
constructed by a surgery on an algebraically split link
$\hcLs\subset S^3$. Then, according to Lemma\rw{l3.3},
Proposition\rw{p3.1} holds for the link $\cLzL$, where a knot $\cL$
is desribed by \ex{2.70}. Thus we start with \ex{3.29} and substitute
it into \ex{2.66} which in view of \ex{2.71}
can be rewritten for the invariants $\Ih$ as
\qq
\IhrabLpM = \smumu\,
e^{i\pi\mu a_1}
\Ihrbas{\ube}{\cLzL,M}{\mu/K,\ua}.
\label{3.30}
\qqq

In view of \eex{3.10*} and\rx{1.11*1},
\qq
\lefteqn{
{1\over \AFLnqmMua}
}
\label{3.31}
\\
&&
\hspace{1in}
= q^{\hlfv\mu\lrbc{1+\lzo+\nup_0(\cL,\hcLs)}}
\slmmn {w_{l,m,n}(\etpua)\,(q^\mu -1)^m
\over
\bigg(\AFLMua (\fdhbas{e}{i\pi\lzo a_1})\bigg)^{2n+2l+1} },
\nonumber
\\
&&\hspace{1.5in}w_{l,m,n}(\ut)\in\ZZihoMuti,
\nonumber
\qqq
where we used a notation $\lzo=\llimzo$. Although we know from
\ex{2.71} that $\lzo=1$,
still we keep it explicitly in \ex{3.31}, because we will use
the same expression in the next lemma when that linking number has a
different value.

It is easy to see that
\qq
\lefteqn{
\qphlkLpzM\,
\PhubLmpMa
}
\label{3.32}\\
&&
 =
q^{\hlfv\lrbc{1 + \lzo +
(\mu-1)\sjoLp\llim{\cLzL,\cLp}{M}_{0j}(\b_j-1)} }
\,e^{-i\pi \lzo a_1}
\,\qphlkLpM\,\PhubLpMa.
\nonumber
\qqq
Since in view of the second relation of\rx{3.6*}
\qq
q^{\hlfv\mu\nup_0(\cL,\hcLs) + {\mu+1\over 2}(\lzo+1)
+{\mu-1\over 2}\sjoLp\llim{\cLzL,\cLp}{M}_{0j}(\b_j-1)
}
%
\in\ZZihoMh,
\label{3.33}
\qqq
then we conclude from \eex{3.29} and\rx{3.30} that there exist
polynomials $$\Qnut\in\ZZtpmoh$$ such that
\qq
\IhrabLpM
& = & \hi\, \hoM\,
\qphlkLpM\,
\PhubLpMa
\nonumber
\\
&&
\qquad
\times
\snzi
{\Qnua\over
\bigg( \AFLMua(\etpao-1)\bigg)^{2n+1}}\,h^n
\label{3.34}
\qqq
(here we used the fact that $\lzo=1$).
Thus we have to prove that $\Qnut$ is divisible by $(t_1-1)^{2n+1}$.
To see this consider the same calculation that led us to \ex{3.34}
but with a different knot $\cLz$. Since $L\geq 2$, we can choose
$\cLz$ which satisfies conditions\rx{2.70} with $\delta_{2i}$ instead
of $\delta_{1i}$ in the expression for $\llim{\cLzL}{S^3}_{0i}$. As a
result, we will get a similar formula for $\IhrabLpM$
\qq
\IhrabLpM
& = & \hi\, \hoM\,
\qphlkLpM\,
\PhubLpMa
\nonumber
\\
&&
\qquad
\times
\snzi
{\Qnuap\over
\bigg( \AFLMua(\etpat-1)\bigg)^{2n+1}}\,h^n,
\label{3.35}
\qqq
with some other polynomials
$\Qnutp\in\ZZtpmoh$.
Comparing \eex{3.34} and (\ref{3.35}) we find that
\qq
(t_2-1)^{2n+1}\Qnut = (t_1-1)^{2n+1}\Qnutp.
\label{3.36}
\qqq
Since the factors $(t_1-1)$ and $(t_2-1)$ are coprime, this means
that $\Qnut$ is indeed divisible by $(t_1-1)^{2n+1}$ so that
\qq
\Qnut = (t_1-1)^{2n+1}\,\PpubnLLpMut,\quad
\PpubnLLpMut\in\ZZtpmoh,
\label{3.36*1}
\qqq
and \ex{3.34}
leads to \ex{1.66**}.

Relation\rx{1.66**4} for $\cL,\cLp$ follows
easily from the same relation for $\cLzL$, $\cLp$. Indeed, if we
substitute \ex{1.66**4} for $\cL,\cLp$ into the \rhs of \ex{3.30},
then it follows easily that
\qq
\Qnut = \sumgz \qumn\,(\ut-1)^{\um},\qquad\qumn\in\IQub.
\label{3.36*}
\qqq
Now \ex{1.66**4} for $\cL,\cLp$ emerges if we divide the formal power
series\rx{3.36*} by $(t_1-1)^{2n+1}$.\qed

\begin{lemma}
\label{l3.6}
Proposition\rw{p3.1} holds for a link $\cL$ in a \Rhs $M$ if
$L=1$.
\end{lemma}
\proof
The proof is similar to that of the previous lemma.
Let $\cL$ be a one-component link in a \Rhs $M$. We add to it an
extra component $\cLz$ which is the meridian of $\cL$, so that
\qq
\llim{\cLzL}{M}_{01} = \lzo = 1/o_1.
\label{3.37}
\qqq
In view of Lemma\rw{l3.5}, Proposition\rw{p3.1} holds for the
link $\cLzL$. Therefore we repeat the same calculations that led us
from \eex{3.29} and\rx{3.30} to \ex{3.34}. This time however we can
use \ex{1.12} for the Alexander-Conway polynomial of the knot $\cL$.
As a result, instead of \ex{3.34} we end up with the formula
\qq
\IhraobLpM
& = & \hi\, \hoM\,
\qphlkLpM\,
\PhubLpMao\,e^{-i\pi a_1/o_1}
\nonumber
\\
&&
\qquad
\times
\snzi
{\Qnuao\over
\APpnoLMa}\,h^n,
\label{3.38}\\
&&
\Qnuto\in\ZZtoohom.
\nonumber
\qqq

Since $L=1$, then according to \ex{2.19}
\qq
\IhraobLpM = \IraobLpM = \IaobLpM.
\label{3.39}
\qqq
The latter invariant is given by a formal power series\rx{2.7*1}.
Since this series is proportional to $a_1$, then
\qq
\Ihrbas{\ube}{\cL,\cLp,M}{0} = 0.
\label{3.40}
\qqq
In view of \ex{1.13}, this implies that
\qq
\Qnuto = (t_1^{1/o_1}-1)\,\PpubnLLpMt
\label{3.41}
\qqq
with $\PpubnLLpMt \in \ZZtopmoh$. Substituting this relation into
\ex{3.38} we come to \ex{1.66**}.

Relation\rx{1.66**4} for $\cL$ follows easily from the same relation
for $\cLzL$ in the same way as in the previous lemma.\qed

\begin{lemma}
\label{l3.7}
Proposition\rw{p3.1} holds for an emply link $\cL$, that is, if
$L=0$.
\end{lemma}
\proof
We place a knot $\cLz$ into $M$ in such a way that $\cLz$ represents
a trivial element in $\HoZ{M}$. According to the previous lemma,
Proposition\rw{p3.1} holds for $\cLz$. On the other hand, according
to \eex{2.7*},\rx{1.49},\rx{1.52**} and\rx{2.1*}
\qq
\IhrabLpMe = \IrabLpMe = \Ihrasub.
\label{3.42}
\qqq
Thus
\qq
\IhrabLpMe & = & \hi\, \hoM\,
\qphlkLzpM\,\PhubLzpMao
\,(\fdh{q})
\label{3.43}
\\
&&\qquad\times
\snzi {\PpubnLzLpMq\over \APpLzq}\,h^n.
\nonumber
\qqq
Now we have to expand the instances of $q=1+h$ in the \rhs of this
equation in powers of $h$. We do this calculation similarly to that
of Lemma\rw{l3.5}, except that we use \ex{3.10**1} instead of
\ex{3.10*}. According to \ex{3.10**1},
\qq
{1\over \APpLzq} = q^{\hlfv \nup(\cLz,\hcLs)} \smzi w_m\,h^m,\qquad
w_m\in\ZZihoM.
\label{3.44}
\qqq
Also it is easy to see that
\qq
\hi\,(\fdh{q})\,\qphlkLzpM\,\PhubLzpMao
=
\qphlkLpMe.
\label{3.45}
\qqq
Substituting \eex{3.44} and\rx{3.45} into \ex{3.43} we obtain an
expression which is in the \rhs of \ex{1.66**} for $L=0$.\qed

Three lemmas\rw{l3.5},\rw{l3.6} and\rw{l3.7} prove
Proposition\rw{p3.1}

\subsection{Final proofs}

Now we can prove the theorems of subsection\rw{ssres}.

\pr{Theorem}{t1.4}
First of all, let us prove the existence of polynomials
$\PubnLLpMut$. Consider the formal power series $\ZhrubLLpMuth$ as
defined by \ex{3.3}. Proposition\rw{p3.1} says that it can be
presented in the form\rx{1.66**},\rx{1.66**3}.
If we combine
$\PhubLpMt$ with the polynomials $\PpubnLLpMut$, expand
the factor $\qphlkLpM$ in powers of $h=q-1$ and write the resulting
expression as a single power series in $h$, then, obviously, we come
to \ex{1.56} with the polynomials $\PubnLLpMut$ satisfying\rx{1.54}
and\rx{1.55}. The property\rx{1.61} follows from the fact that
%
$\PhubLpMt \in \ZZbas{\ut,\ut^{-1}}$ if $\ube$ are odd.
%

In view of the definitions\rx{3.3} and\rx{3.2}, \ex{1.58} follows from
\ex{2.1}, \ex{1.59} follows from \ex{2.19*1} while \ex{1.60} follows
from \ex{2.27}.

The uniqueness of the polynomials $\PubnLLpMut$ and their
properties\rx{1.63}--(\ref{1.65}) follow from Theorem~3.1 of\cx{Ro10}
in exactly the same way as in the case of $M=S^3$, which was
considered in that paper.

The formula\rx{1.62} follows from Propositions\rw{p2.13}
and\rw{p2.14}.\qed

Since we adopted \ex{3.3} as the definition of $\ZhrubLLpMuth$, then
Theorem\rw{t1.6} is obvious and the claim of
Theorem\rw{t1.5} is contained in Proposition\rw{p2.15}.
Finally, note that
Theorem\rw{t1.5} can be extended to the case when a surgery is
performed on a link in a general \Rhsp, rather than in $S^3$.
\begin{theorem}
Let $\cL$ be a non-empty
link in a \Rhs $M\p$
such that
\qq
\AFLMput \not\equiv 0.
\qqq
Suppose
that a \Rhs $M\p$ is constructed by Dehn's surgery on a link
$\hcLs\in M$ such that neither of $\Ls$ functions
\qq
c\est_j(\ua) = - \sioL \cijM a_i,
\quad\mbox{\rm where}\;\;
\cijM = \skoLs \llim{\cL,\cLs}{M}_{ik}
(\llim{\cLs}{M})^{-1}_{kj}
\qqq
is
identically equal to zero. Then
\qq
\lefteqn{
\IrabLpMp  =
\exfruabpKyM
(2K)^{\Ls /2}
}
\\
&&
\qquad\times
\incrs \lrbf{d\uc}
\lrbf{ \sin(\pi\uc) }
\etpclM
\IrabcLpsM,
\nonumber
\qqq
%
where the framing correction $\fruabpKyM$ is defined by the formula
similar to \ex{1.37}
\qq
\hspace{-4pt}\fruaM
 =  - {3\over 8}\,(K-2)\sigl{\hcLs}{M}
-\qrtv\sjoLs\llim{\hcLs}{M}_{jj}
- {1\over 4} \sjoL \llim{\hcL}{M\p}_{jj}(\a_j^2-1).
\qqq
\end{theorem}

\proof
A case of a surgery on a link in a \Rhs is similar to that of a
surgery on a link in $S^3$, because an obvious analog of the
formula\rx{1.51} holds if we substitute $S^3$ by $M$ and $M$ by
$M\p$. The proof of Proposition\rw{p2.9*} also stays the same if we
perform this substitution. \qed


\section*{Appendix}
\appendix
\nsection{Construction of a rational homology sphere by a surgery
on an algebraically split link}

\begin{center}

{\sc by T.~Ohtsuki\footnote{
Department of Mathematical and Computing Sciences,
Tokyo Institute of Technology,
Oh-okayama, Meguro-ku, Tokyo 152-8552, Japan}
}
\end{center}
\def\text#1{ \mbox{#1} }
\vskip 1pc

The aim of this note is to show
Corollary \ref{cor.diagonalize} below,
which implies that
any rational homology 3-sphere can be obtained from $S^3$
by integral surgery along some algebraically split framed link
after adding some lens spaces.
For similar lemmas, see \cite{Oh2,murakami}.

Before showing the corollary we prepare notations of linking pairings.
Let $G$ be a finite Abelian group.
A {\it linking pairing} on $G$ is
a non-singular symmetric bilinear map of $G \times G$ to $\Q/\Z$.
For a non-singular symmetric integral $n \times n$ matrix $A$,
we have an induced linking pairing $\phi$ on $\Z^n/A\Z^n$
defined by $\phi([v],[v'])={}^t v A^{-1} v'$ for $v,v' \in \Z^n$
whose images in $\Z^n/A\Z^n$ are denoted by $[v],[v']$;
note that the right-hand side of this formula is well-defined in $\Q/\Z$.
We denote this linking pairing by $\iota(A)$.
It is known \cite{KP,Dur} that
if two non-singular symmetric integral matrices $A_1$, $A_2$
give the same linking pairing $\iota(A_1)$ and $\iota(A_2)$
then there exists a unimodular integral matrix $P$ such that
$$
{}^t P \cdot \Big( A_1 \oplus (\pm1) \oplus \cdots \oplus (\pm1) \Big) \cdot P
= A_2 \oplus (\pm1) \oplus \cdots \oplus (\pm1).
$$
The set of linking pairing becomes an Abelian semigroup
with respect to direct sum.
Generators and relations of the semigroup are known \cite{Wall,KaKo}.
The generators in \cite{Wall} are:
\begin{eqnarray*}
& [1/p^k], [d_p/p^k] \text{ on }\Z/p^k\Z \\
& \qquad\qquad
    \text{ for $p$ odd primes, $d_p$ a quadratic non-residue modulo $p$} \\
& [1/2] \text{ on } \Z/2\Z, \quad
  [1/2^2], [-1/2^2] \text{ on } \Z/2^2\Z \\
& [1/2^k], [-1/2^k], [3/2^k], [-3/2^k]
   \text{ on } \Z/2^k\Z \text{ for }k \ge 3 \\
& E_0^k \text{ on } \Z/2^k\Z \oplus \Z/2^k\Z \text{ for }k \ge 1 \\
& E_1^k \text{ on } \Z/2^k\Z \oplus \Z/2^k\Z \text{ for }k \ge 2
\end{eqnarray*}
where we denote by $[b/a]$ a linking pairing $\phi$ on $\Z/a\Z$
defined by $\phi([v],[v'])=bvv'/a$ for $v,v' \in \Z$
and we define $E_0^k$ and $E_1^k$ by
$$
E_0^k([v],[v'])=
{}^t v
\pmatrix{ 0 & 2^{-k} \cr 2^{-k} & 0 }
v', \qquad
E_1^k([v],[v'])=
{}^t v
\pmatrix{ 2^{1-k} & 2^{-k} \cr 2^{-k} &2^{1-k} }
v'
$$
for $v, v' \in \Z \oplus \Z$.

\begin{lem}\label{lem.diagonalize}
Let $\phi$ be any linking pairing
on a finite Abelian group $G$.
Then there exists linking pairings
$[b_1/a_1], [b_2/a_2], \cdots, [b_N/a_N]$
and integers $n_1, n_2, \cdots, n_\nu$
such that
$$
\phi \oplus \bigoplus_{i=1}^N [b_i/a_i]
= \iota({\bigoplus_{\xi=1}^{\nu} (n_\xi)})
$$
and the order of $G$ is divisible by each $a_i$.
\end{lem}

\proof
Since any linking pairing $\phi$
is equal to a direct sum of some generators given above,
it is sufficient to show the lemma when $\phi$ is equal to each generator.

If $\phi = [1/p^k]$, then we have $\phi = \iota \big((p^k)\big)$.

If $\phi = [d_p/p^k]$, then we have
$\phi \oplus [d_p/p^k] = 2[1/p^k] = \iota \big(2(p^k)\big)$
using the relation $2[d_p/p^k]=2[1/p^k]$ in \cite{Wall}.

If $\phi = [\pm1/2^k]$, then we have $\phi = \iota \big((\pm2^k) \big)$.

If $\phi = [\pm3/2^k]$, then we have
$\phi \oplus [\pm3/2^k] = 2[\mp1/2^k] = \iota \big(2(\mp 2^k)\big)$
using a relation $2[\pm3/2^k]=2[\mp1/2^k]$ in \cite{KaKo}.

If $\phi=E^k_0$ or $E^k_1$, then we have
\begin{eqnarray*}
E^k_0 \oplus [3/2^k]
&= 3[1/2^k] = \iota \big(3(2^k)\big) \\
E^k_1 \oplus [-1/2^k]
&= [1/2^k] \oplus 2[-1/2^k] = \iota \big((2^k) \oplus 2(-2^k)\big)
\end{eqnarray*}
using relations in \cite{KaKo}, completing the proof.

Let $L$ be a framed link in $S^3$,
$A$ the linking matrix of $L$
and $M$ the closed 3-manifold obtained by integral surgery along $L$.
Here we assume that
$M$ is a rational homology 3-sphere,
which implies that
the matrix $A$ is non-singular.
It is easily checked that
$H_1(M;\Z)$ is isomorphic to $\Z^n/A\Z^n$,
the linking pairing on $H_1(M;\Z)$ is equal to $\iota(A)$ and
the linking pairing on $H_1(L(a,b);\Z)$ is equal to $[b/a]$,
where $L(a,b)$ is the lens space of type $(a,b)$.
Further, for a uni-modular integral matrix $P$,
a framed link whose linking matrix is equal to
$$
{}^t P \cdot \Big( A \oplus (\pm1) \oplus \cdots \oplus (\pm1) \Big) \cdot P
$$
can be obtained from $L$
by applying Kirby moves to $L$;
in particular
the 3-manifold obtained by integral surgery along the new framed link
is homeomorphic to $M$.
By Lemma \ref{lem.diagonalize}
putting $\phi$ to be the linking pairing on $H_1(M;\Z)$, we have

\begin{cor}\label{cor.diagonalize}
For any rational homology 3-sphere $M$,
there exist lens spaces of types $(a_1,b_1), (a_2,b_2), \cdots, (a_N,b_N)$
such that
the connected sum of $M$ and these lens spaces
can be obtained by integral surgery along some algebraically split framed link
and the order of $H_1(M;\Z)$ is divisible by each $a_i$.
\end{cor}

\begin{rem}
As for lens spaces in Corollary \ref{cor.diagonalize}, in fact,
we need only lens spaces of types
$(n,1)$, $({p}^{k},d_{p})$ or $(2^{k}, \pm3)$;
we can see it
by checking the proof of Lemma \ref{lem.diagonalize}.
\end{rem}

\end{document}

\bigskip
\noindent
{\bf Acknowledgements.}
This work was supported by NSF Grant DMS-9704893.

\nappendixsp{1}{1}

In this paper we use the following notations. $\cL$ denotes mostly an
$L$-component link either in $S^3$ or in a rational homology sphere
$M$. Its components are denoted as $\cL_j$, while $\cLk$ is a sublink
of $\cL$ which contains the components $\cL_{k+1},\ldots,\cL_L$. The
linking numbers of $\cL$ in $M$ (or in $S^3$ if $M$ is not considered
at all) are denoted as $\lij$. We assume that the manifold $M$ is
constructed by a surgery on an $L\p$-component link $\cLp$. The
matrices $\laa,\lbb,\lab$ denote the linking numbers between the
components of $\cL$ and $\cL$, $\cLp$ and $\cLp$, $\cL$ and $\cLp$
respectively.

We also use the multi-index notations which vary depending on the
context. Usually when we work with $\cL$, $\tcL$ and
$\cLp$, then $\ux=\brlst{x}{1}{L}$, $\tux=\brlst{\tx}{1}{\tL}$
 and $\uux=\brlst{x}{1}{L\p}$,
however when we consider $\cL$ and $\cLk$, then
$\ux=\brlst{x}{k+1}{L}$ and $\uux=\brlst{x}{1}{k}$. We also use the
following shortcuts: if, say, $\ux=\brlst{x}{1}{L}$ then
\qq
&y\ux  =  \brlst{yx}{1}{L},\qquad
\ux^y  =  \brlst{x^y}{1}{L},\qquad
\ux\uy  =  \{x_1 y_1,\ldots,x_L y_L\},
\nonumber\\
&\prb{f(\ux)}  =  \pjoL f(x_j),\qquad
|\ux|  =  \sjoL x_j,
\nonumber
\qqq
%
%
and $\ux=y$ means that $x_j=y$ for all $1\leq j \leq L$. We also use
the notation
\qq
\ux^{\uy} = \pjoL x_j^{y_j}
\label{a1.6}
\qqq
except for
\qq
\ut^{\uy} = \{t_1^{y_1},\ldots,t_L^{y_L}\}.
\label{a1.7}
\qqq
%


\nappendixsp{2}{2}
\label{a2}

\begin{center}
{\Large \bf Construction of a rational homology sphere by a surgery
on an algebraically split link}

{\sc by T.~Ohtsuki\footnote{
Department of Mathematical and Computing Sciences,
Tokyo Institute of Technology,
Oh-okayama, Meguro-ku, Tokyo 152-8552, Japan}
}
\end{center}
\def\text#1{ \mbox{#1} }
\vskip 1pc

The aim of this note is to show
Corollary \ref{cor.diagonalize} below,
which implies that
any rational homology 3-sphere can be obtained from $S^3$
by integral surgery along some algebraically split framed link
after adding some lens spaces.
For similar lemmas, see \cite{Oh2,murakami}.

Before showing the corollary we prepare notations of linking pairings.
Let $G$ be a finite Abelian group.
A {\it linking pairing} on $G$ is
a non-singular symmetric bilinear map of $G \times G$ to $\Q/\Z$.
For a non-singular symmetric integral $n \times n$ matrix $A$,
we have an induced linking pairing $\phi$ on $\Z^n/A\Z^n$
defined by $\phi([v],[v'])={}^t v A^{-1} v'$ for $v,v' \in \Z^n$
whose images in $\Z^n/A\Z^n$ are denoted by $[v],[v']$;
note that the right-hand side of this formula is well-defined in $\Q/\Z$.
We denote this linking pairing by $\iota(A)$.
It is known \cite{KP,Dur} that
if two non-singular symmetric integral matrices $A_1$, $A_2$
give the same linking pairing $\iota(A_1)$ and $\iota(A_2)$
then there exists a unimodular integral matrix $P$ such that
$$
{}^t P \cdot \Big( A_1 \oplus (\pm1) \oplus \cdots \oplus (\pm1) \Big) \cdot P
= A_2 \oplus (\pm1) \oplus \cdots \oplus (\pm1).
$$
The set of linking pairing becomes an Abelian semigroup
with respect to direct sum.
Generators and relations of the semigroup are known \cite{Wall,KaKo}.
The generators in \cite{Wall} are:
\begin{eqnarray*}
& [1/p^k], [d_p/p^k] \text{ on }\Z/p^k\Z \\
& \qquad\qquad
    \text{ for $p$ odd primes, $d_p$ a quadratic non-residue modulo $p$} \\
& [1/2] \text{ on } \Z/2\Z, \quad
  [1/2^2], [-1/2^2] \text{ on } \Z/2^2\Z \\
& [1/2^k], [-1/2^k], [3/2^k], [-3/2^k]
   \text{ on } \Z/2^k\Z \text{ for }k \ge 3 \\
& E_0^k \text{ on } \Z/2^k\Z \oplus \Z/2^k\Z \text{ for }k \ge 1 \\
& E_1^k \text{ on } \Z/2^k\Z \oplus \Z/2^k\Z \text{ for }k \ge 2
\end{eqnarray*}
where we denote by $[b/a]$ a linking pairing $\phi$ on $\Z/a\Z$
defined by $\phi([v],[v'])=bvv'/a$ for $v,v' \in \Z$
and we define $E_0^k$ and $E_1^k$ by
$$
E_0^k([v],[v'])=
{}^t v
\pmatrix{ 0 & 2^{-k} \cr 2^{-k} & 0 }
v', \qquad
E_1^k([v],[v'])=
{}^t v
\pmatrix{ 2^{1-k} & 2^{-k} \cr 2^{-k} &2^{1-k} }
v'
$$
for $v, v' \in \Z \oplus \Z$.

\begin{lem}\label{lem.diagonalize}
Let $\phi$ be any linking pairing
on a finite Abelian group $G$.
Then there exists linking pairings
$[b_1/a_1], [b_2/a_2], \cdots, [b_N/a_N]$
and integers $n_1, n_2, \cdots, n_\nu$
such that
$$
\phi \oplus \bigoplus_{i=1}^N [b_i/a_i]
= \iota({\bigoplus_{\xi=1}^{\nu} (n_\xi)})
$$
and the order of $G$ is divisible by each $a_i$.
\end{lem}

\proof
Since any linking pairing $\phi$
is equal to a direct sum of some generators given above,
it is sufficient to show the lemma when $\phi$ is equal to each generator.

If $\phi = [1/p^k]$, then we have $\phi = \iota \big((p^k)\big)$.

If $\phi = [d_p/p^k]$, then we have
$\phi \oplus [d_p/p^k] = 2[1/p^k] = \iota \big(2(p^k)\big)$
using the relation $2[d_p/p^k]=2[1/p^k]$ in \cite{Wall}.

If $\phi = [\pm1/2^k]$, then we have $\phi = \iota \big((\pm2^k) \big)$.

If $\phi = [\pm3/2^k]$, then we have
$\phi \oplus [\pm3/2^k] = 2[\mp1/2^k] = \iota \big(2(\mp 2^k)\big)$
using a relation $2[\pm3/2^k]=2[\mp1/2^k]$ in \cite{KaKo}.

If $\phi=E^k_0$ or $E^k_1$, then we have
\begin{eqnarray*}
E^k_0 \oplus [3/2^k]
&= 3[1/2^k] = \iota \big(3(2^k)\big) \\
E^k_1 \oplus [-1/2^k]
&= [1/2^k] \oplus 2[-1/2^k] = \iota \big((2^k) \oplus 2(-2^k)\big)
\end{eqnarray*}
using relations in \cite{KaKo}, completing the proof.

Let $L$ be a framed link in $S^3$,
$A$ the linking matrix of $L$
and $M$ the closed 3-manifold obtained by integral surgery along $L$.
Here we assume that
$M$ is a rational homology 3-sphere,
which implies that
the matrix $A$ is non-singular.
It is easily checked that
$H_1(M;\Z)$ is isomorphic to $\Z^n/A\Z^n$,
the linking pairing on $H_1(M;\Z)$ is equal to $\iota(A)$ and
the linking pairing on $H_1(L(a,b);\Z)$ is equal to $[b/a]$,
where $L(a,b)$ is the lens space of type $(a,b)$.
Further, for a uni-modular integral matrix $P$,
a framed link whose linking matrix is equal to
$$
{}^t P \cdot \Big( A \oplus (\pm1) \oplus \cdots \oplus (\pm1) \Big) \cdot P
$$
can be obtained from $L$
by applying Kirby moves to $L$;
in particular
the 3-manifold obtained by integral surgery along the new framed link
is homeomorphic to $M$.
By Lemma \ref{lem.diagonalize}
putting $\phi$ to be the linking pairing on $H_1(M;\Z)$, we have

\begin{cor}\label{cor.diagonalize}
For any rational homology 3-sphere $M$,
there exist lens spaces of types $(a_1,b_1), (a_2,b_2), \cdots, (a_N,b_N)$
such that
the connected sum of $M$ and these lens spaces
can be obtained by integral surgery along some algebraically split framed link
and the order of $H_1(M;\Z)$ is divisible by each $a_i$.
\end{cor}

\begin{rem}
As for lens spaces in Corollary \ref{cor.diagonalize}, in fact,
we need only lens spaces of types
$(n,1)$, $({p}^{k},d_{p})$ or $(2^{k}, \pm3)$;
we can see it
by checking the proof of Lemma \ref{lem.diagonalize}.
\end{rem}

\end{document}

\bibitem{ohtsuki}
Ohtsuki, T.,
{\it A polynomial invariant of rational homology 3-spheres},
Invent. Math. {\bf 123} (1996) 241--257.

\pub
"The Trivial Connection Contribution to Witten's Invariant and Finite
Type Invariants of Rational Homology Spheres"
\journal Communications in Mathematical Physics & 183, (97) 23-54.

\pub
"Witten's Invariant of Rational Homology Spheres at Prime Values of
$K$ and Trivial Connection Contribution"
\journal Communications in Mathematical Physics & 180 (96) 297-324.


\pub
"On Finite Type Invariants of Links and Rational Homology Spheres
Derived from the Jones Polynomial and Witten-Reshetikhin-Turaev
Invariant" \rm preprint q-alg/9511025.

\pub
"Higher Order Terms in the Melvin-Morton Expansion of the Colored
Jones Polynomial"
\journal Communications in Mathematical Physics & 183 (97) 291-306.

\pub
"On $p$-Adic Convergence of Perturbative Invariants of Some Rational
Homology Spheres" \rm preprint q-alg/9601015, to be published in Duke
Math. Journ.

\pub
"The Universal $R$-Matrix, Burau Representation and the Melvin-Morton
Expansion of the Colored Jones Polynomial" \rm preprint
q-alg/9604005, to be published in Adv. Math.

\end{document}